\documentclass{article} 
\usepackage{nips15submit_e,times}

\usepackage{hyperref}
\usepackage{url}

\usepackage{graphicx} 
\usepackage{epsfig}
\usepackage{subfigure}

\usepackage{amssymb}
\usepackage{amsmath}
\usepackage{amsthm}
\usepackage{mathtools}
\usepackage[raggedright]{titlesec} 

\usepackage{algorithm,algcompatible,amsmath}
\algnewcommand\INPUT{\item[\textbf{Input:}]}%
\algnewcommand\OUTPUT{\item[\textbf{Output:}]}%

\usepackage{arydshln}
\def\mat#1{\mbox{\bf #1}}

\DeclareMathAlphabet\mathbfcal{OMS}{cmsy}{b}{d}

\newcommand{\changeBM}[1]{#1}
\newcommand{\changeHK}[1]{#1}
\newcommand{\changeBMM}[1]{#1}
\newcommand{\changeHKK}[1]{#1}
\newcommand{\changeBMMM}[1]{#1}

\newcommand{\argmin}{\operatornamewithlimits{arg\,min}}
\makeatletter
\def\hlinewd#1{%
  \noalign{\ifnum0=`}\fi\hrule \@height #1 \futurelet
   \reserved@a\@xhline}
\makeatother

\title{Riemannian preconditioning for tensor completion}

\date{Compiled on \today}

\author{
Hiroyuki Kasai\thanks{Technische Universit\"{a}t M\"{u}nchen, Department of Electrical and Computer Engineering, Munich, Germany  (\texttt{hiroyuki.kasai@tum.de}).}\\
Graduate School of Information Systems,\\
The university of Electro-Communications\\
Chofu-shi, Tokyo, 182-8585, Japan\\
\texttt{kasai@is.uec.ac.jp} \\
\And
Bamdev Mishra\thanks{University of Cambridge, Department of Engineering, Cambridge, UK (\texttt{bm458@cam.ac.uk}).}\\
Department of EECS,\\
University of Li\`ege\\
4000 Li\`ege, Belgium \\
\texttt{b.mishra@ulg.ac.be} \\
}

\nipsfinalcopy 

\begin{document}


\maketitle

\begin{abstract}
We propose a novel Riemannian preconditioning approach for the tensor completion problem with rank constraint. A Riemannian metric or inner product is proposed that exploits the least-squares structure of the cost function and takes into account the structured \changeBMMM{symmetry} in Tucker decomposition. The specific metric allows to use the versatile framework of Riemannian optimization on quotient manifolds to develop a \changeBMMM{preconditioned} nonlinear conjugate gradient algorithm for the problem. To this end, concrete matrix representations of various optimization-related ingredients are listed. Numerical comparisons suggest that our proposed algorithm \changeBM{robustly} outperforms state-of-the-art algorithms across different problem instances encompassing various synthetic and real-world datasets\footnote{\textcolor{blue}{An extended version of the paper is in \cite{kasai16a}, which includes a stochastic gradient descent algorithm for low-rank tensor completion. The extended version has been accepted to the 33rd International Conference on Machine Learning (ICML 2016).}}.
\end{abstract}

\section{Introduction}
\label{sec:Introduction}
This paper addresses the problem of low-rank tensor completion when the rank is a priori known or estimated. Without loss of generality, we focus on 3-order tensors. Given a tensor $\mathbfcal{X}^{n_1 \times n_2 \times n_3}$, whose entries $\mathbfcal{X}_{i_1, i_2, i_3}^{\star}$ are only known for some indices $(i_1, i_2, i_3) \in \Omega$, where $\Omega$ is a subset of the complete set of indices 
\changeBMMM{$\{(i_1, i_2, i_3): i_d \in \{1, \ldots, n_d \}, d \in \{1,2,3\}\}$},
the \emph{fixed-rank tensor completion problem} is formulated as 
\begin{equation}
\begin{array}{lll}
\label{Eq:CostFunction}
\displaystyle{\min_{\mathbfcal{X} \in
\mathbb{R}^{n_1 \times n_2 \times n_3}} }&   
\displaystyle{\frac{1}{|\Omega |}
\| \mathbfcal{P}_{\Omega}(\mathbfcal{X}) - 
\mathbfcal{P}_{\Omega}(\mathbfcal{X}^{\star}) \|^2_F} \\
{\rm subject\ to}& {\rm rank}(\mathbfcal{X}) = {\bf r},
\end{array}
\end{equation}
where the operator $\mathbfcal{P}_{\Omega}(\mathbfcal{X})_{i_1 i_2 i_3} = \mathbfcal{X}_{i_1 i_2 i_3}$ if $(i_1, i_2, i_3) \in \Omega$ and $\mathbfcal{P}_{\Omega}(\mathbfcal{X})_{i_1 i_2 i_3}  = 0$ otherwise and (with a slight abuse of notation) $\|\cdot \|_F$ is the Frobenius norm. ${\rm rank}(\mathbfcal{X})$ $ (={\bf r}=(r_1, r_2, r_3))$, called the \emph{multilinear rank} of $\mathbfcal{X}$, is the set of the ranks of for each of mode-$d$ unfolding matrices. \changeBM{$r_d \ll n_d$ enforces a low-rank structure.} The {\it mode} is a matrix obtained by concatenating the mode-$d$ fibers along columns, and mode-$d$ {\it unfolding} of $\mathbfcal{X}$ is $\mat{X}_{d} \in \mathbb{R}^{n_d \times n_{d+1}\cdots n_D n_1 \cdots n_{d-1}}$ for $d=\{1,\ldots,D\}$.

Problem (\ref{Eq:CostFunction}) \changeBM{has} many variants, and one of those is extending the nuclear norm regularization approach from the matrix case \cite{Candes_FoundCompuMath_2009_s} to the tensor case. This results in a summation of nuclear norm regularization terms, each one corresponds to each of the unfolding matrices of $\mathbfcal{X}$. While this generalization leads to good results \cite{Liu_IEEETransPAMI_2013_s, Tomioka_Latent_2011_s, Signoretto_MachineLearning_2014_s}, \changeBMMM{its} applicability to large-scale instances is not trivial, especially due to the necessity of high-dimensional singular value decomposition computations. A different approach exploits \emph{Tucker decomposition} \cite[Section~4]{Kolda_SIAMReview_2009_s} of a low-rank tensor $\mathbfcal{X}$ to develop large-scale algorithms for (\ref{Eq:CostFunction}), e.g., in \cite{Filipovi_MultiSysSigPro_2013_s,Kressner_BIT_2014_s}.

The present paper exploits both the \emph{symmetry} present in Tucker decomposition and the \emph{least-squares} structure of the cost function of (\ref{Eq:CostFunction}) to develop a competitive algorithm. To this end, we use the concept of \emph{preconditioning}. While preconditioning in unconstrained optimization is well studied \cite[Chapter~5]{Nocedal_NumericalOpt_2006_s}, preconditioning on constraints with \emph{symmetries}, owing to non-uniqueness of Tucker decomposition \cite[Section~4.3]{Kolda_SIAMReview_2009_s}, is not straightforward. We build upon the recent work \cite{Bamdev_arXiv_2014_s} that suggests to use \emph{Riemannian preconditioning} with a \emph{tailored metric} (inner product) in the Riemannian optimization framework on quotient manifolds \cite{Absil_OptAlgMatManifold_2008, Smith94a, Edelman98a}. Use of Riemannian preconditioning for the low-rank matrix completion problem is discussed in \cite{Mishra_ICDC_2014_s}, where a \emph{preconditioned nonlinear conjugate gradient} algorithm is proposed. It connects to state-of-the-art algorithms in \cite{Ngo_NIPS_2012_s, Wen_MPC_2012_s} and shows competitive performance. In this paper, \changeBMM{we generalize the work \cite{Mishra_ICDC_2014_s} to tensor completion.}

The paper is organized as follows. Section \ref{sec:Fixed-rankTuckerFactorization} discusses the two fundamental structures of symmetry and least-squares associated with (\ref{Eq:CostFunction}) and proposes a novel metric that captures the relevant second-order information of the problem. The optimization-related ingredients on the Tucker manifold are developed in Section \ref{sec:OptimizationRelatedIngredients}. The cost function specific ingredients are developed in Section \ref{sec:AlgorithmDetails}. The final formulas are listed in Table \ref{tab:FinalFormulas}. \changeBM{In Section \ref{sec:NumericalComparisons},} numerical comparisons with state-of-the-art algorithms on various synthetic (both small and large-scale instances) and real-world benchmarks suggest a \changeBM{superior} performance \changeHK{of} our proposed algorithm. Our proposed preconditioned nonlinear conjugate gradient algorithm is implemented\footnote{The Matlab code is available at \url{http://bamdevmishra.com/codes/tensorcompletion/}.} in the Matlab toolbox Manopt \cite{Boumal_Manopt_2014_s}.

The concrete developments of optimization-related ingredients and additional numerical experiments are shown in Sections \ref{Sup_sec:Derivationofmanifold-relatedingredients} and \ref{sec:AdditionalNumericalComparisons}, respectively, of the supplementary material section. 



\section{Exploiting the problem structure}
\label{sec:Fixed-rankTuckerFactorization}
Construction of efficient algorithms depends on properly exploiting both the structure of constraints and cost function. \changeBMM{To this end, we focus on two fundamental structures in (\ref{Eq:CostFunction}): \emph{symmetry} in the constraints, and the \changeHKK{\emph{least-squares structure}} of the cost function}. Finally, a novel metric is proposed.

{\bf The quotient structure of Tucker decomposition.} The Tucker decomposition of a tensor $\mathbfcal{X} \in \mathbb{R}^{n_1 \times n_2 \times n_3}$ of rank {\bf r} (=$(r_1, r_2, r_3)$) is \cite[Section~4.1]{Kolda_SIAMReview_2009_s}
\begin{equation}
\begin{array}{lll}
\label{Eq:TuckerFactorization}
\mathbfcal{X} &  = & 
\mathbfcal{G} {\times_1} \mat{U}_{1} {\times_2} \mat{U}_{2} {\times_3} \mat{U}_{3}, 
\end{array}
\end{equation}
where $\mat{U}_{d} \in {\rm St}(r_d, n_d)$ \changeHK{for $d \in \{1,2,3\}$} belongs to \changeBMM{the \emph{Stiefel manifold} of} matrices of size $n_d \times r_d$ with orthogonal columns and $\mathbfcal{G} \in \mathbb{R}^{r_1 \times r_2 \times r_3}$. Here, $\mathbfcal{W} \changeHK{\times_d} \mat{V} \in \mathbb{R}^{\changeHK{n_1 \times \cdots n_{d-1} \times m \times n_{d+1}  \times \cdots n_N}}$ \changeBMM{computes} the {\it \changeHK{d}-mode product} of a tensor $\mathbfcal{W} \changeHK{\in \mathbb{R}^{n_1 \times \cdots \times n_N}}$ and a matrix $\mat{V} \in \changeHK{\mathbb{R}}^{\changeHK{m} \times \changeHK{n_d}}$.

Tucker decomposition (\ref{Eq:TuckerFactorization}) is \emph{not
unique} as $\mathbfcal{X}$ remains unchanged under the transformation \begin{equation}
\label{Eq:EquivalenceClass_3}
(\mat{U}_{1}, \mat{U}_{2}, \mat{U}_{3}, \mathbfcal{G}) \mapsto  (  \mat{U}_{1}\mat{O}_{1}, \mat{U}_{2}\mat{O}_{2}, \mat{U}_{3}\mat{O}_{3}, \mathbfcal{G} {\times_1} \mat{O}^T_{1} {\times_2} \mat{O}^T_{2} {\times_3} \mat{O}^T_{3}) 
\end{equation}
for all $\mat{O}_{d} \in \mathcal{O}(r_d)$, the set of orthogonal matrices of size of $r_d \times r_d$. The classical remedy to remove this indeterminacy is to have additional structures on $\mathbfcal{G}$ like sparsity or restricted orthogonal rotations \cite[Section~4.3]{Kolda_SIAMReview_2009_s}. In contrast, we encode the transformation (\ref{Eq:EquivalenceClass_3}) in an abstract search space of \emph{equivalence classes}, defined as,
\begin{equation}
\label{Eq:EquivalenceClass}
[(\mat{U}_{1}, \mat{U}_{2}, \mat{U}_{3}, \mathbfcal{G}) ] : =  \{(  \mat{U}_{1}\mat{O}_{1}, \mat{U}_{2}\mat{O}_{2}, \mat{U}_{3}\mat{O}_{3}, \mathbfcal{G} {\times_1} \mat{O}^T_{1} {\times_2} \mat{O}^T_{2} {\times_3} \mat{O}^T_{3}) : \mat{O}_{d} \in \mathcal{O}(r_d)\}.
\end{equation}

The set of equivalence classes is the quotient manifold \cite[Theorem~9.16]{Lee03a}
\begin{equation}
\begin{array}{lll}
\label{Eq:QuotientSpace}
\mathcal{M}/\!\sim
& := & \mathcal{M}/ 
(\mathcal{O}{(r_1)} \times \mathcal{O}{(r_2)} \times \mathcal{O}{(r_3)}),
\end{array}
\end{equation}
where $\mathcal{M}$ is called the \emph{total space} (computational space) that is the product space
\begin{equation}
\begin{array}{lll}
\label{Eq:TotalSpace}
\mathcal{M}
& := & {\rm St}(r_1, n_1) \times {\rm St}(r_2, n_2) \times {\rm St}(r_3, n_3) \times \mathbb{R}^{r_1 \times r_2 \times r_3}.
\end{array}
\end{equation}

Due to \changeBM{the} invariance (\ref{Eq:EquivalenceClass_3}), the local minima of  (\ref{Eq:CostFunction}) in $\mathcal{M}$ are not isolated, but they become isolated on $\mathcal{M}/ \!\sim$. Consequently, the problem (\ref{Eq:CostFunction}) is an optimization problem on a quotient manifold for which systematic procedures are proposed in \cite{Absil_OptAlgMatManifold_2008, Smith94a, Edelman98a} \changeBM{by endowing $\mathcal{M}/ \!\sim$ with a Riemannian structure}. We call $\mathcal{M}/ \!\sim$, defined in (\ref{Eq:QuotientSpace}), the \emph{Tucker manifold} as it results \changeBMM{from Tucker decomposition}.

{\bf The least-squares structure of the cost function.}
In unconstrained optimization, the Newton method is interpreted as a \emph{scaled} steepest descent method, where the search space is endowed with a metric (inner product) induced by the Hessian of the cost function \cite{Nocedal_NumericalOpt_2006_s}. This induced metric (or its approximation) resolves convergence issues of first-order optimization algorithms. Analogously, finding a good inner product for (\ref{Eq:CostFunction}) is of profound consequence. Specifically for the case of quadratic optimization with rank constraint (matrix case), Mishra and Sepulchre \cite[Section~5]{Bamdev_arXiv_2014_s} propose a family of Riemannian metrics from the Hessian of the cost function. Applying this approach directly for the particular cost function of (\ref{Eq:CostFunction}) is computationally costly. To circumvent the issue, we consider a simplified cost function by assuming that $\Omega$ contains the full set of indices, i.e., we focus on $\| \mathbfcal{X} - \mathbfcal{X}^{\star}\|_F^2$ to propose a metric candidate. Applying the metric tuning approach of \cite[Section~5]{Bamdev_arXiv_2014_s} to the simplified cost function leads to a family of Riemannian metrics. A good trade-off between computational cost and simplicity is by considering only the \emph{block diagonal} elements of the Hessian of $\| \mathbfcal{X} - \mathbfcal{X}^{\star}\|_F^2$. It should be noted that the cost function $\| \mathbfcal{X} - \mathbfcal{X}^{\star}\|_F^2$  is {\it convex and quadratic} in $\mathbfcal{X}$. Consequently, it is also convex and quadratic in the arguments 
$(\mat{U}_{1}, \mat{U}_{2}, \mat{U}_{3}, \mathbfcal{G})$ individually. Equivalently, the block diagonal approximation of the Hessian of $\| \mathbfcal{X} - \mathbfcal{X}^{\star}\|_F^2$ in $(\mat{U}_{1}, \mat{U}_{2}, \mat{U}_{3}, \mathbfcal{G})$ is  
\begin{equation}
\begin{array}{lll}
\label{Eq:BlockApproximation}
((\mat{G}_{1} \mat{G}_{1}^T) \otimes \mat{I}_{n_1}, (\mat{G}_{2} \mat{G}_{2}^T) \otimes \mat{I}_{n_2}, (\mat{G}_{3} \mat{G}_{3}^T) \otimes \mat{I}_{n_3}, 
\mat{I}_{r_1 r_2  r_3}),
\end{array}
\end{equation}
where $\mat{G}_{d}$ is the mode-$d$ unfolding of $\mathbfcal{G}$ and is assumed to be full rank. The terms $\mat{G}_{d} \mat{G}_{d}^T$ for $d  \in \{1, 2,3\}$ are \emph{positive definite} when $r_1 \leq r_2 r_3$, $r_2 \leq r_1 r_3$, and $r_3 \leq r_1 r_2$, which is a reasonable modeling assumption.

{\bf A novel Riemannian metric.} An element $x$ in the total space $\mathcal{M}$ has the matrix representation $(\mat{U}_{1}, \mat{U}_{2}, \mat{U}_{3}, \mathbfcal{G})$. Consequently, the tangent space $T_{{x}} \mathcal{M}$ is the Cartesian product of the tangent spaces of the individual manifolds of (\ref{Eq:TotalSpace}), i.e., $T_{{x}} \mathcal{M}$ has the matrix characterization \cite{Edelman98a}
\begin{equation}
\begin{array}{lll}
\label{Eq:tangent_space}
T_{{x}} {\mathcal{M}} 
& = &  \{ (\mat{Z}_{{\bf U}_{1}}, \mat{Z}_{{\bf U}_{2}}, \mat{Z}_{{\bf U}_{3}}, \mat{Z}_{\mathbfcal{G}}) \in 
\mathbb{R}^{n_1 \times r_1} \times  
\mathbb{R}^{n_2 \times r_2} \times 
\mathbb{R}^{n_3 \times r_3} \times 
\mathbb{R}^{r_1 \times r_2 \times r_3} : \\
&  & \mat{U}_{d}^T \mat{Z}_{{\bf U}_{d}} +  \mat{Z}_{{\bf U}_{d}}^T \mat{U}_{d} = 0,\ {\rm for\ } d \changeHK{\in} \{1,2,3 \} \}.
\end{array}
\end{equation}

From the earlier discussion on symmetry and least-squares structure, we propose the novel metric ${g}_{x}:T_x \mathcal{M} \times T_x \mathcal{M} \rightarrow \mathbb{R}$
\begin{equation}
\begin{array}{lll}
\label{Eq:metric}
{g}_{x}(\xi_{x}, \eta_{x}) 
 =  
\langle \xi_{\scriptsize \mat{U}_{1}},
{\eta}_{\scriptsize\mat{U}_{1}} (\mat{G}_{1} \mat{G}_{1}^T) \rangle +
\langle \xi_{\scriptsize \mat{U}_{2}},
{\eta}_{\scriptsize\mat{U}_{2}} (\mat{G}_{2} \mat{G}_{2}^T) \rangle +
\langle \xi_{\scriptsize \mat{U}_{3}},
{\eta}_{\scriptsize\mat{U}_{3}} (\mat{G}_{3} \mat{G}_{3}^T) \rangle 
+ \langle {\xi}_{\scriptsize \mathbfcal{G}}, {\eta}_{\scriptsize \mathbfcal{G}}\rangle ,
\end{array}
\end{equation}
where ${\xi}_{{x}}, {\eta}_{{x}} \in T_{{x}} {\mathcal{M}}$ are 
tangent vectors with matrix characterizations, shown in (\ref{Eq:tangent_space}), $({\xi}_{\scriptsize \mat{U}_{1}}, {\xi}_{\scriptsize \mat{U}_{2}}, {\xi}_{\scriptsize \mat{U}_{3}}, {\xi}_{\scriptsize \mathbfcal{G}})$ and $({\eta}_{\scriptsize \mat{U}_{1}}, {\eta}_{\scriptsize \mat{U}_{2}}, {\eta}_{\scriptsize \mat{U}_{3}}, {\eta}_{\scriptsize \mathbfcal{G}})$, respectively and $\langle \cdot, \cdot \rangle$ is the Euclidean inner product. It should be emphasized that the proposed metric (\ref{Eq:metric}) is induced from (\ref{Eq:BlockApproximation}).

%
%
%
\section{Notions of optimization on the Tucker manifold}
\label{sec:OptimizationRelatedIngredients}
%
%
%



\begin{figure}[h]
\center
\includegraphics[scale = 0.35]{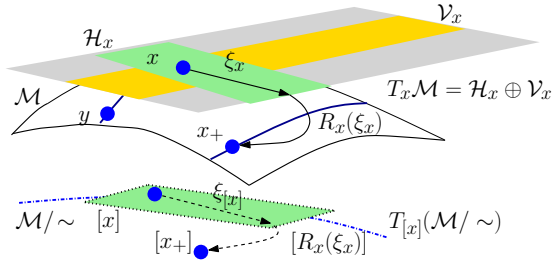}
\caption{Riemannian optimization framework: geometric objects, \changeBMM{shown in dotted lines}, on the quotient manifold $\mathcal{M}/\!\sim$ call for matrix representatives, \changeBMM{shown in solid lines}, in the total space $\mathcal{M}$.}
\label{fig:OptimizationOnSymmetries}
\end{figure}

Each point on a quotient manifold represents an entire equivalence class of matrices in the total space. Abstract geometric objects on a quotient manifold call for matrix representatives in the total space. Similarly, algorithms are run in the total space $\mathcal{M}$, but under appropriate compatibility between the Riemannian structure of $\mathcal{M}$ and the Riemannian structure of the quotient manifold $\mathcal{M}/ \!\sim$, they define algorithms on the quotient manifold. The key is endowing $\mathcal{M}/ \!\sim$ with a Riemannian structure. Once this is the case, a constraint optimization problem, for example (\ref{Eq:CostFunction}), is conceptually transformed into an unconstrained optimization over the Riemannian quotient manifold (\ref{Eq:QuotientSpace}). Below we briefly show the development of various geometric objects that are required to optimize a smooth cost function on the quotient manifold (\ref{Eq:QuotientSpace}) with first-order methods, e.g., conjugate gradients.





{\bf Quotient manifold representation and horizontal lifts.} Figure \ref{fig:OptimizationOnSymmetries} illustrates a schematic view of optimization with equivalence classes, where the points $x$ and $y$ in $\mathcal{M}$ belong to the same equivalence class (shown in solid blue color) and they represent a single point $[x]:=\{ y \in \mathcal{M} : y \sim x\}$ on the quotient manifold $\mathcal{M}/\!\sim$. The abstract tangent space $T_{[x]}(\mathcal{M}/\!\sim)$ at $[x] \in \mathcal{M}/\!\sim$ has the matrix representation in $T_x \mathcal{M}$, but restricted to the directions that do not induce a displacement along the equivalence class $[x]$. This is realized by decomposing $T_x \mathcal{M}$ into two complementary subspaces, the vertical and horizontal subspaces. The vertical space $\mathcal{V}_x$ is the tangent space of the equivalence class $[x]$. On the other hand, the horizontal space $\mathcal{H}_x$ is the \emph{orthogonal subspace} to $\mathcal{V}_x$ in the sense of the metric (\ref{Eq:metric}). Equivalently, $T_x \mathcal{M} =\mathcal{V}_x \oplus \mathcal{H}_x $. The horizontal subspace provides a valid matrix representation to the abstract tangent space $T_{[x]}(\mathcal{M}/\!\sim)$ \cite[Section~3.5.8]{Absil_OptAlgMatManifold_2008}. An abstract tangent vector $\xi_{[x]} \in T_{[x]}(\mathcal{M}/\!\sim)$ at $[x]$ has a unique element $\xi_x \in \mathcal{H}_x$ that is called its \emph{horizontal lift}.

A Riemannian metric $g_x:T_x \mathcal{M} \times  T_x \mathcal{M} \rightarrow \mathbb{R}$ at $x \in \mathcal{M}$ defines a Riemannian metric $g_{[x]}:T_{[x]}(\mathcal{M}/\! \sim) \times T_{[x]}(\mathcal{M}/\!\sim) \rightarrow \mathbb{R}$, i.e., $g_{[x]}(\xi_{[x]},\eta_{[x]}) := g_{x}(\xi_{x},\eta_{x})$ on the quotient manifold $\mathcal{M}/ \!\sim$, if $g_{x}(\xi_{x},\eta_{x})$ does not depend on a specific representation along the equivalence class $[x]$. Here, $\xi_{[x]}$ and $\eta_{[x]}$ are tangent vectors in $T_{[x]}(\mathcal{M}/\!\sim)$, and $\xi_{x}$ and $\eta_{x}$ are their horizontal lifts in $\mathcal{H}_x$ at $x$, respectively. Equivalently, the definition of the Riemannian metric is well posed when $g_x(\xi_x, \zeta_x)=g_x(\xi_y, \zeta_y)$ for all $x, y \in [x]$, where $\xi_x, \zeta_x \in \mathcal{H}_x$ and $\xi_y, \zeta_y \in \mathcal{H}_y$ are the horizontal lifts of $\xi_{[x]}, \zeta_{[x]} \in T_{[x]}(\mathcal{M}/\!\sim)$ along the same equivalence class $[x]$. \changeBM{From \cite[Proposition 3.6.1]{Absil_OptAlgMatManifold_2008}, it suffices to show that the metric (\ref{Eq:metric}) for tangent vectors $\xi_x, \zeta_x \in T_x \mathcal{M}$ does not change under the transformations 
$(\mat{U}_{1}, \mat{U}_{2}, \mat{U}_{3}, \mathbfcal{G}) \mapsto (\mat{U}_{1}\mat{O}_{1}, \mat{U}_{2}\mat{O}_{2}, \mat{U}_{3}\mat{O}_{3}, \mathbfcal{G} {\times_1} \mat{O}^T_{1} {\times_2} \mat{O}^T_{2} {\times_3} \mat{O}^T_{3})$, $(\xi_{\scriptsize \mat{U}_{1}}, \xi_{\scriptsize \mat{U}_{2}}, \xi_{\scriptsize \mat{U}_{3}}, \xi_{\scriptsize \mathbfcal{G}}) \mapsto
(\xi_{\scriptsize \mat{U}_{1}}\mat{O}_{1}, \xi_{\scriptsize \mat{U}_{2}}\mat{O}_{2}, \xi_{\scriptsize \mat{U}_{3}}\mat{O}_{3}, 
\xi_{\scriptsize \mathbfcal{G}} {\times_1} \mat{O}^T_{1} {\times_2} \mat{O}^T_{2} {\times_3} \mat{O}^T_{3})$, and $(\zeta_{\scriptsize \mat{U}_{1}}, \zeta_{\scriptsize \mat{U}_{2}}, \zeta_{\scriptsize \mat{U}_{3}}, \zeta_{\scriptsize \mathbfcal{G}}) \mapsto
(\zeta_{\scriptsize \mat{U}_{1}}\mat{O}_{1}, \zeta_{\scriptsize \mat{U}_{2}}\mat{O}_{2}, \zeta_{\scriptsize \mat{U}_{3}}\mat{O}_{3}, \zeta_{\scriptsize \mathbfcal{G}} {\times_1} \mat{O}^T_{1} {\times_2} \mat{O}^T_{2} {\times_3} \mat{O}^T_{3})$. A few straightforward computations show that this is indeed the case. Endowed with the Riemannian metric (\ref{Eq:metric}), the quotient manifold $\mathcal{M}/ \!\sim$ is a {\it Riemannian submersion} of $\mathcal{M}$. The submersion principle allows to work out concrete matrix representations of abstract object on $\mathcal{M}/ \!\sim$, e.g., the gradient of a smooth cost function \cite[Section~3.62]{Absil_OptAlgMatManifold_2008}.
}

Starting from an arbitrary matrix (with appropriate dimensions), two linear projections are needed: the first projection $\Psi_{x}$ is onto the tangent space $T_x\mathcal{M}$, while the
second projection $\Pi_{{x}}$ is onto the horizontal subspace $\mathcal{H}_x$. \changeBM{The computation cost of these projections is $O(n_1 r_1^2 + n_2 r_2^2 + n_3r_3  ^2)$.}

The tangent space $T_{x} \mathcal{M}$ projection operation is obtained by extracting the component normal to $T_{x} \mathcal{M}$ in the ambient space. The normal space $N_{x} \mathcal{M}$ has the matrix characterization $\{(\mat{U}_{1}\mat{S}_{\scriptsize \mat{U}_{1}}(\mat{G}_{1}  \mat{G}_{1}^T)^{-1},
\mat{U}_{2}\mat{S}_{\scriptsize \mat{U}_{2}}(\mat{G}_{2}  \mat{G}_{2}^T)^{-1},
\mat{U}_{3}\mat{S}_{\scriptsize \mat{U}_{3}}(\mat{G}_{3}  \mat{G}_{3}^T)^{-1}, 0) :\mat{S}_{\scriptsize \mat{U}_{d}}  \in \mathbb{R}^{r_d \times r_d},  \mat{S}_{\scriptsize \mat{U}_{d}}^T  = \mat{S}_{\scriptsize \mat{U}_{d}}, \text{ for\ } d \in \{1, 2, 3\} \}$. Symmetric matrices \changeHK{${\mat S}_{\scriptsize \mat{U}_d}$} for all $d \in \{1,2,3 \}$ parameterize the normal space. Finally, the operator $\Psi_{{x}}: 
\mathbb{R}^{n_1 \times r_1} \times  
\mathbb{R}^{n_2 \times r_2} \times 
\mathbb{R}^{n_3 \times r_3} \times 
\mathbb{R}^{r_1 \times r_2 \times r_3} \rightarrow T_{{x}} {\mathcal{M}} :(\mat{Y}_{\scriptsize \mat{U}_{1}}, \mat{Y}_{\scriptsize \mat{U}_{2}}, \mat{Y}_{\scriptsize \mat{U}_{3}}, \mat{Y}_{\scriptsize \mathbfcal{G}} )$
$\mapsto \Psi_{{x}}(\mat{Y}_{\scriptsize \mat{U}_{1}}, \mat{Y}_{\scriptsize \mat{U}_{2}}, \mat{Y}_{\scriptsize \mat{U}_{3}}, \mat{Y}_{\scriptsize \mathbfcal{G}} )$ has the matrix characterization
\begin{equation}
\label{Eq:B_Requirements}
\begin{array}{lll}
\Psi_{{x}}(\mat{Y}_{\scriptsize \mat{U}_{1}}\!, \mat{Y}_{\scriptsize \mat{U}_{2}}\!, \mat{Y}_{\scriptsize \mat{U}_{3}}\!, \mat{Y}_{\scriptsize \mathbfcal{G}} )
&= &   (\mat{Y}_{\scriptsize \mat{U}_{1}}\!\! - \!\mat{U}_{1} \mat{S}_{\scriptsize \mat{U}_{1}} (\mat{G}_{1}  \mat{G}_{1}^T)^{-1},
\mat{Y}_{\scriptsize \mat{U}_{2}}\!\! - \!\mat{U}_{2} \mat{S}_{\scriptsize \mat{U}_{2}} (\mat{G}_{2}  \mat{G}_{2}^T)^{-1}, \\
&  &\mat{Y}_{\scriptsize \mat{U}_{3}}\!\! -\! \mat{U}_{3}\mat{S}_{\scriptsize \mat{U}_{3}} (\mat{G}_{3}  \mat{G}_{3}^T)^{-1},
\mat{Y}_{\scriptsize \mathbfcal{G}}),
\end{array}
\end{equation}
where $\mat{S}_{\scriptsize \mat{U}_{d}}$ is the solution to the \emph{Lyapunov} equation $\mat{S}_{\scriptsize \mat{U}_{d}} \mat{G}_{d}  \mat{G}_{d}^T + \mat{G}_{d}  \mat{G}_{d}^T \mat{S}_{\scriptsize \mat{U}_{d}}   
 = \mat{G}_{d}  \mat{G}_{d}^T (\mat{Y}_{\scriptsize \mat{U}_{d}}^T \mat{U}_{d} + \mat{U}_{d}^T \mat{Y}_{\scriptsize \mat{U}_{d}}) \mat{G}_{d}  \mat{G}_{d}^T$ for $ d\in \{ 1,2,3\}$,
\changeBM{which are solved efficiently with the Matlab's} \verb+lyap+ \changeBM{routine.} 

The horizontal space projection operator of a tangent vector is obtained by removing the component along the vertical space. In particular, the vertical space $\mathcal{V}_{x} $ has the matrix characterization
$\{ (\mat{U}_1{\bf \Omega}_1, \mat{U}_2 {\bf \Omega}_2, \mat{U}_3 {\bf \Omega}_3, 
 - (\mathbfcal{G}{\times_1} {\bf \Omega}_1  + 
\mathbfcal{G}{{\times_2}} {\bf \Omega}_2 +
\mathbfcal{G}{{\times_3}} {\bf \Omega}_3)): {\bf \Omega}_d \in \mathbb{R}^{r_d \times r_d}, {\bf \Omega}_d ^T = -{\bf \Omega}_d \text{ for\ } d \in \{1, 2, 3\} \}$. Skew symmetric matrices ${\bf \Omega}_{d}$ for all $d \in \{1,2,3 \}$ parameterize the vertical space. Finally, the horizontal projection operator $\Pi_{{x}}: T_{{x}} {\mathcal{M}} : \rightarrow \mathcal{H}_{{x}} : $
$\eta_{{x}} \mapsto \Pi_{{x}}(\eta_{{x}})$ has the expression
\begin{equation*}
\begin{array}{lll}
\Pi_{{x}}(\eta_{{x}} )
= (
\eta_{\scriptsize \mat{U}_{1}} \!\!- \mat{U}_{1} {\bf \Omega}_{1},
\eta_{\scriptsize \mat{U}_{2}} \!\!- \mat{U}_{2} {\bf \Omega}_{2}, 
\eta_{\scriptsize \mat{U}_{3}} \!\!- \mat{U}_{3} {\bf \Omega}_{3}, 
\eta_{\scriptsize \mathbfcal{G}} \!-\! (\! - (\mathbfcal{G}{\times_1} {\bf \Omega}_{1} \! + \!
\mathbfcal{G}{{\times_2}} {\bf \Omega}_{2}\! +\!
\mathbfcal{G}{{\times_3}} {\bf \Omega}_{3})) 
),
\end{array}
\end{equation*}
where $\eta_x = (\eta_{\scriptsize \mat{U}_{1}}, \eta_{\scriptsize \mat{U}_{2}}, \eta_{\scriptsize \mat{U}_{3}}, \eta_{\scriptsize \mathbfcal{G}}) \in T_x \mathcal{M}$ and  ${\bf \Omega}_{d}$ is a skew-symmetric matrix of size $r_d \times r_d$ that is the solution to the \emph{coupled} Lyapunov equations
\begin{equation}
\begin{array}{lll}
\label{Eq:OmegaRequirements}
\left\{
\begin{array}{l}
\mat{G}_{1}  \mat{G}_{1}^T {\bf \Omega}_{1} + {\bf \Omega}_{1} \mat{G}_{1}  \mat{G}_{1}^T 
-\mat{G}_{1}(\mat{I}_{r_3} \otimes {\bf \Omega}_{2}) \mat{G}_{1}^T 
- \mat{G}_{1}( {\bf \Omega}_{3} \otimes \mat{I}_{r_2} )\mat{G}_{1}^T  \\
\hspace{6cm} = {\rm Skew}(\mat{U}_1^T\eta_{\scriptsize \mat{U}_1}\mat{G}_{1}  \mat{G}_{1}^T) + {\rm Skew}(\mat{G}_{1}\eta_{\scriptsize \mat{G}_{1}}^T), \\
\mat{G}_{2}  \mat{G}_{2}^T {\bf \Omega}_{2} + {\bf \Omega}_{2} \mat{G}_{2}  \mat{G}_{2}^T 
-\mat{G}_{2}(\mat{I}_{r_3} \otimes {\bf \Omega}_{1}) \mat{G}_{2}^T 
- \mat{G}_{2}( {\bf \Omega}_{3} \otimes \mat{I}_{r_1} )\mat{G}_{2}^T  \\
\hspace{6cm} = {\rm Skew}(\mat{U}_2^T\eta_{\scriptsize \mat{U}_2}\mat{G}_{2}  \mat{G}_{2}^T) + {\rm Skew}(\mat{G}_{2}\eta_{\scriptsize \mat{G}_{2}}^T), \\
\mat{G}_{3}  \mat{G}_{3}^T {\bf \Omega}_{3} + {\bf \Omega}_{3} \mat{G}_{3}  \mat{G}_{3}^T 
-\mat{G}_{3}(\mat{I}_{r_2} \otimes {\bf \Omega}_{1}) \mat{G}_{3}^T 
- \mat{G}_{3}( {\bf \Omega}_{2} \otimes \mat{I}_{r_1} )\mat{G}_{3}^T  \\
\hspace{6cm} = {\rm Skew}(\mat{U}_3^T\eta_{\scriptsize \mat{U}_3}\mat{G}_{3}  \mat{G}_{3}^T) + {\rm Skew}(\mat{G}_{3}\eta_{\scriptsize \mat{G}_{3}}^T),
\end{array}
\right.
\end{array}
\end{equation}
where ${\rm Skew}(\cdot)$ extracts the skew-symmetric part of a square matrix, i.e., ${\rm Skew}(\mat{D})=(\mat{D}-\mat{D}^T)/2$. \changeBM{The coupled Lyapunov equations (\ref{Eq:OmegaRequirements}) are solved efficiently with the Matlab's} \verb+pcg+ \changeBM{routine that is combined with a specific preconditioner resulting from the Gauss-Seidel approximation of (\ref{Eq:OmegaRequirements}).}

{\bf Retraction.}
A retraction is a mapping that maps vectors in the horizontal space to points on the search space $\mathcal{M}$ and 
satisfies the local rigidity condition \cite[Definition~4.1]{Absil_OptAlgMatManifold_2008}. It provides a natural way to move on the manifold along a search direction. Because the total space  $\mathcal{M}$ has the product nature, we can choose a retraction by combining retractions on the individual manifolds, i.e.,
\begin{eqnarray*}
\label{eq:Retraction}
R_{x} (\xi_x)  =   ({\rm uf}(\mat{U}_{1}+\xi_{\scriptsize \mat{U}_{1}}), {\rm uf}(\mat{U}_{2}+\xi_{\scriptsize \mat{U}_{2}}), 
{\rm uf}(\mat{U}_{3}+\xi_{\scriptsize \mat{U}_{3}}), \mathbfcal{G}+\xi_{\scriptsize \mathbfcal{G}}),
\end{eqnarray*}
where $\xi_x \in \mathcal{H}_x$ and ${\rm uf}(\cdot)$ extracts the orthogonal factor of a full column rank matrix, i.e., 
${\rm uf}(\mat{A})=\mat{A}(\mat{A}^T\mat{A})^{-1/2}$. The retraction $ {R}_{  x}$ defines a retraction ${R}_{[x]}( {\xi}_{[x]}) : =[R_x (\xi_x)] $ on the quotient manifold ${\mathcal{M}}/\sim$, as the equivalence class $[R_x (\xi_x)] $ does not depend on specific matrix representations of $[x]$ and ${\xi}_{[x]}$, where $ {\xi}_{  x}$ is the horizontal lift of the abstract tangent vector $\xi_{[x]} \in T_{[x]} (\mathcal{M} /\sim)$.

\changeBM{{\bf Vector transport.}
A vector transport $\mathcal{T}_{\eta_x} \xi_x$ on a manifold $\mathcal{M}$ is a smooth mapping that transports a tangent vector $\xi_x \in T_x \mathcal{M}$ at $x \in \mathcal{M}$ to a vector in the tangent space at $R_x(\eta_x)$ \cite[Section~8.1.4]{Absil_OptAlgMatManifold_2008}. It generalizes the classical concept of translation of vectors in the Euclidean space to manifolds. The horizontal lift of the abstract vector transport $\mathcal{T}_{\eta_{[x]}} \xi_{[x]}$ on $\mathcal{M}/\!\sim$ has the matrix characterization $\Pi_{R_x(\eta_x)}(\mathcal{T}_{\eta_x} \xi_x) = \Pi_{R_x(\eta_x)}(\Psi_{R_x(\eta_x)}(\xi_x))$, where $\xi_x$ and $\eta_x$ are the horizontal lifts in $\mathcal{H}_x$ of $\xi_{[x]}$ and $\eta_{[x]}$ that belong to 
$T_{[x]}(\mathcal{M}/ \!\sim)$. The computational cost of transporting a vector solely depends on the projection and retraction operations.}

\section{\changeBM{Preconditioned conjugate gradient algorithm for (\ref{Eq:CostFunction})}}
\label{sec:AlgorithmDetails}
We propose a Riemannian nonlinear conjugate gradient algorithm for the tensor completion problem (\ref{Eq:CostFunction}) that is based on the developments in Section \ref{sec:OptimizationRelatedIngredients}. The preconditioning effect follows from the specific choice of the metric (\ref{Eq:metric}). The earlier developments allow to use the off-the-shelf conjugate gradient implementation of Manopt for any smooth cost function \cite{Boumal_Manopt_2014_s}. A complete description of the Riemannian nonlinear conjugate gradient method is in \cite[Chapter~8]{Absil_OptAlgMatManifold_2008}. The convergence analysis of the Riemannian conjugate gradient method follows from \cite{Sato15a, Ring_SIAMJOptim_2012_s}. The only remaining ingredients are the cost function specific ingredients. To this end, we show the computation of the Riemannian gradient as well as a way to compute an initial guess for the step-size, which is used in the conjugate gradient method. The \changeHK{concrete} formulas are shown in Table \ref{tab:FinalFormulas}. The total computational cost per iteration of our proposed algorithm is \changeBM{$O(|\Omega| r_1 r_2 r_3)$}, where $|\Omega|$ is the number of known entries.

\begin{table}[h]
\begin{center}  \small 
\caption{Ingredients to implement an off-the-shelf conjugate gradient algorithm for (\ref{Eq:CostFunction}).}
\label{tab:FinalFormulas}
\begin{tabular}{l|l}
\hline
Matrix representation  & $x =  
(\mat{U}_{1}, \mat{U}_{2}, \mat{U}_{3}, \mathbfcal{G})$
\\ 
\hdashline
Computational space $\mathcal{M}$ &
${\rm St}(r_1, n_1) \times {\rm St}(r_2, n_2) \times {\rm St}(r_3, n_3) \times \mathbb{R}^{r_1 \times r_2 \times r_3}$
\\ \hdashline
Group action & 
$\{ (\mat{U}_{1}\mat{O}_{1}, \mat{U}_{2}\mat{O}_{2}, \mat{U}_{3}\mat{O}_{3}, \mathbfcal{G}{\times_1} \mat{O}^T_{1}{{\times_2}} \mat{O}^T_{2}{{\times_3}} \mat{O}^T_{3}): $
\\
& $ \mat{O}_{d} \in \mathcal{O}{(r_d)}, \text{for\ }d \in \{1,2,3\} \}$
\\ \hdashline
Quotient space $\mathcal{M}/\!\sim$ & 
 ${\rm St}(r_1, n_1) \times {\rm St}(r_2, n_2) \times {\rm St}(r_3, n_3) \times \mathbb{R}^{r_1 \times r_2 \times r_3}$\\
 & $/ (\mathcal{O}{(r_1)} \times \mathcal{O}{(r_2)} \times \mathcal{O}{(r_3)})$
\\ 
\hdashline
\changeBM{Ambient space} &
$\mathbb{R}^{ n_1 \times r_1} \times \mathbb{R}^{ n_2 \times r_2} \times \mathbb{R}^{ n_3 \times r_3} \times \mathbb{R}^{r_1 \times r_2 \times r_3}$
\\
\hdashline
Tangent vectors in & 
$
\{ (\mat{Z}_{{\bf U}_{1}}, \mat{Z}_{{\bf U}_{2}}, \mat{Z}_{{\bf U}_{3}}, \mat{Z}_{\mathbfcal{G}}) \in 
\mathbb{R}^{n_1 \times r_1} \times  
\mathbb{R}^{n_2 \times r_2} \times 
\mathbb{R}^{n_3 \times r_3} \times 
\mathbb{R}^{r_1 \times r_2 \times r_3} $
\\
$T_x \mathcal{M}$ & : $\mat{U}_{d}^T \mat{Z}_{{\bf U}_{d}} +  \mat{Z}_{{\bf U}_{d}}^T \mat{U}_{d} = 0, \text{ for }d \in \{1,2,3\} \}
$
\\ \hdashline
Metric ${g}_{x}(\xi_{x}, \eta_{x})$ for & 
$\langle \xi_{\scriptsize \mat{U}_{1}},
{\eta}_{\scriptsize\mat{U}_{1}} (\mat{G}_{1} \mat{G}_{1}^T) \rangle \! + \!
\langle \xi_{\scriptsize \mat{U}_{2}},
{\eta}_{\scriptsize\mat{U}_{2}} (\mat{G}_{2} \mat{G}_{2}^T) \rangle \! + \!
\langle \xi_{\scriptsize \mat{U}_{3}},
{\eta}_{\scriptsize\mat{U}_{3}} (\mat{G}_{3} \mat{G}_{3}^T) \rangle 
\! + \! \langle {\xi}_{\scriptsize \mathbfcal{G}}, {\eta}_{\scriptsize \mathbfcal{G}}\rangle $
\\
any $\xi_x, \eta_x \in T_x \mathcal{M}$ & 
\\ \hdashline
Vertical tangent  & 
$\{ (\mat{U}_{1} {\bf \Omega}_{1}, \mat{U}_{2}  {\bf \Omega}_{2}, \mat{U}_{3}  {\bf \Omega}_{3}, 
 - (\mathbfcal{G}{\times_1} {\bf \Omega}_{1}  + 
\mathbfcal{G}{{\times_2}} {\bf \Omega}_{2} +
\mathbfcal{G}{{\times_3}} {\bf \Omega}_{3})):$
\\
vectors in $\mathcal{V}_x$ & $ {\bf \Omega}_{d} \in \mathbb{R}^{r_d \times r_d}, {\bf \Omega}_{d}^T = -{\bf \Omega}_{d}, \text{for }d \in \{1,2,3\} \}$
\\ \hdashline
Horizontal tangent  & 
$\{(\zeta_{\mat{U}_{1}}, \zeta_{\mat{U}_{2}}, \zeta_{\mat{U}_{3}}, \zeta_{\mathbfcal{G}}) \in T_x \mathcal{M} : $
\\
vectors in $\mathcal{H}_x$& $(\mat{G}_{d}  \mat{G}_{d}^T) \zeta_{\scriptsize \mat{U}_{d}}^T \mat{U}_{d}  
+ \zeta_{\scriptsize \mat{G}_{d}}  \mat{G}_{d}^T \text{ is symmetric}, \text{for }d \in \{1,2,3\} \}$
\\ \hdashline
$\Psi(\cdot)$ projects \changeBM{an ambient} & 
$(
\mat{Y}_{\scriptsize \mat{U}_{1}} - \mat{U}_{1} \mat{S}_{\scriptsize \mat{U}_{1}} (\mat{G}_{1}  \mat{G}_{1}^T)^{-1},
\mat{Y}_{\scriptsize \mat{U}_{2}} - \mat{U}_{2} \mat{S}_{\scriptsize \mat{U}_{2}} (\mat{G}_{2}  \mat{G}_{2}^T)^{-1}, 
$
\\
 vector  $(\mat{Y}_{\scriptsize \mat{U}_{1}}, \mat{Y}_{\scriptsize \mat{U}_{2}}, \mat{Y}_{\scriptsize \mat{U}_{3}}, \mat{Y}_{\scriptsize \mathbfcal{G}} )$& $\mat{Y}_{\scriptsize \mat{U}_{3}} - \mat{U}_{3} \mat{S}_{\scriptsize \mat{U}_{3}} (\mat{G}_{3}  \mat{G}_{3}^T)^{-1},
\mat{Y}_{\scriptsize \mathbfcal{G}}  
)$, where $\mat{S}_{\scriptsize \mat{U}_{d}}$ for $d \in \{1, 2,3 \}$ are computed
\\
onto $T_x \mathcal{M}$&  by solving Lyapunov equations as in (\ref{Eq:B_Requirements}).
%
\\ \hdashline
$\Pi(\cdot)$ projects a tangent  &  $(
\xi_{\scriptsize \mat{U}_{1}} - \mat{U}_{1} {\bf \Omega}_{1},
\xi_{\scriptsize \mat{U}_{2}} - \mat{U}_{2} {\bf \Omega}_{2}, 
\xi_{\scriptsize \mat{U}_{3}} - \mat{U}_{3} {\bf \Omega}_{3}, $ 
\\
vector $\xi$ onto $\mathcal{H}_x$ & 
$\xi_{\scriptsize \mathbfcal{G}} - ( - (\mathbfcal{G}{\times_1} {\bf \Omega}_{1}  + 
\mathbfcal{G}{{\times_2}} {\bf \Omega}_{2} +
\mathbfcal{G}{{\times_3}} {\bf \Omega}_{3})) 
)$, ${\bf \Omega}_{d}$ is computed in (\ref{Eq:OmegaRequirements}).
\\ \hdashline
First-order derivative & 
$(\mat{S}_{1} (\mat{U}_{3} \otimes \mat{U}_{2}) \mat{G}_{1}^T, \mat{S}_{2} (\mat{U}_{3} \otimes \mat{U}_{1}) \mat{G}_{2}^T, \mat{S}_{3} (\mat{U}_{2} \otimes \mat{U}_{1}) \mat{G}_{3}^T), $ 
\\
of $f(x)$ & $\mathbfcal{S} \times_1 \mat{U}_{1}^T \times_2 \mat{U}_{2}^T \times_3 \mat{U}_{3}^T),$
\\
 & where $\mathbfcal{S} = \frac{2}{|{\Omega} |} (\mathbfcal{P}_{\Omega}(\mathbfcal{G}{\times_1} {\mat{U}_{1}}{\times_2} {\mat{U}_{2}}{\times_3} {\mat{U}_{3}}) - 
\mathbfcal{P}_{\Omega}(\mathbfcal{X}^{\star}))$.
\\ \hdashline
Retraction $R_x(\xi_x)$ & 
$({\rm uf}(\mat{U}_{1}+\xi_{\scriptsize \mat{U}_{1}}), {\rm uf}(\mat{U}_{2}+\xi_{\scriptsize \mat{U}_{2}}), {\rm uf}(\mat{U}_{3}+\xi_{\scriptsize \mat{U}_{3}}), \mathbfcal{G}+\xi_{\scriptsize \mathbfcal{G}})$
\\ \hdashline
\changeBM{Horizontal lift of the} & 
$\Pi_{R_x(\eta_x)}(\Psi_{R_x(\eta_x)}(\xi_x))$
\\
vector transport $\mathcal{T}_{\eta_{[x]}} \xi_{[x]}$  & \\
\hline
\end{tabular}
\end{center}
\end{table}

{\bf Riemannian gradient computation.} Let $f(\mathbfcal{X})=\| \mathbfcal{P}_{\Omega}(\mathbfcal{X}) - \mathbfcal{P}_{\Omega}(\mathbfcal{X}^{\star}) \|^2_F/|\Omega |$ be the mean square error function of (\ref{Eq:CostFunction}), and 
$\mathbfcal{S} = 2 (\mathbfcal{P}_{\Omega}(\mathbfcal{G}{\times_1} \mat{U}_{1} {\times_2} \mat{U}_{2}{\times_3} \mat{U}_{3}) - 
\mathbfcal{P}_{\Omega}(\mathbfcal{X}^{\star}))/|{\Omega}|$ 
be an auxiliary sparse tensor variable that is interpreted as the Euclidean gradient of $f$ in $\mathbb{R}^{n_1 \times n_2 \times n_3}$. The partial derivatives of the function $f$ with respect to $(\mat{U}_{1}, \mat{U}_{2}, \mat{U}_{3}, \mathbfcal{G})$ are computed in terms of the unfolding matrices $\mat{S}_{d}$. Due to the specific scaled metric (\ref{Eq:metric}), the partial derivatives are further scaled by $((\mat{G}_{1}\mat{G}_{1}^T)^{-1}, (\mat{G}_{2}\mat{G}_{2}^T)^{-1}, (\mat{G}_{3}\mat{G}_{3}^T)^{-1}, \mathbfcal{I})$, denoted as ${\rm egrad}_{x} f$ (after scaling). Finally, from the Riemannian submersion theory \cite[Section ~3.6.2]{Absil_OptAlgMatManifold_2008}, the horizontal lift of ${\rm grad}_{[x]}f $ is equal to ${\rm grad}_{x}f \ =\ \Psi({\rm egrad}_{x} f )$. Subsequently,
\begin{equation*}
\begin{array}{lll}
\text{the horizontal lift of\ }{\rm grad}_{[x]} f & = &  
(\mat{S}_{1} (\mat{U}_{3} \otimes \mat{U}_{2}) \mat{G}_{1}^T (\mat{G}_{1}\mat{G}_{1}^T)^{-1} 
- \mat{U}_{1} \mat{B}_{\scriptsize \mat{U}_{1}}(\mat{G}_{1}\mat{G}_{1}^T)^{-1}, \\
&  & \hspace{0cm} \mat{S}_{2} (\mat{U}_{3} \otimes \mat{U}_{1}) \mat{G}_{2}^T (\mat{G}_{2}\mat{G}_{2}^T)^{-1}
- \mat{U}_{2} \mat{B}_{\scriptsize \mat{U}_{2}}(\mat{G}_{2}\mat{G}_{2}^T)^{-1}, \\
&  & \hspace{0cm} \mat{S}_{3} (\mat{U}_{2} \otimes \mat{U}_{1}) \mat{G}_{3}^T(\mat{G}_{3}\mat{G}_{3}^T)^{-1}
- \mat{U}_{3} \mat{B}_{\scriptsize \mat{U}_{3}}(\mat{G}_{3}\mat{G}_{3}^T)^{-1}, \\
& & \mathbfcal{S} \times_1 \mat{U}_{1}^T \times_2 \mat{U}_{2}^T \times_3 \mat{U}_{3}^T), 
\end{array}
\end{equation*}
where $\mat{B}_{\scriptsize \mat{U}_{d}}$ for $d \in \{1, 2,3\} $ are the solutions to the Lyapunov equations
\begin{equation*}
\left\{
\begin{array}{lll}
\mat{B}_{\scriptsize \mat{U}_{1}} \mat{G}_{1}  \mat{G}_{1}^T + \mat{G}_{1}  \mat{G}_{1}^T \mat{B}_{\scriptsize \mat{U}_{1}}
& = & 2 {\rm Sym} (\mat{G}_{1}  \mat{G}_{1}^T\mat{U}_{1}^T (\mat{S}_{1} (\mat{U}_{3} \otimes \mat{U}_{2}) \mat{G}_{1}^T), \\ 
\mat{B}_{\scriptsize \mat{U}_{2}} \mat{G}_{2}  \mat{G}_{2}^T + \mat{G}_{2}  \mat{G}_{2}^T \mat{B}_{\scriptsize \mat{U}_{2}}
& = & 2 {\rm Sym} (\mat{G}_{2}  \mat{G}_{2}^T\mat{U}_{2}^T (\mat{S}_{2} (\mat{U}_{3} \otimes \mat{U}_{1}) \mat{G}_{2}^T), \\ 
\mat{B}_{\scriptsize \mat{U}_{3}} \mat{G}_{3}  \mat{G}_{3}^T + \mat{G}_{3}  \mat{G}_{3}^T \mat{B}_{\scriptsize \mat{U}_{3}}
& = & 2 {\rm Sym} (\mat{G}_{3}  \mat{G}_{3}^T\mat{U}_{3}^T (\mat{S}_{3} (\mat{U}_{2} \otimes \mat{U}_{1}) \mat{G}_{3}^T),
\end{array}
\right.
\end{equation*}
which are solved efficiently with the Matlab's \verb+lyap+ routine. ${\rm Sym}(\cdot)$ extracts the symmetric part of a square matrix, i.e., ${\rm Sym}(\mat{D})=(\mat{D}+\mat{D}^T)/2$. \changeBM{The total numerical cost of computing the Riemannian gradient depends on computing the partial derivatives, which is $O(|\Omega| r_1 r_2 r_3)$.}

{\bf Initial guess for the step size.} The least-squares structure of the cost function in (\ref{Eq:CostFunction}) also allows to compute a  \emph{linearized} step-size guess efficiently along a search direction by considering a polynomial approximation of degree $2$ over the manifold \cite{Mishra_ICDC_2014_s, Vandereycken_SIAMOpt_2013_s}. Given a search direction $\xi_x \in \mathcal{H}_x$, the step-size guess is 
$\argmin_{s \in \mathbb{R}_{+}}
\| \mathbfcal{P}_{\Omega}(
\mathbfcal{G}{\times_1} \mat{U}_{1}{\times_2} \mat{U}_{2}{\times_3} \mat{U}_{3} 
+ s \mathbfcal{G}{\times_1} {\xi_{\scriptsize \mat{U}_{1}}}{\times_2} \mat{U}_{2}{\times_3} \mat{U}_{3} 
+ s \mathbfcal{G}{\times_1} \mat{U}_{1}{\times_2} {\xi_{\scriptsize \mat{U}_{2}}}{\times_3} \mat{U}_{3} 
+ s \mathbfcal{G}{\times_1} \mat{U}_{1}{\times_2}  \mat{U}_{2}{\times_3} {\xi_{\scriptsize \mat{U}_{3} }}
+ s {\xi_{\scriptsize \mathbfcal{G}}}{\times_1} \mat{U}_{1}{\times_2} \mat{U}_{2}{\times_3} \mat{U}_{3} 
) 
- \mathbfcal{P}_{\Omega}(\mathbfcal{X}^{\star}) \|^2_F$, which has a closed-form expression and the numerical cost of computing it is \changeBM{$O(|\Omega| r_1 r_2 r_3)$}.

\section{Numerical comparisons}\label{sec:NumericalComparisons}
We show a number of numerical comparisons of our proposed Riemannian preconditioned nonlinear conjugate algorithm with state-of-the-art algorithms that include TOpt \cite{Filipovi_MultiSysSigPro_2013_s} and geomCG \cite{Kressner_BIT_2014_s}, for comparisons with Tucker decomposition based algorithms, and  HaLRTC \cite{Liu_IEEETransPAMI_2013_s}, Latent \cite{Tomioka_Latent_2011_s}, and Hard \cite{Signoretto_MachineLearning_2014_s} as nuclear norm minimization algorithms. All simulations are performed in Matlab on a 2.6 GHz Intel Core i7 machine with 16 GB RAM. For specific operations with unfoldings of $\mathbfcal{S}$, we use the \verb+mex+ interfaces for Matlab that are provided by the authors of geomCG. For large-scale instances, our algorithm is only compared with geomCG as other algorithms cannot handle these instances.

Since the dimension of the space of a tensor $\in \mathbb{R}^{n_1 \times n_2 \times n_3}$ of rank ${\bf r} = (r_1, r_2, r_3)$ is $\text{dim}(\mathcal{M}/\!\sim) = \sum_{d=1}^3 (n_d r_d - r_d^2) +  r_1 r_2 r_3$, we randomly and uniformly select known entries based on a multiple of the dimension, called the \emph{over-sampling} (OS) ratio, to create the training set $\Omega$. Algorithms (and problem instances) are initialized randomly, as in \cite{Kressner_BIT_2014_s}, and are stopped when either the mean square error (MSE) on the training set $\Omega$ is below $10^{-12}$ or the number of iterations exceed\changeHK{s} $250$. We also evaluate the mean square error on a test set $\Gamma$, which is different from $\Omega$. Five runs are performed in each scenario and the plots show all of them. The time plots are shown with standard deviations. 



\begin{figure}[h]
\begin{tabular}{cc}
\begin{minipage}{0.32\hsize}
\begin{center}
\includegraphics[width=\hsize]{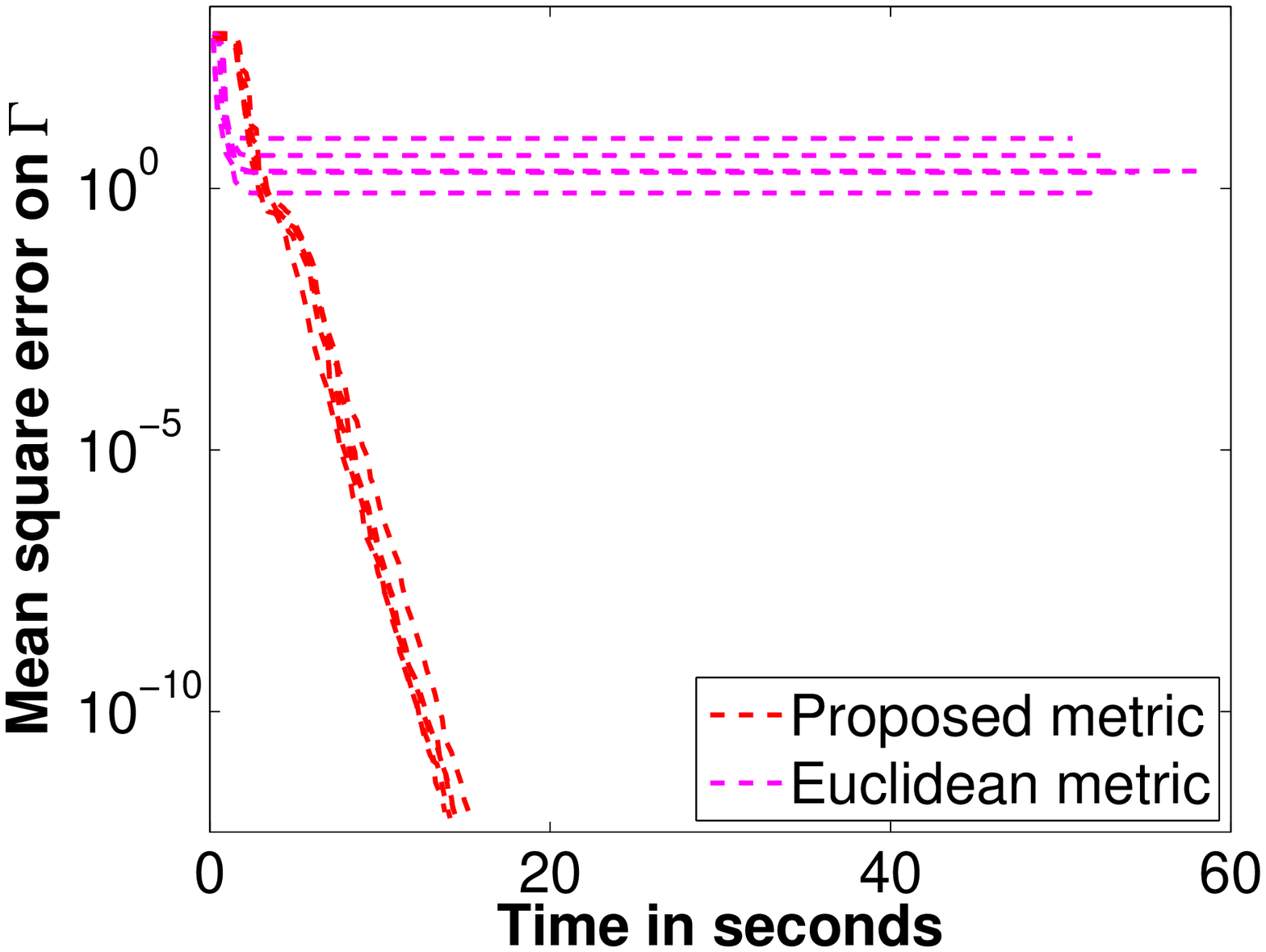}\\
{\scriptsize(a) {\bf Case S1:} comparison between metrics.}
\end{center}
\end{minipage}
\begin{minipage}{0.32\hsize}
\begin{center}
\includegraphics[width=\hsize]{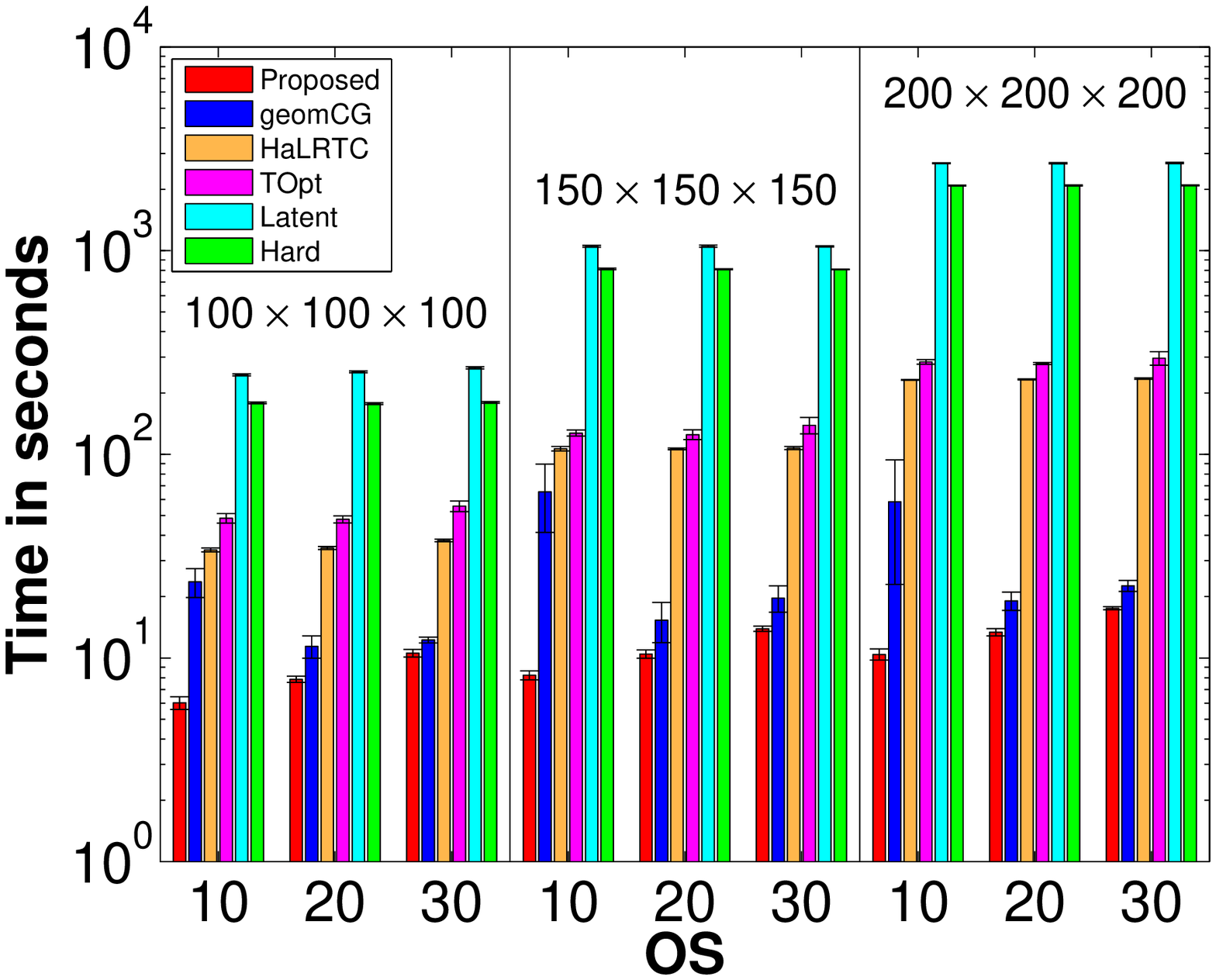}\\
{\scriptsize(b) {\bf Case S2:} {\bf r} = $(10,10,10)$.}
\end{center}
\end{minipage}
\begin{minipage}{0.32\hsize}
\begin{center}
\includegraphics[width=\hsize]{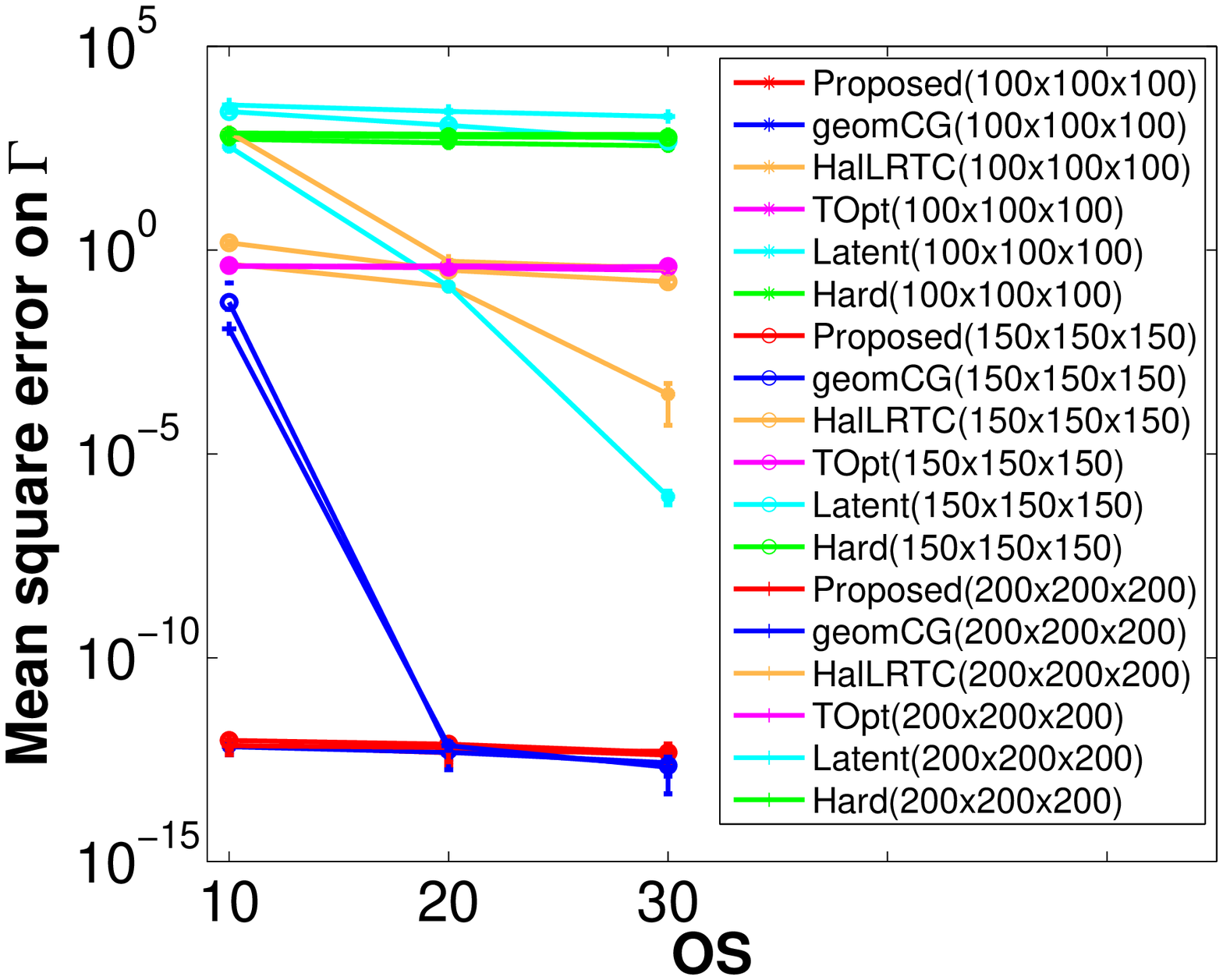}\\
{\scriptsize(c) {\bf Case S2:} {\bf r} = $(10,10,10)$.}
\end{center}
\end{minipage}\\
\begin{minipage}{0.32\hsize}
\begin{center}
\includegraphics[width=\hsize]{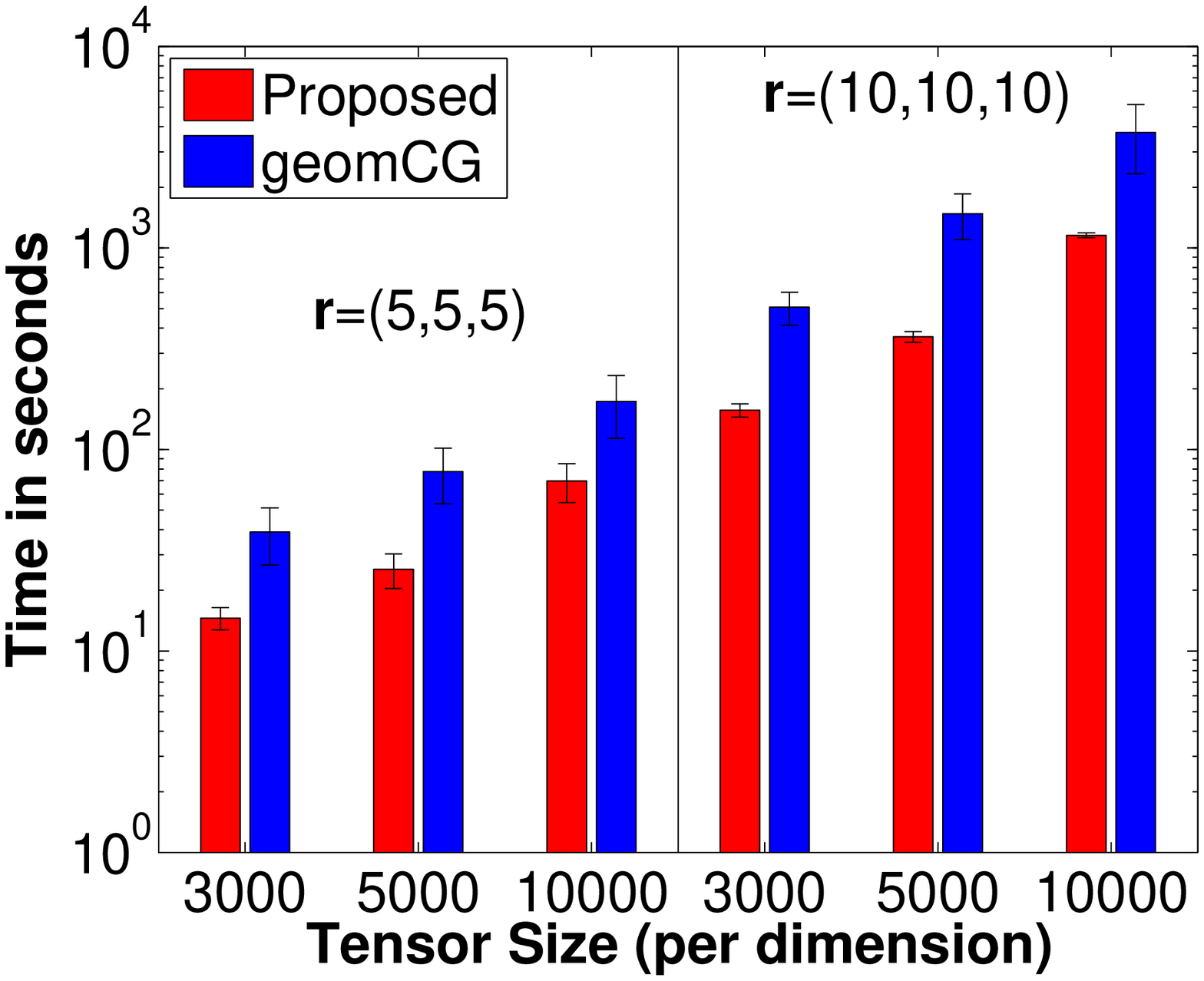}\\
{\scriptsize(d) {\bf Case S3}.}
\end{center}
\end{minipage}
\begin{minipage}{0.32\hsize}
\begin{center}
\includegraphics[width=\hsize]{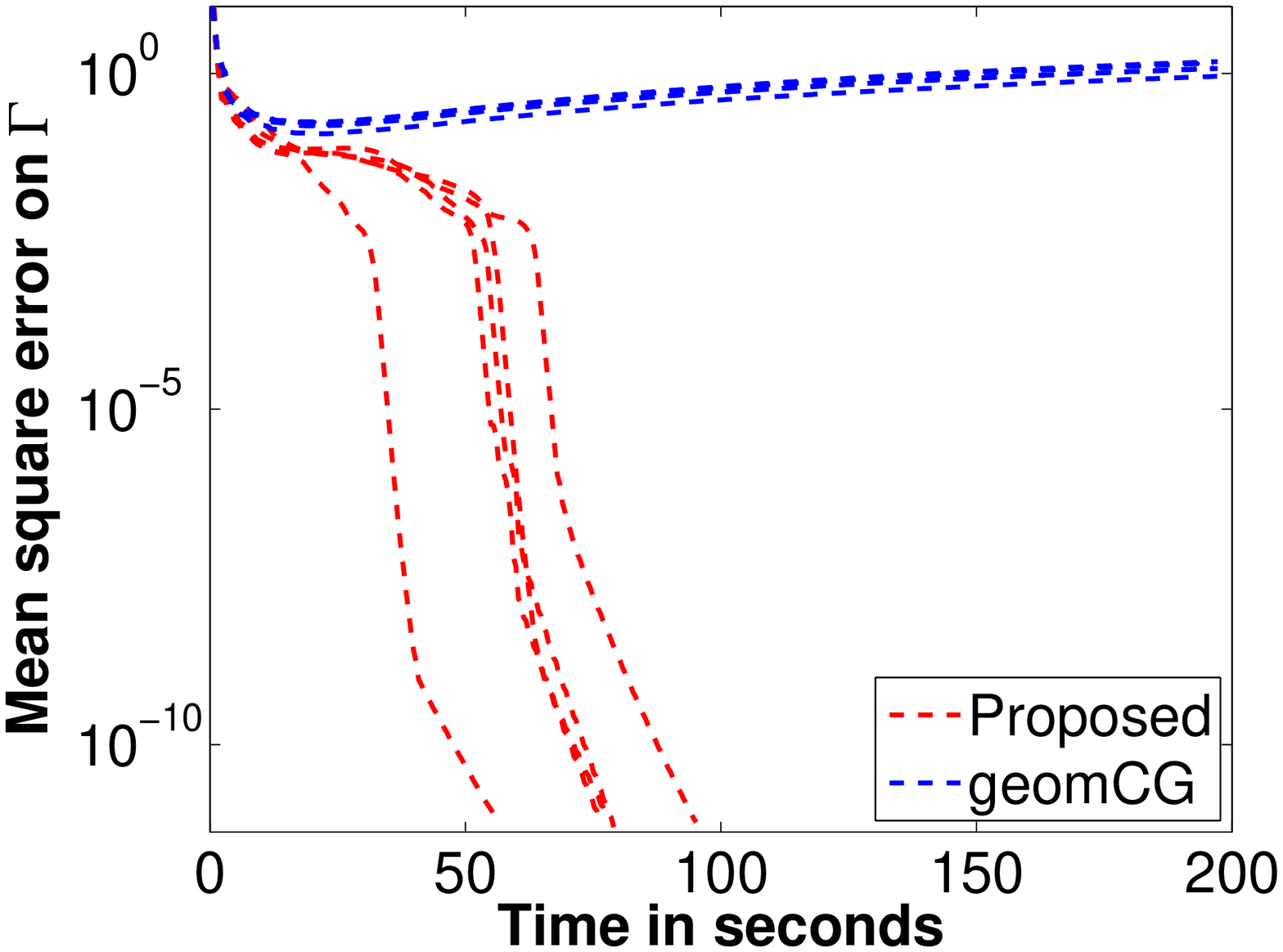}\\
{\scriptsize(e) {\bf Case S4:} OS = $4$.}
\end{center}
\end{minipage}
\begin{minipage}{0.32\hsize}
\begin{center}
\includegraphics[width=\hsize]{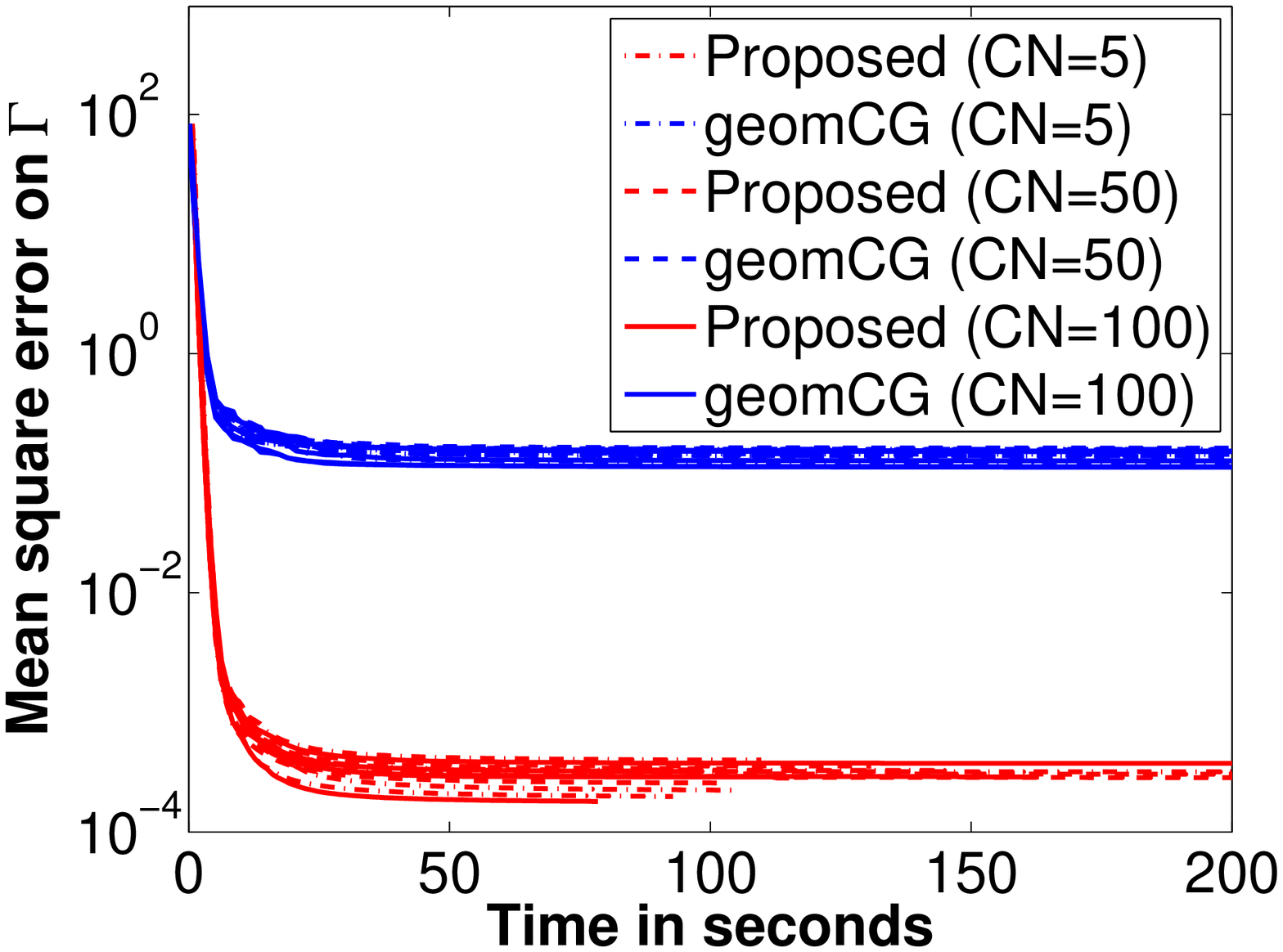}\\
{\scriptsize \changeHK{(f)} {\bf Case S5:} CN = $\{5,50,100\}$.}
\end{center}
\end{minipage}\\
\begin{minipage}{0.32\hsize}
\begin{center}
\includegraphics[width=\hsize]{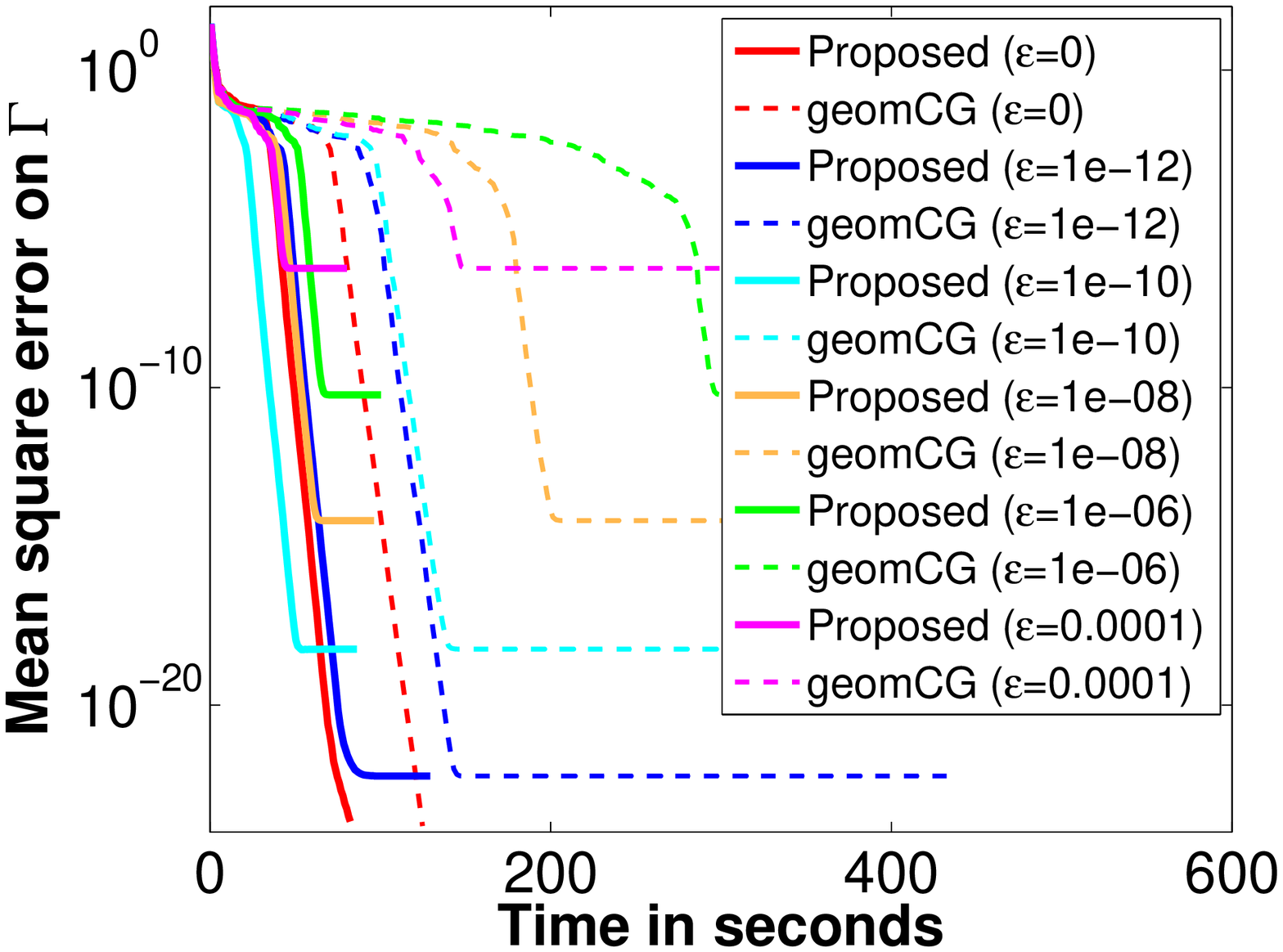}\\
{\scriptsize \changeHK{(g)} {\bf Case S6:} noisy data.}
\end{center}
\end{minipage}
\begin{minipage}{0.32\hsize}
\begin{center}
\includegraphics[width=\hsize]{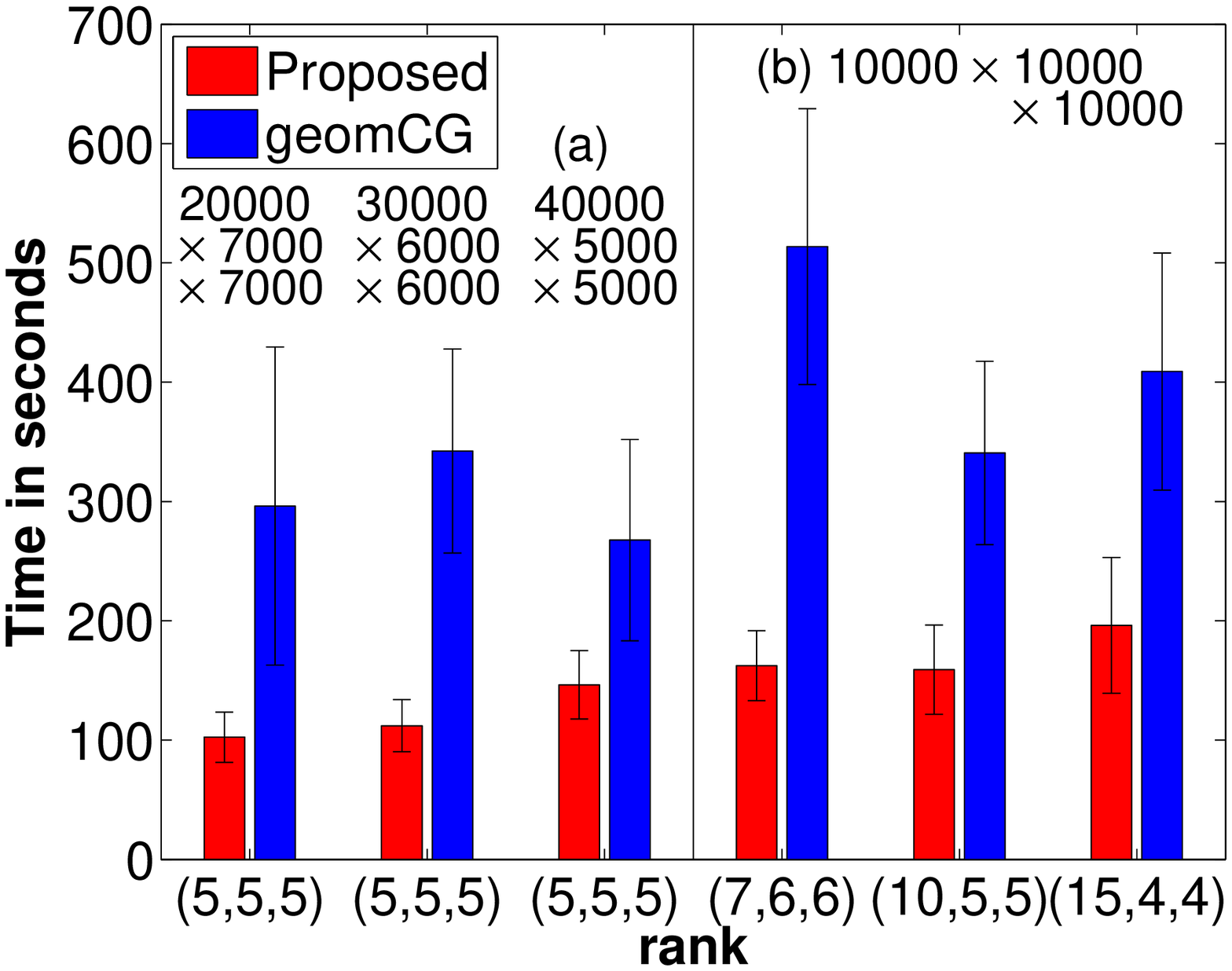}\\
{\scriptsize \changeHK{(h)} {\bf Case S7:} asymmetric tensors.}
\end{center}
\end{minipage}
\begin{minipage}{0.32\hsize}
\begin{center}
\includegraphics[width=\hsize]{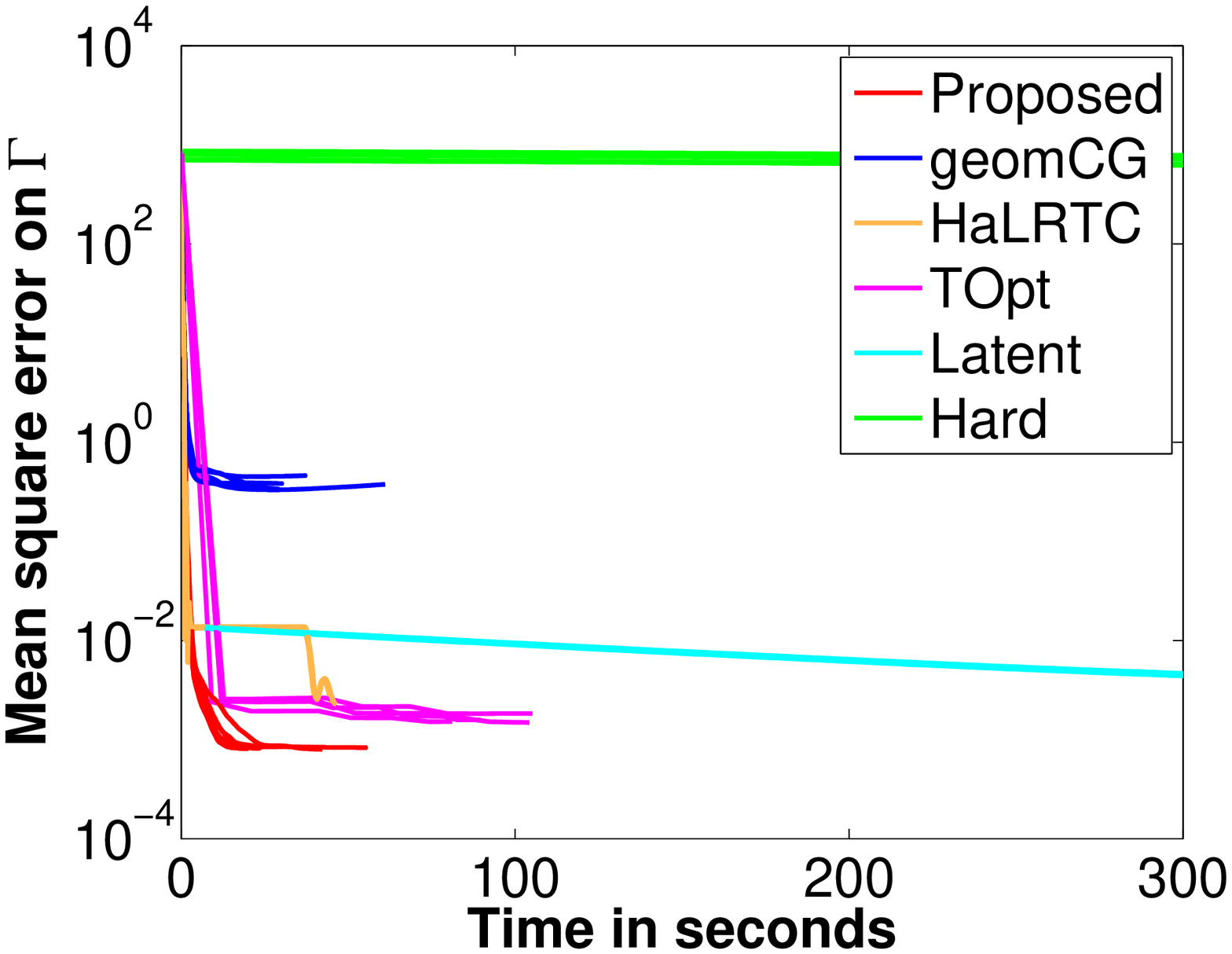}\\
{\scriptsize(i) {\bf Case R1:} Ribeira, OS = 11.}
\end{center}
\end{minipage}
\end{tabular}
\caption{Experiments on synthetic and real datasets.}
\label{fig:syntheticReal}
\end{figure}

{\bf Case S1: comparison with the Euclidean metric.} We first show the benefit of the proposed metric (\ref{Eq:metric}) over the conventional choice of the Euclidean metric that exploits the product structure of $\mathcal{M}$ and symmetry (\ref{Eq:EquivalenceClass_3}). This is defined by combining the individual natural metrics for ${\rm St}(r_d,n_d)$ and $\mathbb{R}^{r_1\times r_2\times r_3}$. For simulations, we randomly generate a tensor of size $200 \times 200 \times 200$ and rank  ${\bf r}=(10,10,10)$. OS is $10$. For simplicity, we compare {\it steepest descent} algorithms with Armijo backtracking linesearch for both the metric choices. Figure \ref{fig:syntheticReal}(a) shows that the algorithm with the metric (\ref{Eq:metric}) \changeHK{gives} a superior performance than that of the conventional metric choice. 

{\bf Case S2: small-scale instances.} Small-scale tensors of size 
$100 \times100 \times100$, $150 \times150\times150$, and $200\times200\times200$ and rank ${\bf r}=(10,10,10)$ are considered. OS is $\{10,20,30\}$. Figure \ref{fig:syntheticReal}(b) shows that the convergence behavior of our proposed algorithm is either competitive or faster than the others. In Figure \ref{fig:syntheticReal}(c), the lowest test errors are obtained by our proposed algorithm and geomCG.

{\bf Case S3: large-scale instances.} We consider large-scale tensors of size $3000 \times 3000 \times 3000$, $5000 \times 5000 \times 5000$, and $10000 \times 10000 \times 10000$ and ranks ${\bf r}=(5,5,5)$ and $(10,10,10)$.  OS is $10$. Our proposed algorithm outperforms geomCG in Figure \ref{fig:syntheticReal}(d).


{\bf Case S4: influence of low sampling.} We look into problem instances from scarcely sampled data, e.g., OS is \changeHK{$4$}. The test requires completing a tensor of size $10000 \times 10000 \times 10000$ and rank ${\bf r}=(5,5,5)$. \changeHK{Figure} \ref{fig:syntheticReal}(e) show\changeHK{s} the superior performance of the proposed algorithm against geomCG. Whereas the test error increases for geomCG, it decreases for the proposed algorithm.  

{\bf Case S5: influence of ill-conditioning and low sampling.} We consider the problem instance of {\bf Case S4} with $\text{OS\ } = 5$. Additionally, for generating the instance, we impose a diagonal core $\mathbfcal{G}$ with exponentially decaying \emph{positive} values of condition numbers (CN)  $5$, $50$, and $100$. Figure \ref{fig:syntheticReal}\changeHK{(f)} shows that the proposed algorithm outperforms geomCG for all the considered CN values. 

{\bf Case S6: influence of noise.} We evaluate the convergence properties of algorithms under the presence of noise by adding \emph{scaled} Gaussian noise $\mathbfcal{P}_{\Omega}(\mathbfcal{E})$ to $\mathbfcal{P}_{\Omega}(\mathbfcal{X}^\star)$ as in \cite[Section~4.2.1]{Kressner_BIT_2014_s}. The different noise levels are $\epsilon= \{10^{-4}, 10^{-6}, 10^{-8}, 10^{-10}, 10^{-12}\}$. In order to evaluate for $\epsilon=10^{-12}$, the stopping threshold on the MSE of the train set is lowered to $10^{-24}$. The tensor size and rank are same as in {\bf Case S4} and OS is $10$. Figure \ref{fig:syntheticReal}\changeHK{(g)} shows that the test error for each $\epsilon$ is almost identical to the $\epsilon^2 \| \mathbfcal{P}_{\Omega} (\mathbfcal{X}^{\star})\|_F ^2$ \cite[Section~4.2.1]{Kressner_BIT_2014_s}, but our proposed algorithm converges faster than geomCG.

{\bf Case S7: asymmetric instances.} We consider instances where the dimensions and ranks along certain modes are different than others. Two cases are considered. Case \changeBM{(7.a)} considers tensors size 
$20000 \times7000 \times 7000$, $30000 \times 6000 \times 6000$, and $40000 \times 5000\times 5000$ with \changeHK{rank} ${\bf r}=(5,5,5)$. Case \changeBM{(7.b)} considers a tensor of size $10000 \times10000 \times10000$ with ranks 
$(7,6,6)$, $(10,5,5)$, and $ (15,4,4)$. In all the cases, the proposed algorithm converges faster than geomCG as shown in Figure \ref{fig:syntheticReal}\changeHK{(h)}. 


{\bf Case R1: hyperspectral image.} We consider the hyperspectral image ``Ribeira" \cite{Foster_VisNero_2007_s} discussed in \cite{Signoretto_IEEESigProLetter_2011_s, Kressner_BIT_2014_s}. The tensor size is $1017 \times 1340 \times 33$, where each slice corresponds to an image of a particular scene measured at a different wavelength. As suggested in \cite{Signoretto_IEEESigProLetter_2011_s, Kressner_BIT_2014_s}, we resize it to $203 \times 268 \times 33$. We compare all the algorithms, and perform five random samplings of the pixels based on the OS values $11$ and $22$, corresponding to the rank \changeHKK{{\bf r=}}$(15,15,6)$ adopted in \cite{Kressner_BIT_2014_s}. This set is further randomly split into $80$/$10$/$10$--train/validation/test partitions. The algorithms are stopped when the MSE on the validation set starts to increase. While $\text{OS}=22$ corresponds to the observation ratio of $10\%$ studied in \cite{Kressner_BIT_2014_s}, $\text{OS}=11$ considers a challenging scenario with the observation ratio of $5\%$. Figures \ref{fig:syntheticReal}(i) shows the good performance of our proposed algorithm. Table \ref{tab:R1R2} compiles the results.

{\bf Case R2: MovieLens-10M\footnote{\url{http://grouplens.org/datasets/movielens/}.}.} This dataset contains $10000054$ ratings corresponding to $71567$ users and $10681$ movies. We split the time into $7$-days wide bins results, and finally, get a tensor of size $71567 \times 10681 \times 731$. The fraction of known entries is less than $0.002\%$. The tensor completion task on this dataset reveals {\it periodicity} of the {\it latent} genres. \changeBM{We perform five random $80$/$10$/$10$--train/validation/test partitions.} \changeHK{The maximum iteration threshold is set to $500$.} As shown in Table \ref{tab:R1R2}, our proposed algorithm consistently gives lower test errors than geomCG across different ranks.


\begin{table}[h]
\caption{\changeBM{{\bf Cases R1 and R2:} \changeHK{t}est MSE on $\Gamma$ and time in seconds.}}
\label{tab:R1R2}
\begin{center}
{ \scriptsize 
\changeHK{
\begin{tabular}{l|l|l|l|l}
\hlinewd{1.0pt}
\textbf{Ribeira}& \multicolumn{2}{c|}{OS  = $11$}  & \multicolumn{2}{c}{OS = $22$}  \\
\hdashline
\quad Algorithm & \qquad Time & \qquad \qquad  MSE on $\Gamma$ & \qquad Time & \qquad \qquad  MSE on $\Gamma$  \\
\hline 
Proposed &${\bf 33 \pm 13}$ & ${\bf 8.2095\cdot 10^{-4} \pm 1.7\cdot 10^{-5} }$ & $67 \pm 43$ & ${\bf 6.9516 \cdot 10^{-4} \pm 1.1 \cdot  10^{-5} }$\\
\hdashline
geomCG & $36 \pm 14$ & $3.8342 \cdot  10^{-1} \pm 4.2 \cdot  10^{-2} $ & $150 \pm 48$ & $6.2590 \cdot 10^{-3} \pm 4.5 \cdot  10^{-3} $\\
\hdashline
HaLRTC & $46 \pm  0$ & $2.2671 \cdot  10^{-3} \pm 3.6 \cdot  10^{-5} $ & $ 48 \pm  0$ & $1.3880 \cdot 10^{-3} \pm 2.7 \cdot  10^{-5} $\\
\hdashline
TOpt & $80 \pm 32$ & $1.7854 \cdot  10^{-3} \pm 3.8 \cdot  10^{-4} $ & ${\bf 27 \pm 21}$ & $2.1259 \cdot 10^{-3} \pm 3.8 \cdot  10^{-4} $\\
\hdashline
Latent & $553 \pm 3$ & $2.9296 \cdot 10^{-3} \pm 6.4 \cdot  10^{-5} $ & $558 \pm 3$ & $1.6339  \cdot  10^{-3} \pm 2.3 \cdot  10^{-5} $\\
\hdashline
Hard & $400 \pm 5$ & $6.5090 \cdot 10^{2}\pm 6.1 \cdot  10^{1} $ & $402 \pm 4$ & $6.5989 \cdot  10^{2} \pm 9.8 \cdot  10^{1}$\\
\hlinewd{1.0pt}
\textbf{MovieLens-10M} & \multicolumn{2}{c|}{Proposed}  & \multicolumn{2}{c}{geomCG}  \\
\hdashline
\qquad  {\bf r} & \qquad Time & \qquad \qquad  MSE on $\Gamma$ &  \qquad  Time & \qquad \qquad MSE on $\Gamma$  \\
\hline
$(4,4,4)$ & ${\bf 1748 \pm 441}$ & ${\bf 0.6762 \pm 1.5 \cdot  10^{-3}}$ &  $2981 \pm 40$ & $ 0.6956 \pm?2.8 \cdot 10^{-3}$ \\
\hdashline
$(6,6,6)$ & ${\bf 6058 \pm 47}$ & ${\bf 0.6913 \pm 3.3 \cdot 10^{-3}}$ &  $6554 \pm 655$ & $0.7398 \pm?7.1 \cdot  10^{-3}$ \\
\hdashline
$(8,8,8)$ & ${\bf 11370 \pm 103}$ & ${\bf 0.7589 \pm 7.1 \cdot 10^{-3}}$ &  $13853 \pm 118 $ & $0.8955 \pm?3.3 \cdot 10^{-2}$ \\
\hdashline
$(10,10,10)$ & ${\bf 32802 \pm 52}$ & ${\bf 1.0107 \pm 2.7 \cdot 10^{-2}}$ &  $38145 \pm 36$ & $1.6550 \pm?8.7 \cdot 10^{-2}$ \\
\hlinewd{1.0pt}
\end{tabular}
}
}
\end{center}
\end{table}

\section{Conclusion and future work}
\label{sec:Conclusion}
We have proposed a preconditioned nonlinear conjugate gradient algorithm for the tensor completion problem. The algorithm stems from the Riemannian preconditioning approach that exploits the fundamental structures of symmetry, due to non-uniqueness of Tucker decomposition, and least-squares of the cost function. A novel Riemannian metric (inner product) is proposed that enables to use the versatile Riemannian optimization framework. Concrete matrix expressions are worked out. Numerical comparisons suggest that our proposed algorithm has a superior performance on different benchmarks. As future research directions, we intend to look into ways of updating ranks in tensors as well as look into the issue of preconditioning on other tensor decomposition models, e.g., hierarchical Tucker decomposition \cite{Hackbusch_JRAMFA_2009_s} and tensor networks \cite{Oseledets_SIAMJSC_2009_s}.

\subsubsection*{Acknowledgments}
We thank Rodolphe Sepulchre, Paul Van Dooren, and Nicolas Boumal for useful discussions on the paper. This paper presents research results of the Belgian Network DYSCO (Dynamical Systems, Control, and Optimization), funded by the Interuniversity Attraction Poles Programme, initiated by the Belgian State, Science Policy Office. The scientific responsibility rests with its authors. BM is a research fellow (aspirant) of the Belgian National Fund for Scientific Research (FNRS).




\small
\bibliographystyle{unsrt}  
\bibliography{KMarXiv2016}



\clearpage

\normalsize

\appendix

\renewcommand\thefigure{A.\arabic{figure}}  
\setcounter{figure}{0} 

\renewcommand\thetable{A.\arabic{table}}  
\setcounter{table}{0} 

\renewcommand{\theequation}{A.\arabic{equation}}
\setcounter{equation}{0}

\hrule height 1mm depth 0mm width 140mm
\vspace{0.2cm}

\begin{center}
\changeBMM{{\bf\LARGE Riemannian preconditioning for tensor completion:
\bf\LARGE supplementary material}}
\end{center}
\vspace{0.3cm}

\hrule height 0.2mm depth 0mm width 140mm

\section{Derivation of manifold-related ingredients}
\label{Sup_sec:Derivationofmanifold-relatedingredients}
The concrete computations of the optimization-related ingredientspresented in the paper are discussed below. 

\changeBM{The total space is $\mathcal{M}: = {\rm St}(r_1, n_1) \times {\rm St}(r_2, n_2) \times {\rm St}(r_3, n_3) \times \mathbb{R}^{r_1 \times r_2 \times r_3}$. Each element $x \in \mathcal{M}$ has the matrix representation $(\mat{U}_{1}, \mat{U}_{2}, \mat{U}_{3}, \mathbfcal{G})$. Invariance of Tucker decomposition under the transformation $(\mat{U}_{1}, \mat{U}_{2}, \mat{U}_{3}, \mathbfcal{G}) \mapsto  (  \mat{U}_{1}\mat{O}_{1}, \mat{U}_{2}\mat{O}_{2}, \mat{U}_{3}\mat{O}_{3}, \mathbfcal{G} {\times_1} \mat{O}^T_{1} {\times_2} \mat{O}^T_{2} {\times_3} \mat{O}^T_{3})$ for all $\mat{O}_{d} \in \mathcal{O}(r_d)$, the set of orthogonal matrices of size of $r_d \times r_d$ results in equivalence classes of the form $[x] = [(\mat{U}_{1}, \mat{U}_{2}, \mat{U}_{3}, \mathbfcal{G})] := \{(  \mat{U}_{1}\mat{O}_{1}, \mat{U}_{2}\mat{O}_{2}, \mat{U}_{3}\mat{O}_{3}, \mathbfcal{G} {\times_1} \mat{O}^T_{1} {\times_2} \mat{O}^T_{2} {\times_3} \mat{O}^T_{3}) : \mat{O}_{d} \in \mathcal{O}(r_d)\}$.}

\subsection{\changeHK{Tangent space \changeBMM{characterization} and the Riemannian metric}}
\label{Sup_Sec:TangentspaceandanewRiemannian metric}
\changeHK{
The tangent space, $T_{{x}} \mathcal{M}$, at $x$ given by $(\mat{U}_{1}, \mat{U}_{2}, \mat{U}_{3}, \mathbfcal{G})$ in the total space $\mathcal{M}$ is the product space of the tangent spaces of the individual manifolds. From \cite[Example~3.5.2]{Absil_OptAlgMatManifold_2008}, the tangent space has the matrix characterization 
\begin{equation}
\begin{array}{lll}
\label{Sup_Eq:tangent_space}
T_{{x}} {\mathcal{M}} 
& = &  \{ (\mat{Z}_{{\bf U}_{1}}, \mat{Z}_{{\bf U}_{2}}, \mat{Z}_{{\bf U}_{3}}, \mat{Z}_{\mathbfcal{G}}) \in 
\mathbb{R}^{n_1 \times r_1} \times  
\mathbb{R}^{n_2 \times r_2} \times 
\mathbb{R}^{n_3 \times r_3} \times 
\mathbb{R}^{r_1 \times r_2 \times r_3}  \\
& &
: \mat{U}_{d}^T \mat{Z}_{{\bf U}_{d}} +  \mat{Z}_{{\bf U}_{d}}^T \mat{U}_{d} = 0,\ {\rm for\ } d \in \{1,2,3 \} \}.
\end{array}
\end{equation}
}
\changeHK{
The proposed metric ${g}_{x}:T_x \mathcal{M} \times T_x \mathcal{M} \rightarrow \mathbb{R}$ is
\begin{eqnarray}
\label{Sup_Eq:metric}
{g}_{x}(\xi_{x}, \eta_{x}) 
 =  
\langle \xi_{\scriptsize \mat{U}_{1}},
{\eta}_{\scriptsize\mat{U}_{1}} (\mat{G}_{1} \mat{G}_{1}^T) \rangle +
\langle \xi_{\scriptsize \mat{U}_{2}},
{\eta}_{\scriptsize\mat{U}_{2}} (\mat{G}_{2} \mat{G}_{2}^T) \rangle +
\langle \xi_{\scriptsize \mat{U}_{3}},
{\eta}_{\scriptsize\mat{U}_{3}} (\mat{G}_{3} \mat{G}_{3}^T) \rangle 
+ \langle {\xi}_{\scriptsize \mathbfcal{G}}, {\eta}_{\scriptsize \mathbfcal{G}}\rangle ,
\end{eqnarray}
where ${\xi}_{{x}}, {\eta}_{{x}} \in T_{{x}} {\mathcal{M}}$ are 
tangent vectors with matrix characterizations 
$({\xi}_{\scriptsize \mat{U}_{1}}, {\xi}_{\scriptsize \mat{U}_{2}}
, {\xi}_{\scriptsize \mat{U}_{3}}, {\xi}_{\scriptsize \mathbfcal{G}})$ and
$({\eta}_{\scriptsize \mat{U}_{1}}, {\eta}_{\scriptsize \mat{U}_{2}}
, {\eta}_{\scriptsize \mat{U}_{3}}, {\eta}_{\scriptsize \mathbfcal{G}})$, respectively and $\langle \cdot, \cdot \rangle$ is the Euclidean inner product.
}

\subsection{\changeHK{Characterization of the normal space}}
\label{Sup_Sec:normalspace}
\changeHK{
Given a \changeBM{vector in $\mathbb{R}^{ n_1 \times r_1} \times \mathbb{R}^{ n_2 \times r_2} \times \mathbb{R}^{ n_3 \times r_3} \times \mathbb{R}^{r_1 \times r_2 \times r_3}$}, its projection onto the tangent space $T_{{x}} \mathcal{M}$ is obtained by extracting the component {\it normal}, in the metric sense, to the tangent space. This section describes the characterization of the {\it normal space}, $N_x \mathcal{M}$.
}

Let $\zeta_{x} = (\zeta_{\scriptsize \mat{U}_1}, \zeta_{\scriptsize \mat{U}_2}, \zeta_{\scriptsize \mat{U}_3}, \zeta_{\scriptsize \mathbfcal{G}}) \in N_{x} \mathcal{M}$, and $\eta_{{x}} = (\eta_{\scriptsize \mat{U}_1}, \eta_{\scriptsize \mat{U}_2}, \eta_{\scriptsize \mat{U}_3}, \eta_{\scriptsize \mathbfcal{G}}) \in T_{{x}} \mathcal{M}$. Since $\zeta_{{x}}$ is orthogonal to $\eta_{{x}}$, \changeHK{i.e., $g_{{x}}(\zeta_{{x}}, \eta_{{x}})=0$}, the conditions
\begin{eqnarray}
\label{Sup_Eq:TraceConditions}
{\rm Trace}( \mat{G}_{d}  \mat{G}_{d}^T  \zeta_{\scriptsize \mat{U}_d}^T  \eta_{\scriptsize \mat{U}_d}) = 0, \ {\rm for\ } d \in \{1,2,3 \}
\end{eqnarray}
must hold for all $\eta_{{x}}$ in the tangent space. Additionally from \cite[Example 3.5.2]{Absil_OptAlgMatManifold_2008}, $\eta_{\scriptsize \mat{U}_d}$ has the characterization
\begin{eqnarray}
\label{Sup_Eq:alternateTangentSpace}
\eta_{\scriptsize \mat{U}_d} & = & \mat{U}_d {\bf \Omega} + {\mat{U}_d}_{\perp} \mat{K},
\end{eqnarray}
where ${\bf \Omega}$ is any skew-symmetric matrix, $\mat{K}$ is a any matrix of size ${(n_d-r_d) \times r_d}$, and ${\mat{U}_d}_{\perp}$ is any $n_d \times (n_d - r_d)$ that is orthogonal complement of $\mat{U}_d$. Let $\tilde{\zeta}_{\scriptsize \mat{U}_d} =  \zeta_{\scriptsize \mat{U}_d} \mat{G}_{d}  \mat{G}_{d}^T $ and  \changeHK{let $\tilde{\zeta}_{\scriptsize \mat{U}_d}$} is defined as 
\begin{eqnarray}
\label{Sup_Eq:DefVariableChange}
\tilde{\zeta}_{\scriptsize \mat{U}_d} = \mat{U}_d \mat{A} + {\mat{U}_d}_{\perp} \mat{B}
\end{eqnarray}
without loss of generality, where $\mat{A} \in \mathbb{R}^{r_d \times r_d}$ and $\mat{B} \in \mathbb{R}^{(n_d-r_d) \times r_d}$ are to be characterized from (\ref{Sup_Eq:TraceConditions}) and (\ref{Sup_Eq:alternateTangentSpace}). A few standard computations show that $\mat{A}$ has to be symmetric and $\mat{B} = \mat{0}$. Consequently,
$\tilde{\zeta}_{\scriptsize \mat{U}_d} = \mat{U}_d \mat{S}_{\scriptsize \mat{U}_d}$, where $\mat{S}_{\scriptsize \mat{U}_d} = \mat{S}_{\scriptsize \mat{U}_d}^T$. Equivalently, $\zeta_{\scriptsize \mat{U}_d}  = \mat{U}_d \mat{S}_{\scriptsize \mat{U}_d} (\mat{G}_{d}  \mat{G}_{d}^T)^{-1}$ for a symmetric matrix $\mat{S}_{\scriptsize \mat{U}_d}$. Finally, the normal space $N_x \mathcal{M}$ has the characterization
\begin{equation}
\begin{array}{lll}
\label{Sup_Eq:NormalSpace}
N_x \mathcal{M}
 & = & \{(\mat{U}_{1}\mat{S}_{\scriptsize \mat{U}_{1}}(\mat{G}_{1}  \mat{G}_{1}^T)^{-1},
\mat{U}_{2}\mat{S}_{\scriptsize \mat{U}_{2}}(\mat{G}_{2}  \mat{G}_{2}^T)^{-1},
\mat{U}_{3}\mat{S}_{\scriptsize \mat{U}_{3}}(\mat{G}_{3}  \mat{G}_{3}^T)^{-1}, 0) : \\
&  & \mat{S}_{\scriptsize \mat{U}_{d}}  \in \mathbb{R}^{r_d \times r_d},  \mat{S}_{\scriptsize \mat{U}_{d}}^T  = \mat{S}_{\scriptsize \mat{U}_{d}}, \text{ for\ } d \in \{1, 2, 3\} \}.
\end{array}
\end{equation}

\subsection{\changeHK{Characterization of the vertical space}}
\label{Sup_Sec:verticalspace}

\changeHK{
The horizontal space projector of a tangent vector is obtained by removing the component along the vertical direction. This section shows the matrix characterization of the vertical space $\mathcal{V}_x$.
}

$\mathcal{V}_{x}$ is the defined as the linearization of the equivalence class $ [(\mat{U}_1, \mat{U}_2, \mat{U}_3, \mathbfcal{G})]$ at $x = [(\mat{U}_1, \mat{U}_2, \mat{U}_3, \mathbfcal{G})]$. Equivalently, $\mathcal{V}_{x}$ is the linearization of $(
 \mat{U}_1\mat{O}_{1}, \mat{U}_2\mat{O}_{2}, \mat{U}_3\mat{O}_{3}, \mathbfcal{G}_{\times 1} \mat{O}^T_{1}{_{\times 2}} \mat{O}^T_{2}{_{\times 3}} \mat{O}^T_{3})$ along $\mat{O}_{d} \in \mathcal{O}(r_d)$ at the \emph{identity element} for $d \in \{1, 2,3 \}$. From the characterization of linearization of an orthogonal matrix \cite[Example~3.5.3]{Absil_OptAlgMatManifold_2008}, we have the characterization for the vertical space as
\begin{equation}
\begin{array}{lll}
\label{Sup_Eq:VerticalComponents}
\mathcal{V}_{x} & = &\{ (\mat{U}_1{\bf \Omega}_1, \mat{U}_2 {\bf \Omega}_2, \mat{U}_3 {\bf \Omega}_3, 
 - (\mathbfcal{G}{\times_1} {\bf \Omega}_1  + 
\mathbfcal{G}{{\times_2}} {\bf \Omega}_2 +
\mathbfcal{G}{{\times_3}} {\bf \Omega}_3)):\\
&   & {\bf \Omega}_d \in \mathbb{R}^{r_d \times r_d}, {\bf \Omega}_d ^T = -{\bf \Omega}_d \text{ for\ } d \in \{1, 2, 3\} \}.
\end{array}
\end{equation}

\subsection{\changeHK{Characterization of the horizontal space}}
\label{Sup_Sec:horizontalspace}

\changeHK{
The characterization of the horizontal space $\mathcal{H}_{{x}}$ is derived from its orthogonal relationship with the vertical space $\mathcal{V}_{{x}}$.
}

Let $\xi_{{x}}=(\xi_{{\bf U}_1}, \xi_{{\bf U}_2}, \xi_{{\bf U}_3}, \xi_{\mathbfcal{G}})$  $\in \mathcal{H}_{{x}}$, and $\zeta_{{x}}=(\zeta_{{\bf U}_1}, \zeta_{{\bf U}_2}, \zeta_{{\bf U}_3}, \zeta_{\mathbfcal{G}})$ $\in \mathcal{V}_{{x}}$. Since $\xi_{{x}}$ must be orthogonal to $\zeta_{{x}}$, which \changeHK{is equivalent to} $g_{{x}}(\xi_{{x}}, \zeta_{{x}})=0$ \changeHK{in (\ref{Sup_Eq:metric})}, the characterization for $\xi_{{x}}$ is derived from (\ref{Sup_Eq:metric}) and  (\ref{Sup_Eq:VerticalComponents}).

\begin{eqnarray*}
{g}_{x}(\xi_{{x}}, \zeta_{{x}}) & = & 
\langle \xi_{\scriptsize \mat{U}_1},
\zeta_{\scriptsize \mat{U}_1} (\mat{G}_{1} \mat{G}_{1}^T) \rangle
+\langle \xi_{\scriptsize \mat{U}_2}, 
\zeta_{\scriptsize \mat{U}_2} (\mat{G}_{2} \mat{G}_{2}^T) \rangle 
+ \langle \xi_{\scriptsize \mat{U}_3}, 
\zeta_{\scriptsize \mat{U}_3} (\mat{G}_{3} \mat{G}_{3}^T) \rangle
+ \langle \xi_{\mathbfcal{G}},  \zeta_{\mathbfcal{G}} \rangle
 \\
& = & \langle {\xi}_{\scriptsize \mat{U}},
{\eta}_{\scriptsize \mat{U}_1} (\mat{G}_{1} \mat{G}_{1}^T) \rangle
+\langle {\xi}_{\scriptsize \mat{U}_2}, 
{\eta}_{\scriptsize \mat{U}_2} (\mat{G}_{2} \mat{G}_{2}^T) \rangle 
+ \langle {\xi}_{\scriptsize \mat{U}_3}, 
{\eta}_{\scriptsize \mat{U}_3} (\mat{G}_{3} \mat{G}_{3}^T) \rangle  \\
&  & \hspace{1.5cm}+ \langle \xi_{\mathbfcal{G}}, 
 - (\mathbfcal{G}{\times_1} {\bf \Omega}_1  + 
\mathbfcal{G}{{\times_2}} {\bf \Omega}_2 +
\mathbfcal{G}{{\times_3}} {\bf \Omega}_3) \rangle \\
 & = & 
\langle {\xi}_{\scriptsize \mat{U}_1},
{\eta}_{\scriptsize \mat{U}_1} (\mat{G}_{1} \mat{G}_{1}^T) \rangle
+\langle {\xi}_{\scriptsize  \mat{U}_2}, 
{\eta}_{\scriptsize \mat{U}_2} (\mat{G}_{2} \mat{G}_{2}^T) \rangle 
+ \langle {\xi}_{\scriptsize  \mat{U}_3}, 
{\eta}_{\scriptsize \mat{U}_3} (\mat{G}_{3} \mat{G}_{3}^T) \rangle  \\
&& \hspace{1.5cm} + \langle \xi_{\mathbfcal{G}}, - \mathbfcal{G}{\times_1} {\bf \Omega}_1  \rangle
+ \langle \xi_{\mathbfcal{G}}, - \mathbfcal{G}{\times_2} {\bf \Omega}_2 \rangle 
+ \langle \xi_{\mathbfcal{G}}, - \mathbfcal{G}{\times_3} {\bf \Omega}_3 \rangle \\
 & & \text{(We switch to unfoldings of\ } \mathbfcal{G}.) \\
 & = & 
{\rm Trace}( (\mat{G}_{1}  \mat{G}_{1}^T) \xi_{\scriptsize  \mat{U}_1}^T ( \mat{U}_1 {\bf \Omega}_1) ) 
+{\rm Trace}( (\mat{G}_{2}  \mat{G}_{2}^T) \xi_{\scriptsize  \mat{U}_2}^T ( \mat{U}_2 {\bf \Omega}_2) )  \\
&&+{\rm Trace}( (\mat{G}_{3}  \mat{G}_{3}^T) \xi_{\scriptsize  \mat{U}_3}^T ( \mat{U}_3 {\bf \Omega}_3) )  \\
&& + {\rm Trace} (\xi_{\scriptsize \mat{G}_{1}} (- {\bf \Omega}_1 \mat{G}_{1})^T ) 
+ {\rm Trace} (\xi_{\scriptsize \mat{G}_{2}} (- {\bf \Omega}_2 \mat{G}_{2})^T ) + {\rm Trace} (\xi_{\scriptsize \mat{G}_{3}} (- {\bf \Omega}_3 \mat{G}_{3})^T ) \\
 & = &
 {\rm Trace} \left[ \left\{(\mat{G}_{1}  \mat{G}_{1}^T) \xi_{\scriptsize  \mat{U}_1}^T \mat{U}_1 + \xi_{\scriptsize \mat{G}_{1}}  \mat{G}_{1}^T \right\} {\bf \Omega}_1 \right] \\
&& \hspace{0cm} +{\rm Trace} \left[ \left\{(\mat{G}_{2}  \mat{G}_{2}^T) \xi_{\scriptsize  \mat{U}_2}^T  \mat{U}_2  + \xi_{\scriptsize \mat{G}_{2}}  \mat{G}_{2}^T \right\} {\bf \Omega}_2 \right]    \\
&&  \hspace{0cm} +{\rm Trace} \left[ \left\{(\mat{G}_{3}  \mat{G}_{3}^T) \xi_{\scriptsize  \mat{U}_3}^T  \mat{U}_3 + \xi_{\scriptsize \mat{G}_{3}}  \mat{G}_{3}^T \right\} {\bf \Omega}_3 \right] , 
\end{eqnarray*}
where $\xi_{\scriptsize \mat{G}_{d}}$ is the mode-$d$ unfolding of $\xi_{\mathbfcal{G}}$.
Since ${g}_{{x}}(\xi_{{x}}, \zeta_{{x}})$ above should be zero for all skew-matrices ${\bf \Omega}_d$, $\xi_{{x}}=(\xi_{{\bf U}_1}, \xi_{{\bf U}_2}, \xi_{{\bf U}_3}, \xi_{\mathbfcal{G}})$  $\in \mathcal{H}_{{x}}$ must satisfy 
\begin{eqnarray}
\label{Sup_Eq:horizontal_space_reqrements}
(\mat{G}_{d}  \mat{G}_{d}^T) \xi_{\scriptsize \mat{U}_d}^T \mat{U}_d  + \xi_{\scriptsize \mat{G}_{d}}  \mat{G}_{d}^T\quad \text{is symmetric for } d \in  \{ 1, 2, 3 \} .
\end{eqnarray}

%
%
%
%




\subsection{\changeHK{Derivation of the tangent space projector}}
\label{Sup_Sec:tangentspaceprojector}

\changeHK{
The tangent space $T_{x} \mathcal{M}$ projector is obtained by extracting the component normal to $T_{x} \mathcal{M}$ in the ambient space. The normal space $N_{x} \mathcal{M}$ has the matrix characterization shown in (\ref{Sup_Eq:NormalSpace}).
The operator $\Psi_{{x}}: 
\mathbb{R}^{n_1 \times r_1} \times  
\mathbb{R}^{n_2 \times r_2} \times 
\mathbb{R}^{n_3 \times r_3} \times 
\mathbb{R}^{r_1 \times r_2 \times r_3} \rightarrow T_{{x}} {\mathcal{M}} :(\mat{Y}_{\scriptsize \mat{U}_{1}}, \mat{Y}_{\scriptsize \mat{U}_{2}}, \mat{Y}_{\scriptsize \mat{U}_{3}}, \mat{Y}_{\scriptsize \mathbfcal{G}} )$
$\mapsto \Psi_{{x}}(\mat{Y}_{\scriptsize \mat{U}_{1}}, \mat{Y}_{\scriptsize \mat{U}_{2}}, \mat{Y}_{\scriptsize \mat{U}_{3}}, \mat{Y}_{\scriptsize \mathbfcal{G}} )$ has the expression
\begin{equation}
\begin{array}{lll}
\label{Sup_Eq:Tangent_Projection}
\Psi_{{x}}(\mat{Y}_{\scriptsize \mat{U}_{1}}\!, \mat{Y}_{\scriptsize \mat{U}_{2}}\!, \mat{Y}_{\scriptsize \mat{U}_{3}}\!, \mat{Y}_{\scriptsize \mathbfcal{G}} )
=  (\mat{Y}_{\scriptsize \mat{U}_{1}}\!\! - \!\mat{U}_{1} \mat{S}_{\scriptsize \mat{U}_{1}} (\mat{G}_{1}  \mat{G}_{1}^T)^{-1},
\mat{Y}_{\scriptsize \mat{U}_{2}}\!\! - \!\mat{U}_{2} \mat{S}_{\scriptsize \mat{U}_{2}} (\mat{G}_{2}  \mat{G}_{2}^T)^{-1},  \\
\mat{Y}_{\scriptsize \mat{U}_{3}}\!\! -\! \mat{U}_{3}\mat{S}_{\scriptsize \mat{U}_{3}} (\mat{G}_{3}  \mat{G}_{3}^T)^{-1},
\mat{Y}_{\scriptsize \mathbfcal{G}}).
\end{array}
\end{equation}
}

From the definition of the tangent space in (\ref{Sup_Eq:tangent_space}), 
$\mat{U}_d$ should satisfy
\begin{eqnarray*}
 \eta_{{\bf U}_d}^T \mat{U}_d + \mat{U}_d^T \eta_{{\bf U}_d} & = & 
( \mat{Y}_{\scriptsize \mat{U}_d} - \mat{U}_d \mat{S}_{\scriptsize \mat{U}_d} (\mat{G}_{d}  \mat{G}_{d}^T)^{-1})^T \mat{U}_d
+ \mat{U}_d^T (\mat{Y}_{\scriptsize \mat{U}_d} - \mat{U}_d \mat{S}_{\scriptsize \mat{U}_d} (\mat{G}_{d}  \mat{G}_{d}^T)^{-1}) \\
& = & 
\mat{Y}_{\scriptsize \mat{U}_d}^T \mat{U}_d  - (\mat{G}_{d}  \mat{G}_{d}^T)^{-1} 
 \mat{S}_{\scriptsize \mat{U}_d}^T \mat{U}_d ^T \mat{U}_d
+ \mat{U}_d^T \mat{Y}_{\scriptsize \mat{U}_d} - \mat{U}_d^T \mat{U}_d \mat{S}_{\scriptsize \mat{U}_d} (\mat{G}_{d}  \mat{G}_{d}^T)^{-1}  \\
& = & 
\mat{Y}_{\scriptsize \mat{U}_d}^T \mat{U}_d  - (\mat{G}_{d}  \mat{G}_{d}^T)^{-1} 
 \mat{S}_{\scriptsize \mat{U}_d} 
+ \mat{U}_d^T \mat{Y}_{\scriptsize \mat{U}_d} - \mat{S}_{\scriptsize \mat{U}_d} (\mat{G}_{d}  \mat{G}_{d}^T)^{-1} \changeHK{\ =\ 0}.
\end{eqnarray*}
Multiplying $(\mat{G}_{d}  \mat{G}_{d}^T)$ from the right and left sides results in 
\begin{eqnarray*}
(\mat{G}_{d}  \mat{G}_{d}^T)^{-1} \mat{S}_{\scriptsize \mat{U}_d}  + \mat{S}_{\scriptsize \mat{U}_d} (\mat{G}_{d}  \mat{G}_{d}^T)^{-1}  
& = & \mat{Y}_{\scriptsize \mat{U}_d}^T \mat{U}_d + \mat{U}_d^T \mat{Y}_{\scriptsize \mat{U}_d} \\
\mat{S}_{\scriptsize \mat{U}_d} \mat{G}_{d}  \mat{G}_{d}^T + \mat{G}_{d}  \mat{G}_{d}^T \mat{S}_{\scriptsize \mat{U}_d}   
& = & \mat{G}_{d}  \mat{G}_{d}^T (\mat{Y}_{\scriptsize \mat{U}_d}^T \mat{U}_d + \mat{U}_d^T \mat{Y}_{\scriptsize \mat{U}_d} ) \mat{G}_{d}  \mat{G}_{d}^T.
\end{eqnarray*}

Finally, we obtain the \emph{Lyapunov} equation as
\begin{equation}
\begin{array}{lcl} \label{Sup_Eq:Req_horizontal_space}
\mat{S}_{\scriptsize \mat{U}_{d}} \mat{G}_{d}  \mat{G}_{d}^T + \mat{G}_{d}  \mat{G}_{d}^T \mat{S}_{\scriptsize \mat{U}_{d}}   
=\mat{G}_{d}  \mat{G}_{d}^T (\mat{Y}_{\scriptsize \mat{U}_{d}}^T \mat{U}_{d} + \mat{U}_{d}^T \mat{Y}_{\scriptsize \mat{U}_{d}}) \mat{G}_{d}  \mat{G}_{d}^T\ \ {\rm for\ } d\in \{ 1,2,3\},
\end{array}
\end{equation}
that are solved efficiently with the Matlab's \verb+lyap+ routine.

%
%

\subsection{\changeHK{Derivation of the horizontal space projector}}
\label{Sup_Sec:horizontalspaceprojector}

We consider the projection of a tangent vector $\eta_{{x}}=(\eta_{\scriptsize \mat{U}_1}, \eta_{\scriptsize \mat{U}_2}, \eta_{\scriptsize \mat{U}_3}, \eta_{\scriptsize \mathbfcal{G}}) \in T_{{x}} \mathcal{M}$ into a vector $\xi_{{x}}=(\xi_{{\bf U}_1}, \xi_{{\bf U}_2}, \xi_{{\bf U}_3}, \xi_{\mathbfcal{G}}) \in H_{{x}}$. This is achieved by subtracting the component in the vertical space $\mathcal{V}_{{x}}$ in (\ref{Sup_Eq:VerticalComponents}) as
\begin{eqnarray*}
\left\{
\begin{array}{lll}
\eta_{\scriptsize \mat{U}_1} & = &  \underbrace{\eta_{\scriptsize \mat{U}_1} - \mat{U}_1 {\bf \Omega}_1}
_{=\xi_{\scriptsize \mat{U}_1} \in \mathcal{H}_{{x}}}
+ \underbrace{\mat{U}_1 {\bf \Omega}_1}_{ \in \mathcal{V}_{{x}}},  \\
\eta_{\scriptsize \mat{U}_2} & = &  \eta_{\scriptsize \mat{U}_2} - \mat{U}_2 {\bf \Omega}_2 + \mat{U}_2 {\bf \Omega}_2, \\
\eta_{\scriptsize \mat{U}_3} & = &  \eta_{\scriptsize \mat{U}_3} - \mat{U}_3 {\bf \Omega}_3 + \mat{U}_3 {\bf \Omega}_3, \\
\eta_{\scriptsize \mathbfcal{G}} & = &  \eta_{\scriptsize \mathbfcal{G}} - ( - (\mathbfcal{G}{\times_1} {\bf \Omega}_1  + 
\mathbfcal{G}{{\times_2}} {\bf \Omega}_2 +
\mathbfcal{G}{{\times_3}} {\bf \Omega}_3)) + ( - (\mathbfcal{G}{\times_1} {\bf \Omega}_1  + 
\mathbfcal{G}{{\times_2}} {\bf \Omega}_2 +
\mathbfcal{G}{{\times_3}} {\bf \Omega}_3)).
\end{array}
\right.
\end{eqnarray*}

As a result, the horizontal operator $\Pi_{{x}}: T_x \mathcal{M} \rightarrow \mathcal{H}_x: \eta_x \mapsto \Pi_{{x}}(\eta_x)$ has the expression
\begin{equation}
\begin{array}{lll}
\label{Sup_Eq:Horizontal_Projection}
\Pi_{{x}}(\eta_{{x}} )
& = &(
\eta_{\scriptsize \mat{U}_{1}} - \mat{U}_{1} {\bf \Omega}_{1},
\eta_{\scriptsize \mat{U}_{2}} - \mat{U}_{2} {\bf \Omega}_{2}, 
\eta_{\scriptsize \mat{U}_{3}} - \mat{U}_{3} {\bf \Omega}_{3}, \\
& &\eta_{\scriptsize \mathbfcal{G}} \!-\! (\! - (\mathbfcal{G}{\times_1} {\bf \Omega}_{1} + 
\mathbfcal{G}{{\times_2}} {\bf \Omega}_{2} +
\mathbfcal{G}{{\times_3}} {\bf \Omega}_{3})) 
),
\end{array}
\end{equation}
where $\eta_x = (\eta_{\scriptsize \mat{U}_{1}}, \eta_{\scriptsize \mat{U}_{2}}, \eta_{\scriptsize \mat{U}_{3}}, \eta_{\scriptsize \mathbfcal{G}}) \in T_x \mathcal{M}$ and  ${\bf \Omega}_{d}$ is a skew-symmetric matrix of size $r_d \times r_d$. The skew-matrices ${\bf \Omega}_{d}$ for $d = \{1, 2,3\}$ that are identified based on the conditions (\ref{Sup_Eq:horizontal_space_reqrements}).

It should be noted that the tensor $\mathbfcal{G}{\times_1} {\bf \Omega}_1  + \mathbfcal{G}{{\times_2}} {\bf \Omega}_2 + \mathbfcal{G}{{\times_3}} {\bf \Omega}_3$ in (\ref{Sup_Eq:VerticalComponents}) has the following equivalent unfoldings.
\begin{eqnarray*}
\begin{array}{lll}
\label{Sup_Eq:VerticalComponentsMatrixRep}
\mathbfcal{G}{\times_1} {\bf \Omega}_1  + 
\mathbfcal{G}{{\times_2}} {\bf \Omega}_2 +
\mathbfcal{G}{{\times_3}} {\bf \Omega}_3 & \xLeftrightarrow{\text{mode\ }-1} & 
{\bf \Omega}_{1} \mat{G}_{1} + \mat{G}_{1}(\mat{I}_{r_3} \otimes {\bf \Omega}_2)^T
+ \mat{G}_{1}( {\bf \Omega}_3 \otimes \mat{I}_{r_2} )^T \\
& \xLeftrightarrow{\text{mode\ }-2}  & 
\mat{G}_{2}(\mat{I}_{r_3} \otimes {\bf \Omega}_1)^T + {\bf \Omega}_{2} \mat{G}_{2}  
+ \mat{G}_{2}( {\bf \Omega}_3 \otimes \mat{I}_{r_1} )^T \\
& \xLeftrightarrow{\text{mode\ }-3}  & 
\mat{G}_{3}(\mat{I}_{r_2} \otimes {\bf \Omega}_1)^T 
+ \mat{G}_{3}( {\bf \Omega}_2 \otimes \mat{I}_{r_1} )^T + {\bf \Omega}_{3} \mat{G}_{3}.
\end{array}
\end{eqnarray*}

Plugging $\xi_{\scriptsize \mat{U}_1} = \eta_{\scriptsize \mat{U}_{1}} - \mat{U}_{1} {\bf \Omega}_{1}$ and $\xi_{\scriptsize \mat{G}_1} = \changeHK{\eta{\scriptsize \mat{G}_{1}} +}{\bf \Omega}_{1} \mat{G}_{1} + \mat{G}_{1}(\mat{I}_{r_3} \otimes {\bf \Omega}_2)^T
+ \mat{G}_{1}( {\bf \Omega}_3 \otimes \mat{I}_{r_2} )^T$ into (\ref{Sup_Eq:horizontal_space_reqrements}) and using the relation $(\mat{A} \otimes \mat{B})^T = \mat{A}^T \otimes \mat{B}^T$ results in
\begin{eqnarray}
\begin{array}{l}
\label{Sup_Eq:ConditionA}
(\mat{G}_{1}  \mat{G}_{1}^T) \xi_{\scriptsize \mat{U}_1}^T \mat{U}_  + \xi_{\scriptsize \mat{G}_{1}}  \mat{G}_{1}^T    \\
\hspace{1cm}= (\mat{G}_{1}  \mat{G}_{1}^T) (\eta_{\scriptsize \mat{U}_1} - \mat{U}_1 {\bf \Omega}_1)^T \mat{U}_1  \\
\hspace{3cm} + \left\{ \eta{\scriptsize \mat{G}_{1}} +({\bf \Omega}_{1} \mat{G}_{1} + \mat{G}_{1}(\mat{I}_{r_3} \otimes {\bf \Omega}_2)^T
+ \mat{G}_{1}( {\bf \Omega}_3 \otimes \mat{I}_{r_2} )^T)\right\} \mat{G}_{1}^T  \\
\hspace{1cm}= (\mat{G}_{1}  \mat{G}_{1}^T) \eta_{\scriptsize \mat{U}_1} ^T \mat{U}_1 
- (\mat{G}_{1}  \mat{G}_{1}^T) ( \mat{U}_1 {\bf \Omega}_1)^T \mat{U}_1  \\
\hspace{3cm} + \eta_{\scriptsize \mat{G}_{1}}\mat{G}_{1}^T  
+ {\bf \Omega}_{1} \mat{G}_{1}\mat{G}_{1}^T 
+\mat{G}_{1}(\mat{I}_{r_3} \otimes {\bf \Omega}_2)^T \mat{G}_{1}^T 
+ \mat{G}_{1}( {\bf \Omega}_3 \otimes \mat{I}_{r_2} )^T\mat{G}_{1}^T\\
\hspace{1cm}= (\mat{G}_{1}  \mat{G}_{1}^T) \eta_{\scriptsize \mat{U}_1} ^T \mat{U}_1 
+ (\mat{G}_{1}  \mat{G}_{1}^T) {\bf \Omega}_1  \\
\hspace{3cm} + \eta_{\scriptsize \mat{G}_{1}}\mat{G}_{1}^T  
+  {\bf \Omega}_{1} \mat{G}_{1}\mat{G}_{1}^T
- \mat{G}_{1}(\mat{I}_{r_3} \otimes {\bf \Omega}_2) \mat{G}_{1}^T 
- \mat{G}_{1}( {\bf \Omega}_3 \otimes \mat{I}_{r_2} ) \mat{G}_{1}^T,
\end{array}
\end{eqnarray}
which should be a symmetric matrix \changeHK{due to (\ref{Sup_Eq:horizontal_space_reqrements})}, i.e., $
(\mat{G}_{1}  \mat{G}_{1}^T) \xi_{\scriptsize \mat{U}_1}^T \mat{U}_  + \xi_{\scriptsize \mat{G}_{1}}  \mat{G}_{1}^T =  
((\mat{G}_{1}  \mat{G}_{1}^T) \xi_{\scriptsize \mat{U}_1}^T \mat{U}_  + \xi_{\scriptsize \mat{G}_{1}}  \mat{G}_{1}^T)^T$.


Subsequently,  
\begin{eqnarray*}
(\mat{G}_{1}  \mat{G}_{1}^T) \eta_{\scriptsize \mat{U}_1} ^T \mat{U}_1 
+ (\mat{G}_{1}  \mat{G}_{1}^T) {\bf \Omega}_1 
+ \eta_{\scriptsize \mat{G}_{1}}\mat{G}_{1}^T  
+ {\bf \Omega}_{1} \mat{G}_{1}\mat{G}_{1}^T 
- \mat{G}_{1}(\mat{I}_{r_3} \otimes {\bf \Omega}_2) \mat{G}_{1}^T 
- \mat{G}_{1}( {\bf \Omega}_3 \otimes \mat{I}_{r_2} ) \mat{G}_{1}^T \\
= \mat{U}_1^T \eta_{\scriptsize \mat{U}_1} (\mat{G}_{1}  \mat{G}_{1}^T) 
- {\bf \Omega}_1 \mat{G}_{1}  \mat{G}_{1}^T  
+ \mat{G}_{1} \eta_{\scriptsize \mat{G}_{1}}^T  
-  \mat{G}_{1}\mat{G}_{1}^T {\bf \Omega}_{1}
+\mat{G}_{1}(\mat{I}_{r_3} \otimes {\bf \Omega}_2) \mat{G}_{1}^T 
+ \mat{G}_{1}( {\bf \Omega}_3 \otimes \mat{I}_{r_2} ) \mat{G}_{1}^T,
\end{eqnarray*}
which is equivalent to
\begin{eqnarray*}
\mat{G}_{1}  \mat{G}_{1}^T {\bf \Omega}_1 
+ {\bf \Omega}_{1} \mat{G}_{1}\mat{G}_{1}^T 
- \mat{G}_{1}(\mat{I}_{r_3} \otimes {\bf \Omega}_2) \mat{G}_{1}^T 
- \mat{G}_{1}( {\bf \Omega}_3 \otimes \mat{I}_{r_2} ) \mat{G}_{1}^T \\ 
\hspace{2cm}= {\rm Skew}(\mat{U}_1^T \eta_{\scriptsize \mat{U}_1} \mat{G}_{1}  \mat{G}_{1}^T)
+ {\rm Skew}(\mat{G}_{1} \eta_{\scriptsize \mat{G}_{1}}^T ).
\end{eqnarray*}
Here ${\rm Skew}(\cdot)$ extracts the skew-symmetric part of a square matrix, i.e., ${\rm Skew}(\mat{D})=(\mat{D}-\mat{D}^T)/2$.

Finally, we obtain the \emph{coupled} Lyapunov equations 
\begin{equation}
\begin{array}{lll}
\label{Sup_Eq:OmegaRequirements}
\left\{
\begin{array}{l}
\mat{G}_{1}  \mat{G}_{1}^T {\bf \Omega}_{1} + {\bf \Omega}_{1} \mat{G}_{1}  \mat{G}_{1}^T 
-\mat{G}_{1}(\mat{I}_{r_3} \otimes {\bf \Omega}_{2}) \mat{G}_{1}^T 
- \mat{G}_{1}( {\bf \Omega}_{3} \otimes \mat{I}_{r_2} )\mat{G}_{1}^T  \\
\hspace{6cm} = {\rm Skew}(\mat{U}_1^T\eta_{\scriptsize \mat{U}_1}\mat{G}_{1}  \mat{G}_{1}^T) + {\rm Skew}(\mat{G}_{1}\eta_{\scriptsize \mat{G}_{1}}^T), \\
\mat{G}_{2}  \mat{G}_{2}^T {\bf \Omega}_{2} + {\bf \Omega}_{2} \mat{G}_{2}  \mat{G}_{2}^T 
-\mat{G}_{2}(\mat{I}_{r_3} \otimes {\bf \Omega}_{1}) \mat{G}_{2}^T 
- \mat{G}_{2}( {\bf \Omega}_{3} \otimes \mat{I}_{r_1} )\mat{G}_{2}^T  \\
\hspace{6cm} = {\rm Skew}(\mat{U}_2^T\eta_{\scriptsize \mat{U}_2}\mat{G}_{2}  \mat{G}_{2}^T) + {\rm Skew}(\mat{G}_{2}\eta_{\scriptsize \mat{G}_{2}}^T), \\
\mat{G}_{3}  \mat{G}_{3}^T {\bf \Omega}_{3} + {\bf \Omega}_{3} \mat{G}_{3}  \mat{G}_{3}^T 
-\mat{G}_{3}(\mat{I}_{r_2} \otimes {\bf \Omega}_{1}) \mat{G}_{3}^T 
- \mat{G}_{3}( {\bf \Omega}_{2} \otimes \mat{I}_{r_1} )\mat{G}_{3}^T  \\
\hspace{6cm} = {\rm Skew}(\mat{U}_3^T\eta_{\scriptsize \mat{U}_3}\mat{G}_{3}  \mat{G}_{3}^T) + {\rm Skew}(\mat{G}_{3}\eta_{\scriptsize \mat{G}_{3}}^T),
\end{array}
\right.
\end{array}
\end{equation}
that are solved efficiently with the Matlab's \verb+pcg+ routine that is combined with a specific preconditioner resulting from the Gauss-Seidel approximation of (\ref{Sup_Eq:OmegaRequirements}).

\subsection{\changeHK{Derivation of the Riemannian gradient formula}}
\label{Sup_Sec:Riemanniangradientformula}

\changeBM{Let $\changeBM{f(\mathbfcal{X})}=\| \mathbfcal{P}_{\Omega}(\mathbfcal{X}) - \mathbfcal{P}_{\Omega}(\mathbfcal{X}^{\star}) \|^2_F/|\Omega |$ and 
$\mathbfcal{S} = 2 (\mathbfcal{P}_{\Omega}(\mathbfcal{G}{\times_1} \mat{U}_{1} {\times_2} \mat{U}_{2}{\times_3} \mat{U}_{3}) - 
\mathbfcal{P}_{\Omega}(\mathbfcal{X}^{\star}))/|{\Omega}|$ 
be an auxiliary sparse tensor variable that is interpreted as the Euclidean gradient of $f$ in $\mathbb{R}^{n_1 \times n_2 \times n_3}$. }

The partial derivatives of $f(\mat{U}_1, \mat{U}_2, \mat{U}_3, \mathbfcal{G})$ are
\begin{eqnarray*}
\label{Sup_Eq:Egradient}
 \left\{
\begin{array}{lll}
\displaystyle{\frac{\partial f_{1}(\mat{U}_1, \mat{U}_2, \mat{U}_3, \mat{G}_{1})}{\partial \mat{U}_1}}
&  = & \displaystyle{\frac{2}{|{\Omega} |} ( {\mathbfcal{P}}_{\Omega} (\mat{U}_1 \mat{G}_{1} (\mat{U}_3 \otimes \mat{U}_2)^T) -  {\mathbfcal{P}}_{\Omega}(\mat{X}^{\star}_{1}) ) (\mat{U}_3 \otimes \mat{U}_2) \mat{G}_{1}^T} \\ 
&  = & \mat{S}_{1} (\mat{U}_3 \otimes \mat{U}_2) \mat{G}_{1}^T, \\
\displaystyle{\frac{\partial f_{2}(\mat{U}_1, \mat{U}_2, \mat{U}_3, \mat{G}_{2})}{\partial \mat{U}_2}}
&  = & \displaystyle{\frac{2}{|{\Omega} |} ( {\mathbfcal{P}}_{\Omega} (\mat{U}_2 \mat{G}_{2} (\mat{U}_3 \otimes \mat{U}_1)^T) - {\mathbfcal{P}}_{\Omega} (\mat{X}^{\star}_{2}) ) (\mat{U}_3 \otimes \mat{U}_1) \mat{G}_{2}^T} \\ 
&  = & \mat{S}_{2} (\mat{U}_2 \otimes \mat{U}_1) \mat{G}_{2}^T, \\
\displaystyle{\frac{\partial f_{3}(\mat{U}_1, \mat{U}_2, \mat{U}_3, \mat{G}_{3})}{\partial \mat{U}_3}}
&  = & \displaystyle{\frac{2}{|{\Omega} |} ( {\mathbfcal{P}}_{\Omega} (\mat{U}_3 \mat{G}_{3} (\mat{U}_2 \otimes \mat{U}_1)^T) - {\mathbfcal{P}}_{\Omega} (\mat{X}^{\star}_{3}) ) (\mat{U}_2 \otimes \mat{U}_1) \mat{G}_{3}^T} \\ 
&  = & \mat{S}_{3} (\mat{U}_2 \otimes \mat{U}_1) \mat{G}_{3}^T, \\
\displaystyle{\frac{\partial f(\mat{U}_1, \mat{U}_2, \mat{U}_3, \mathbfcal{G})}{\partial \mathbfcal{G}}}
&  = & \displaystyle{\frac{2}{|{\Omega} |} (\mathbfcal{P}_{\Omega}(\mathbfcal{G}_{\times 1}\mat{U}_1 {\times_2} \mat{U}_2 {\times_3} \mat{U}_3) - 
\mathbfcal{P}_{\Omega}(\mathbfcal{X}^{\star}))} \\ 
& & \hspace{4cm} \times_1 \mat{U}_1^T \times_2 \mat{U}_2^T \times_3 \mat{U}_3^T\\
& = & \mathbfcal{S} \times_1 \mat{U}_1^T \times_2 \mat{U}_2T \times_3 \mat{U}_3^T,
\end{array}
\right.
\end{eqnarray*}
where $\mat{X}^{\star}_{d}$ is mode-$d$ unfolding of $\mathbfcal{X}^{\star}$ and 
\begin{eqnarray*}
 \left\{
\begin{array}{lll}
\mat{S}_{1} & =  & \displaystyle{\frac{2}{|{\Omega} |} ( {\mathbfcal{P}}_{\Omega}(\mat{U}_1 \mat{G}_{1} (\mat{U}_3 \otimes \mat{U}_2)^T) - {\mathbfcal{P}}_{\Omega}(\mat{X}^{\star}_{1}) )}\\ 
\mat{S}_{2} & =  & \displaystyle{\frac{2}{|{\Omega} |} ( {\mathbfcal{P}}_{\Omega}(\mat{U}_2 \mat{G}_{2} (\mat{U}_3 \otimes \mat{U}_1)^T) - {\mathbfcal{P}}_{\Omega}(\mat{X}^{\star}_{2}) )}\\ 
\mat{S}_{3} & =  & \displaystyle{\frac{2}{|{\Omega} |} ( {\mathbfcal{P}}_{\Omega}(\mat{U}_3 \mat{G}_{3} (\mat{U}_2 \otimes \mat{U}_1)^T) - {\mathbfcal{P}}_{\Omega}(\mat{X}^{\star}_{3}) )}\\ 
\mathbfcal{S} & =  & \displaystyle{\frac{2}{|{\Omega} |} (\mathbfcal{P}_{\Omega}(\mathbfcal{G} {\times_1} \mat{U}_1 {\times_2} \mat{U}_2 {\times_3} \mat{U}_3) - 
\mathbfcal{P}_{\Omega}(\mathbfcal{X}^{\star}))}.
\end{array}
\right.
\end{eqnarray*}


Due to the specific scaled metric (\ref{Sup_Eq:metric}), the partial derivatives of $f$ are further scaled by $((\mat{G}_{1}\mat{G}_{1}^T)^{-1}, (\mat{G}_{2}\mat{G}_{2}^T)^{-1}, (\mat{G}_{3}\mat{G}_{3}^T)^{-1}, \mathbfcal{I})$, denoted as ${\rm egrad}_{x} f$ (after scaling), i.e.,

\[
\begin{array}{lcl}
{\rm egrad}_{x} f & = &( \mat{S}_{1} (\mat{U}_3 \otimes \mat{U}_2) \mat{G}_{1}^T (\mat{G}_{1}\mat{G}_{1}^T)^{-1}, 
\mat{S}_{2} (\mat{U}_3 \otimes \mat{U}_1) \mat{G}_{2}^T (\mat{G}_{2}\mat{G}_{2}^T)^{-1},  \\
&  & \hspace{0cm} \mat{S}_{3} (\mat{U}_2 \otimes \mat{U}_1) \mat{G}_{3}^T(\mat{G}_{3}\mat{G}_{3}^T)^{-1},  \mathbfcal{S} \times_1 \mat{U}_1^T \times_2 \mat{U}_2^T \times_3 \mat{U}_3^T). \\
\end{array}
\]

Consequently, from the relationship that horizontal lift of ${\rm grad}_{[x]}f$ is equal to ${\rm grad}_{x}f \ =\ \Psi({\rm egrad}_{x} f )$, we obtain that, using (\ref{Sup_Eq:Tangent_Projection}),
\begin{eqnarray*}
\label{Sup_Eq:RiemannianGradientDerive}
\text{the horizontal lift of\ }{\rm grad}_{[x]} f & = &  
\Psi(
\mat{S}_{1} (\mat{U}_3 \otimes \mat{U}_2) \mat{G}_{1}^T (\mat{G}_{1}\mat{G}_{1}^T)^{-1}, 
\mat{S}_{2} (\mat{U}_3 \otimes \mat{U}_1) \mat{G}_{2}^T (\mat{G}_{2}\mat{G}_{2}^T)^{-1},  \\
&  & \hspace{0cm} \mat{S}_{3} (\mat{U}_2 \otimes \mat{U}_1) \mat{G}_{3}^T(\mat{G}_{3}\mat{G}_{3}^T)^{-1},  \mathbfcal{S} \times_1 \mat{U}_1^T \times_2 \mat{U}_2^T \times_3 \mat{U}_3^T)  \\
& = & 
(\mat{S}_{1} (\mat{U}_3 \otimes \mat{U}_2) \mat{G}_{1}^T (\mat{G}_{1}\mat{G}_{1}^T)^{-1} 
- \mat{U}_1\mat{B}_{\scriptsize \mat{U}_1}(\mat{G}_{1}\mat{G}_{1}^T)^{-1}, \\
&  & \hspace{0cm} \mat{S}_{2} (\mat{U}_3\otimes \mat{U}_1) \mat{G}_{2}^T (\mat{G}_{2}\mat{G}_{2}^T)^{-1}
- \mat{U}_2 \mat{B}_{\scriptsize \mat{U}_2}(\mat{G}_{2}\mat{G}_{2}^T)^{-1}, \\
&  & \hspace{0cm} \mat{S}_{3} (\mat{U}_2 \otimes \mat{U}_1) \mat{G}_{3}^T(\mat{G}_{3}\mat{G}_{3}^T)^{-1} 
- \mat{U}_3\mat{B}_{\scriptsize \mat{U}_3}(\mat{G}_{3}\mat{G}_{3}^T)^{-1}, \\
& & \mathbfcal{S} \times_1 \mat{U}_1^T \times_2 \mat{U}_2^T \times_3 \mat{U}_3^T). 
\end{eqnarray*}


From the requirements in (\ref{Sup_Eq:Req_horizontal_space}) for a vector to be in the tangent space, we have the following relationship for mode-$1$.
\begin{eqnarray*}
\mat{B}_{\scriptsize \mat{U}_1} \mat{G}_{1}  \mat{G}_{1}^T + \mat{G}_{1}  \mat{G}_{1}^T \mat{B}_{\scriptsize \mat{U}_1}   
& = & \mat{G}_{1}  \mat{G}_{1}^T (\mat{Y}_{\scriptsize \mat{U}_1}^T \mat{U}_1 + \mat{U}_1^T \mat{Y}_{\scriptsize \mat{U}_1} ) \mat{G}_{1}  \mat{G}_{1}^T,
\end{eqnarray*}
where $\mat{Y}_{\scriptsize \mat{U}_1}  = (\mat{S}_{1} (\mat{U}_3 \otimes \mat{U}_2) \mat{G}_{1}^T (\mat{G}_{1}\mat{G}_{1}^T)^{-1}$.

Subsequently,
\[
\begin{array}{lll}
\mat{G}_{1}  \mat{G}_{1}^T (\mat{Y}_{\scriptsize \mat{U}_1}^T \mat{U}_1 + \mat{U}_1^T \mat{Y}_{\scriptsize \mat{U}_1} ) \mat{G}_{1}  \mat{G}_{1}^T 
& = & 
\mat{G}_{1}  \mat{G}_{1}^T \left\{ 
((\mat{S}_{1} (\mat{U}_3 \otimes \mat{U}_2) \mat{G}_{1}^T (\mat{G}_{1}\mat{G}_{1}^T)^{-1})^T  \mat{U}_1 \right.  \\ 
& & \left. + \mat{U}_1^T (\mat{S}_{1} (\mat{U}_3 \otimes \mat{U}_2) \mat{G}_{1}^T (\mat{G}_{1}\mat{G}_{1}^T)^{-1}  \right\} \mat{G}_{1}  \mat{G}_{1}^T \\ 
& = & ((\mat{S}_{1} (\mat{U}_3 \otimes \mat{U}_2) \mat{G}_{1}^T)^T \mat{U}_1 \mat{G}_{1}  \mat{G}_{1}^T+ \mat{G}_{1}  \mat{G}_{1}^T\mat{U}^T (\mat{S}_{1} (\mat{U}_3\otimes \mat{U}_2) \mat{G}_{1}^T \\ 
& = & (\mat{G}_{1}  \mat{G}_{1}^T\mat{U}_1^T (\mat{S}_{1} (\mat{U}_3 \otimes \mat{U}_2) \mat{G}_{1}^T)^T+ \mat{G}_{1}  \mat{G}_{1}^T\mat{U}_1^T (\mat{S}_{1} (\mat{U}_3 \otimes \mat{U}_2) \mat{G}_{1}^T \\ 
& = & 2 {\rm Sym} (\mat{G}_{1}  \mat{G}_{1}^T\mat{U}_1^T (\mat{S}_{1} (\mat{U}_3 \otimes \mat{U}_2) \mat{G}_{1}^T).
\end{array}
\]

Finally, $\mat{B}_{\scriptsize \mat{U}_{d}}$ for $d \in \{1, 2,3\} $ are obtained by solving the Lyapunov equations
\begin{eqnarray*}
\label{BUBVBWRequirementGradient}
\left\{
\begin{array}{lll}
\mat{B}_{\scriptsize \mat{U}_{1}} \mat{G}_{1}  \mat{G}_{1}^T + \mat{G}_{1}  \mat{G}_{1}^T \mat{B}_{\scriptsize \mat{U}_{1}}
& = & 2 {\rm Sym} (\mat{G}_{1}  \mat{G}_{1}^T\mat{U}_{1}^T (\mat{S}_{1} (\mat{U}_{3} \otimes \mat{U}_{2}) \mat{G}_{1}^T), \\ 
\mat{B}_{\scriptsize \mat{U}_{2}} \mat{G}_{2}  \mat{G}_{2}^T + \mat{G}_{2}  \mat{G}_{2}^T \mat{B}_{\scriptsize \mat{U}_{2}}
& = & 2 {\rm Sym} (\mat{G}_{2}  \mat{G}_{2}^T\mat{U}_{2}^T (\mat{S}_{2} (\mat{U}_{3} \otimes \mat{U}_{1}) \mat{G}_{2}^T), \\ 
\mat{B}_{\scriptsize \mat{U}_{3}} \mat{G}_{3}  \mat{G}_{3}^T + \mat{G}_{3}  \mat{G}_{3}^T \mat{B}_{\scriptsize \mat{U}_{3}}
& = & 2 {\rm Sym} (\mat{G}_{3}  \mat{G}_{3}^T\mat{U}_{3}^T (\mat{S}_{3} (\mat{U}_{2} \otimes \mat{U}_{1}) \mat{G}_{3}^T).
\end{array}
\right.
\end{eqnarray*}
where ${\rm Sym}(\cdot)$ extracts the symmetric part of a square matrix, i.e., ${\rm Sym}(\mat{D})=(\mat{D}+\mat{D}^T)/2$. The above Lyapunov equations are solved efficiently with the Matlab's \verb+lyap+ routine.

\section{Additional numerical comparisons}
\label{sec:AdditionalNumericalComparisons}

In addition to the representative numerical comparisons in the paper, we show additional numerical experiments spanning synthetic and real-world datasets. 

{\bf Experiments on synthetic datasets:}

\begin{figure}[htbp]
\begin{tabular}{cc}
\begin{minipage}{0.32\hsize}
\begin{center}
\includegraphics[width=\hsize]{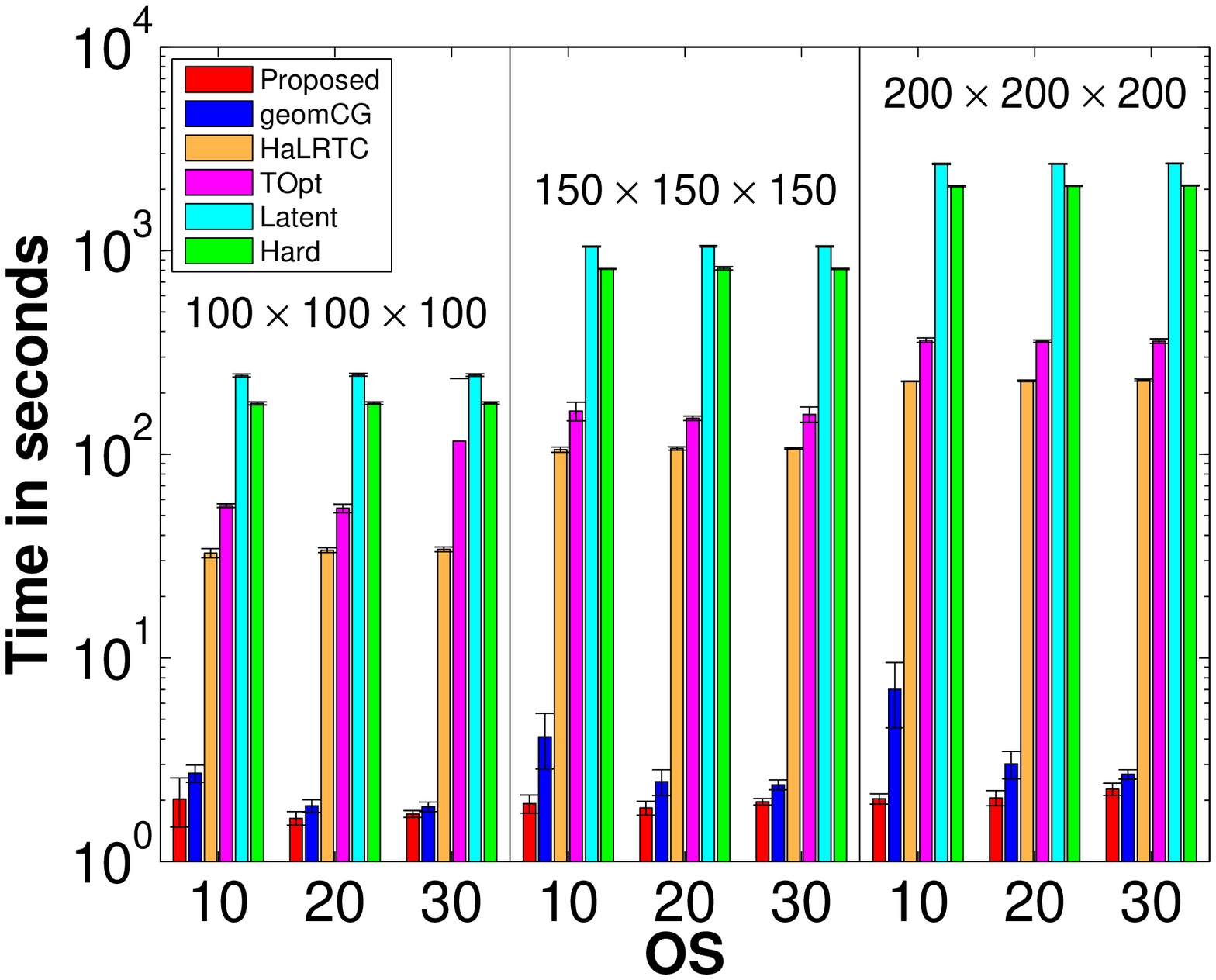}\\
{\scriptsize(a) {\bf r} = ($5,5,5$).}
\end{center}
\end{minipage}
\begin{minipage}{0.32\hsize}
\begin{center}
\includegraphics[width=\hsize]{figures/caseS2_core_10_time_bar.eps}\\
{\scriptsize(b) {\bf r} = ($10,10,10$).}
\end{center}
\end{minipage}
\begin{minipage}{0.32\hsize}
\begin{center}
\includegraphics[width=\hsize]{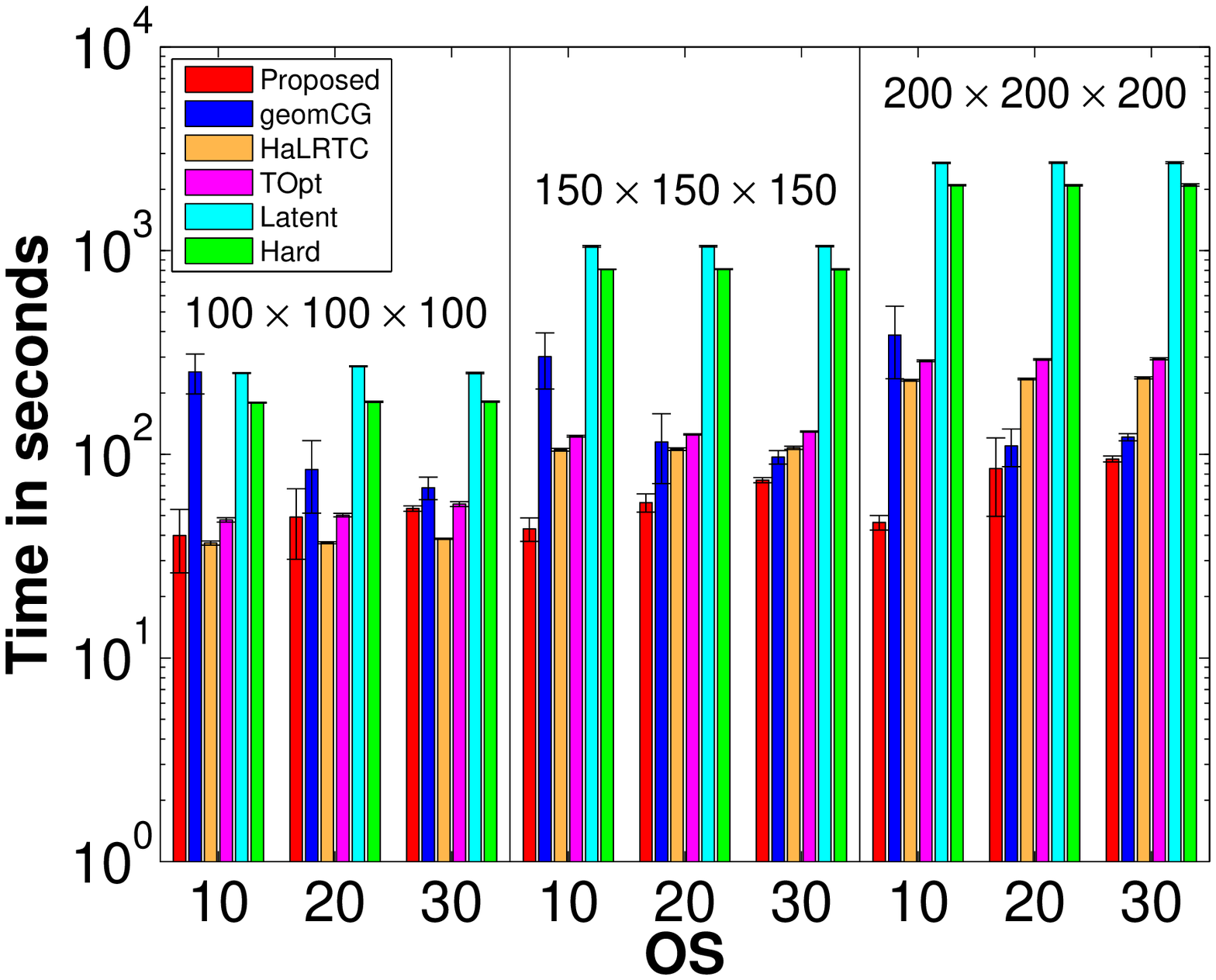}\\
{\scriptsize(c) {\bf r} = ($15,15,15$).}
\end{center}
\end{minipage}\\
\begin{minipage}{0.32\hsize}
\begin{center}
\includegraphics[width=\hsize]{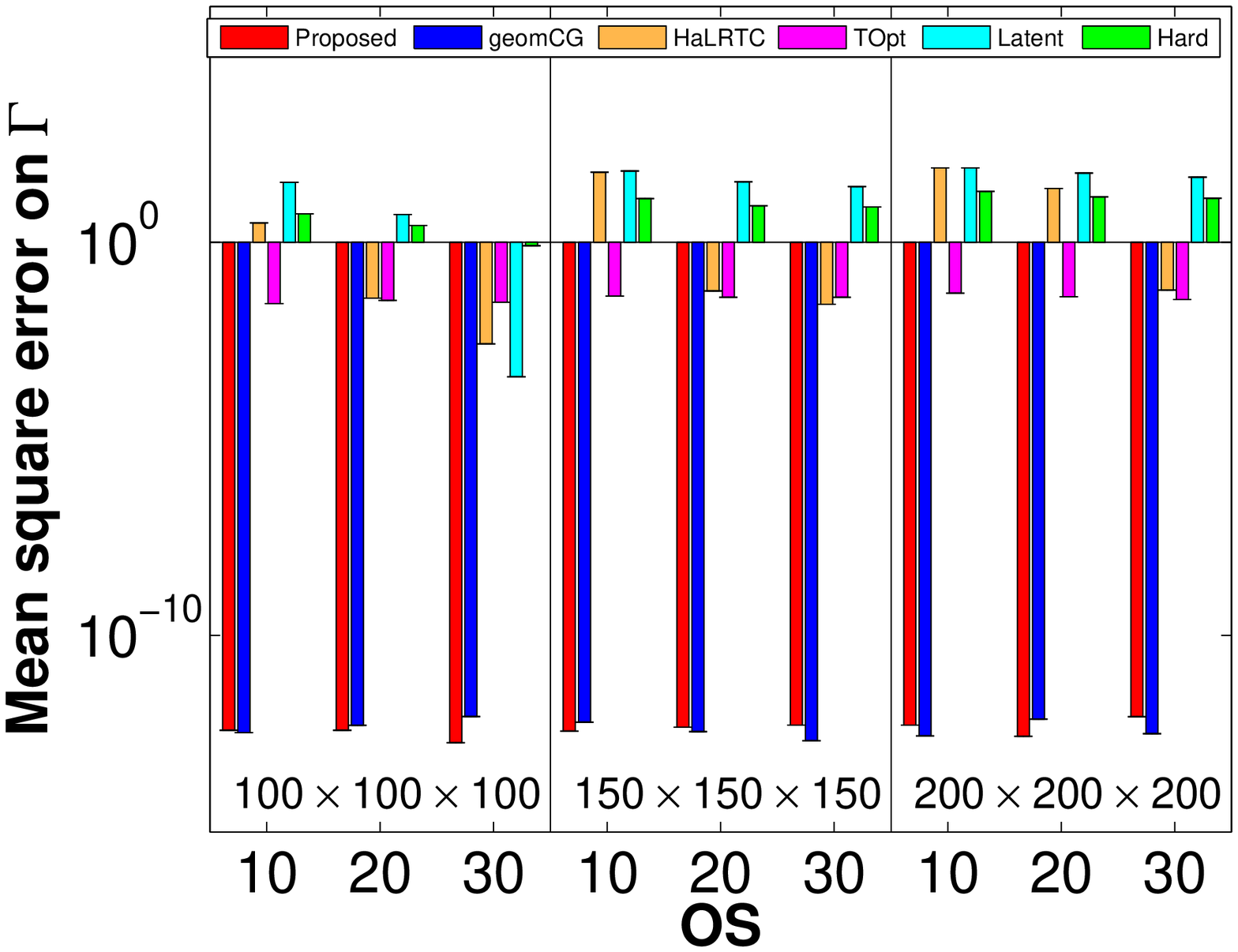}\\
{\scriptsize(d) {\bf r} = ($5,5,5$).}
\end{center}
\end{minipage}
\begin{minipage}{0.32\hsize}
\begin{center}
\includegraphics[width=\hsize]{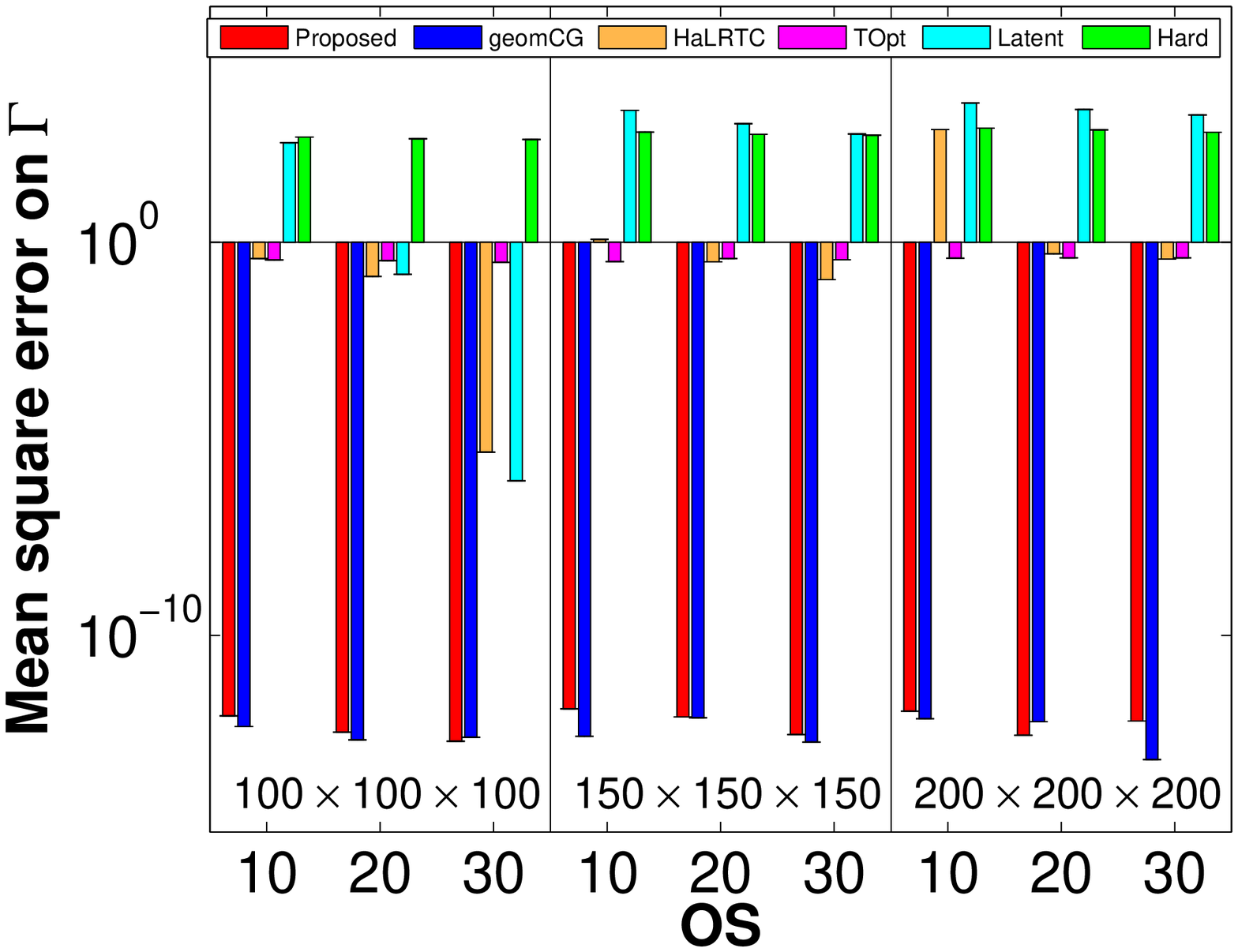}\\
{\scriptsize(e) {\bf r} = ($10,10,10$).}
\end{center}
\end{minipage}
\begin{minipage}{0.32\hsize}
\begin{center}
\includegraphics[width=\hsize]{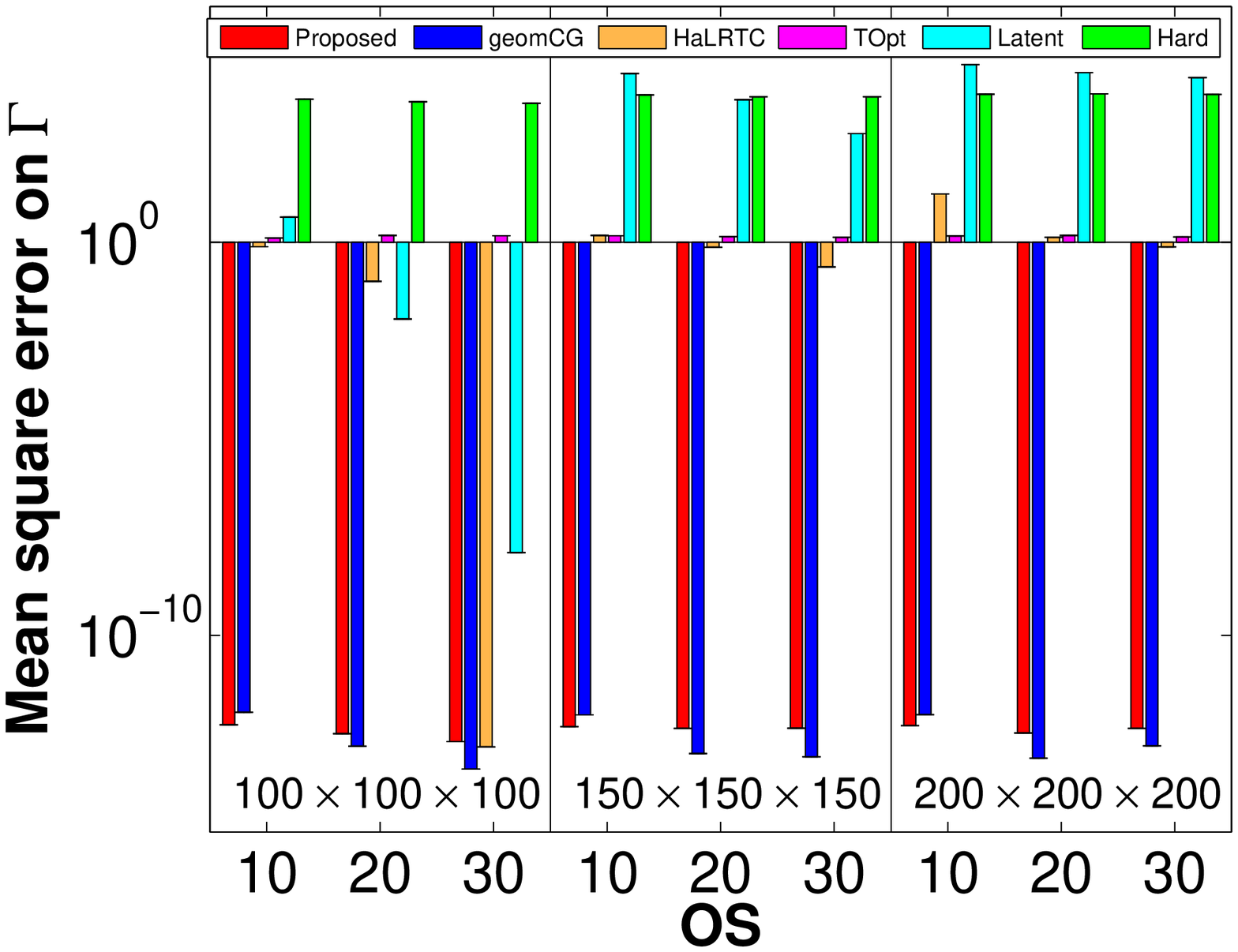}\\
{\scriptsize(f) {\bf r} = ($15,15,15$).}
\end{center}
\end{minipage}\\
\begin{minipage}{0.32\hsize}
\begin{center}
\includegraphics[width=\hsize]{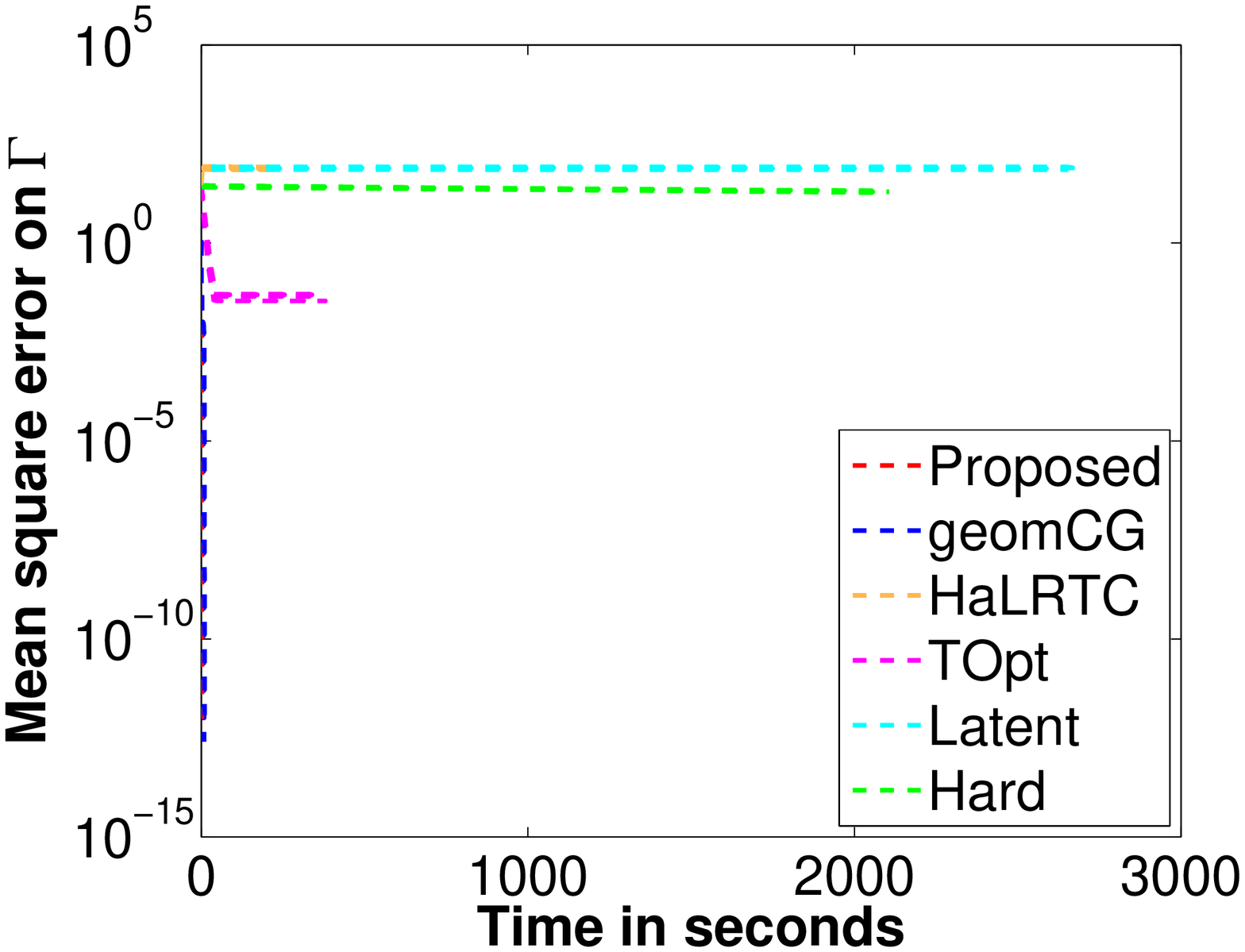}\\
{\scriptsize(g) $200\times 200\times 200$, OS = $10$, \\ {\bf r} = ($5,5,5$).}
\end{center}
\end{minipage}
\begin{minipage}{0.32\hsize}
\begin{center}
\includegraphics[width=\hsize]{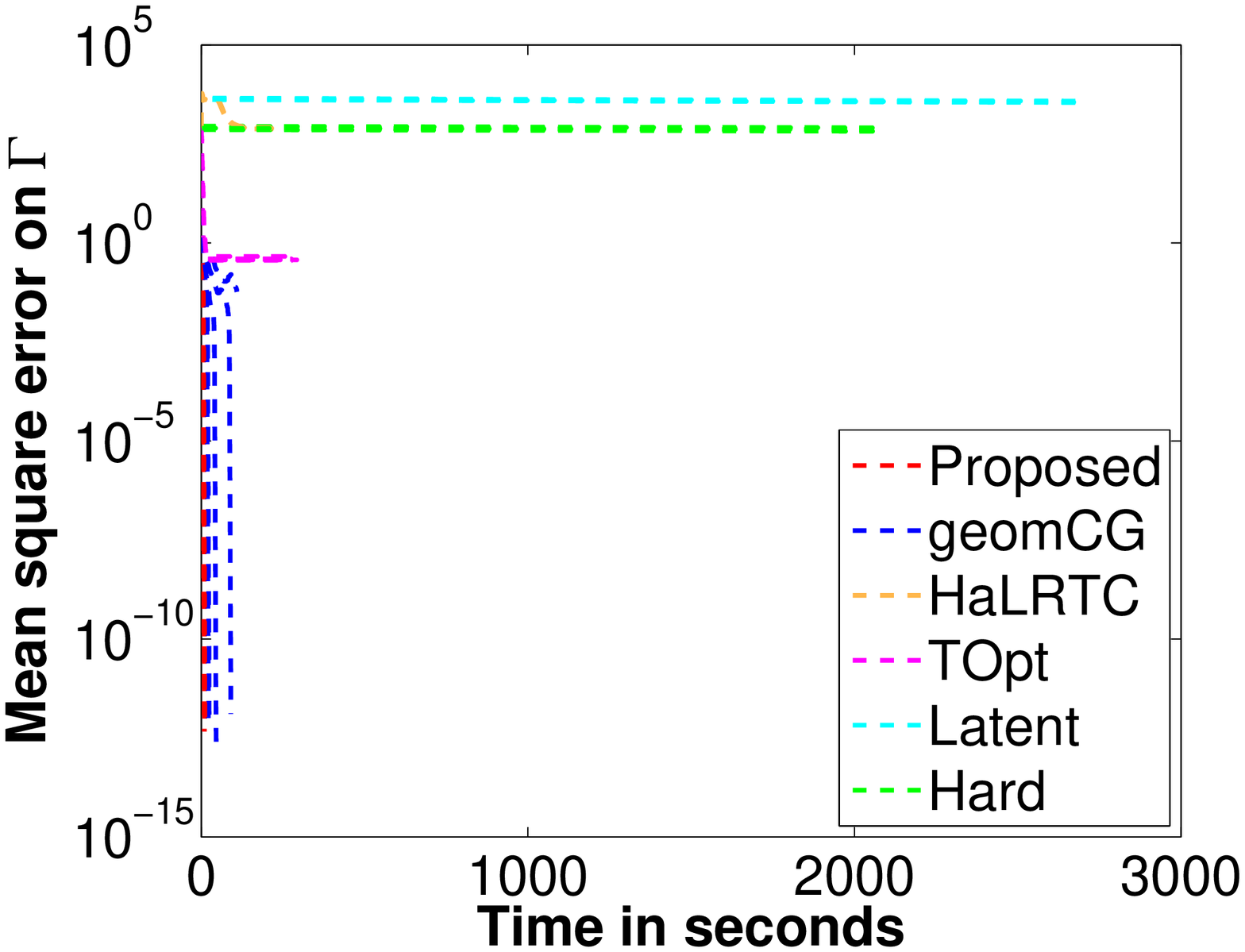}\\
{\scriptsize(h) $200\times 200\times 200$, OS = $10$, \\ {\bf r} = ($10,10,10$).}
\end{center}
\end{minipage}
\begin{minipage}{0.32\hsize}
\begin{center}
\includegraphics[width=\hsize]{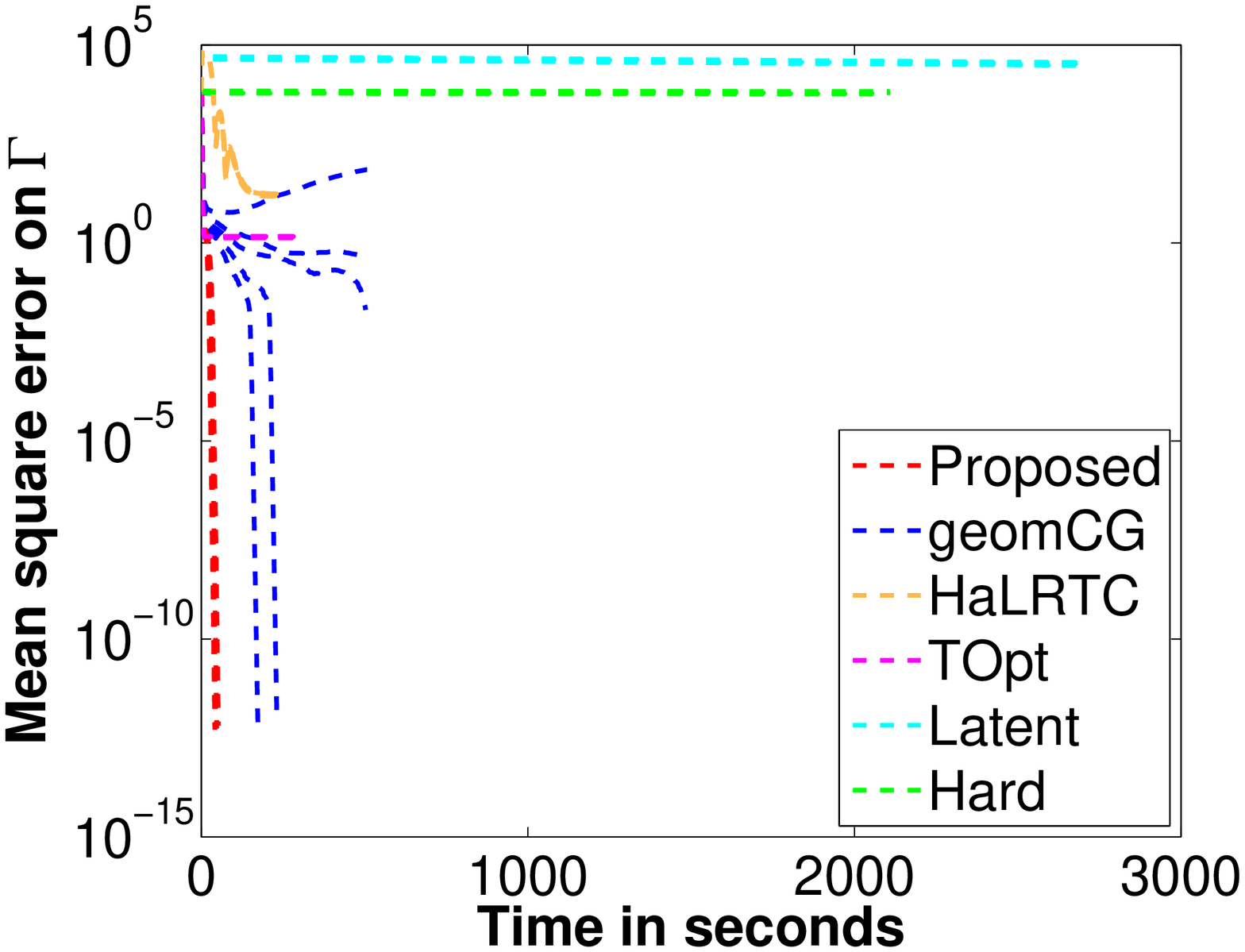}\\
{\scriptsize(i) $200\times 200\times 200$, OS = $10$, \\ {\bf r} = ($15,15,15$).}
\end{center}
\end{minipage}\\
\end{tabular}
\caption{\changeHK{\bf Case S2:} small-scale comparisons.}
\label{appnfig:small-scale}
\end{figure}

{\bf Case S2: small-scale instances.} We consider tensors of size $100 \times 100 \times 100$, $150 \times 150 \times 150$, and  $200 \times 200 \times 200$ and ranks $(5,5,5)$, $(10,10,10)$, and  $(15,15,15)$. OS is $\{10,20,30\}$. \changeHK{Figures \ref{appnfig:small-scale}(a)-(c) show the convergence behavior of different algorithms, 
where (b) is identical to the figure in the manuscript paper. Figures \ref{appnfig:small-scale}(d)-(f) show the mean square error on $\Gamma$ on each algorithm. Furthermore, Figure \ref{appnfig:small-scale}(g)-(i) show the mean square error on $\Gamma$ when OS is $10$ in all the five runs. From Figures \ref{appnfig:small-scale}, our proposed algorithm is consistently competitive or faster than geomCG, HalRTC, and TOpt. In addition, the mean square error on a test set $\Gamma$ is consistently competitive or lower than that of geomCG and HalRTC, especially for lower sampling ratios, e.g, for OS $10$.}

\begin{figure}[htbp]
\begin{tabular}{ccc}
\begin{minipage}{0.32\hsize}
\begin{center}
\includegraphics[width=\hsize]{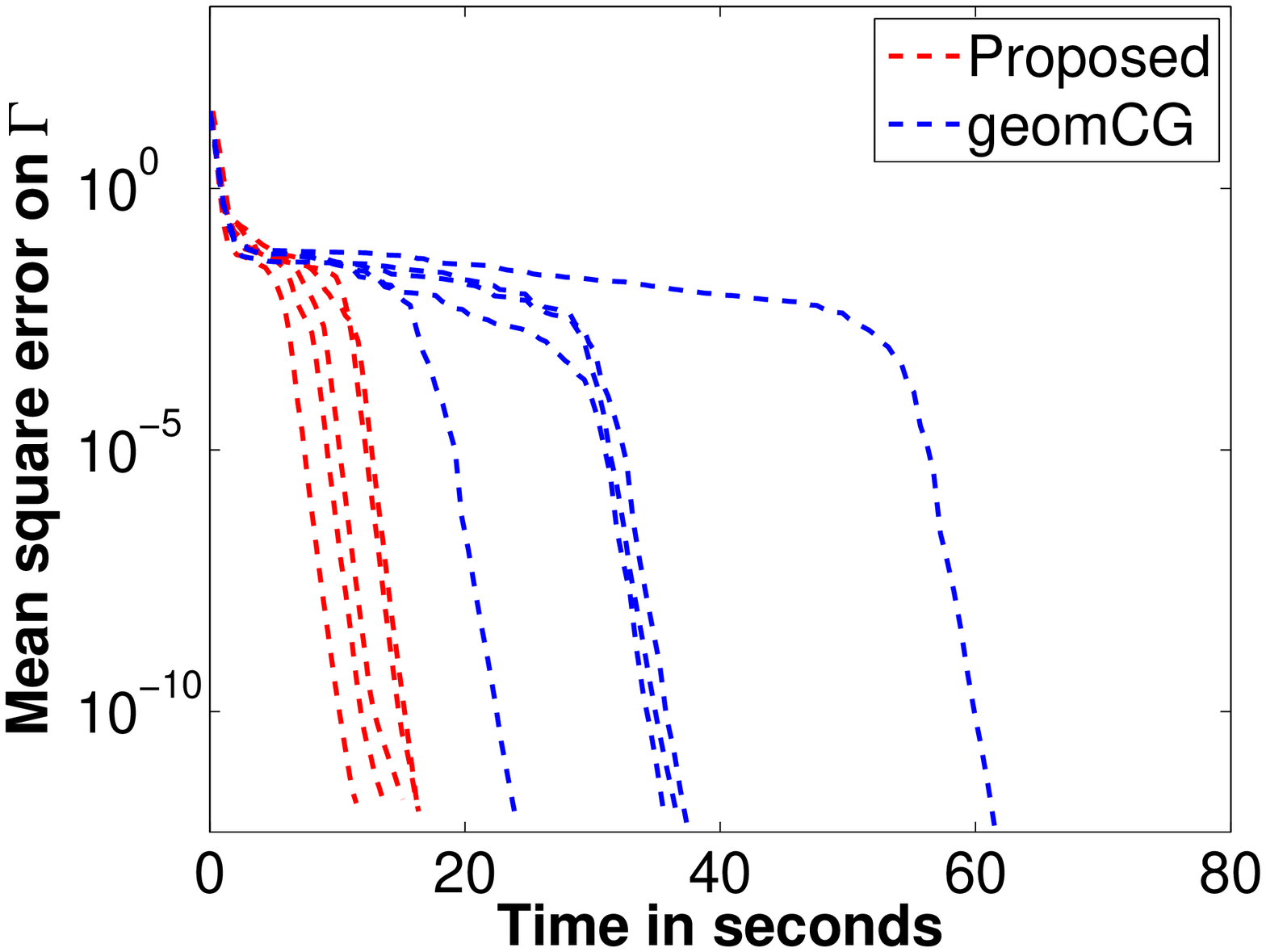}\\
{\scriptsize(a) $3000\times 3000\times 3000$, \\ {\bf r} = ($5\times 5\times 5$).}
\end{center}
\end{minipage}
\begin{minipage}{0.32\hsize}
\begin{center}
\includegraphics[width=\hsize]{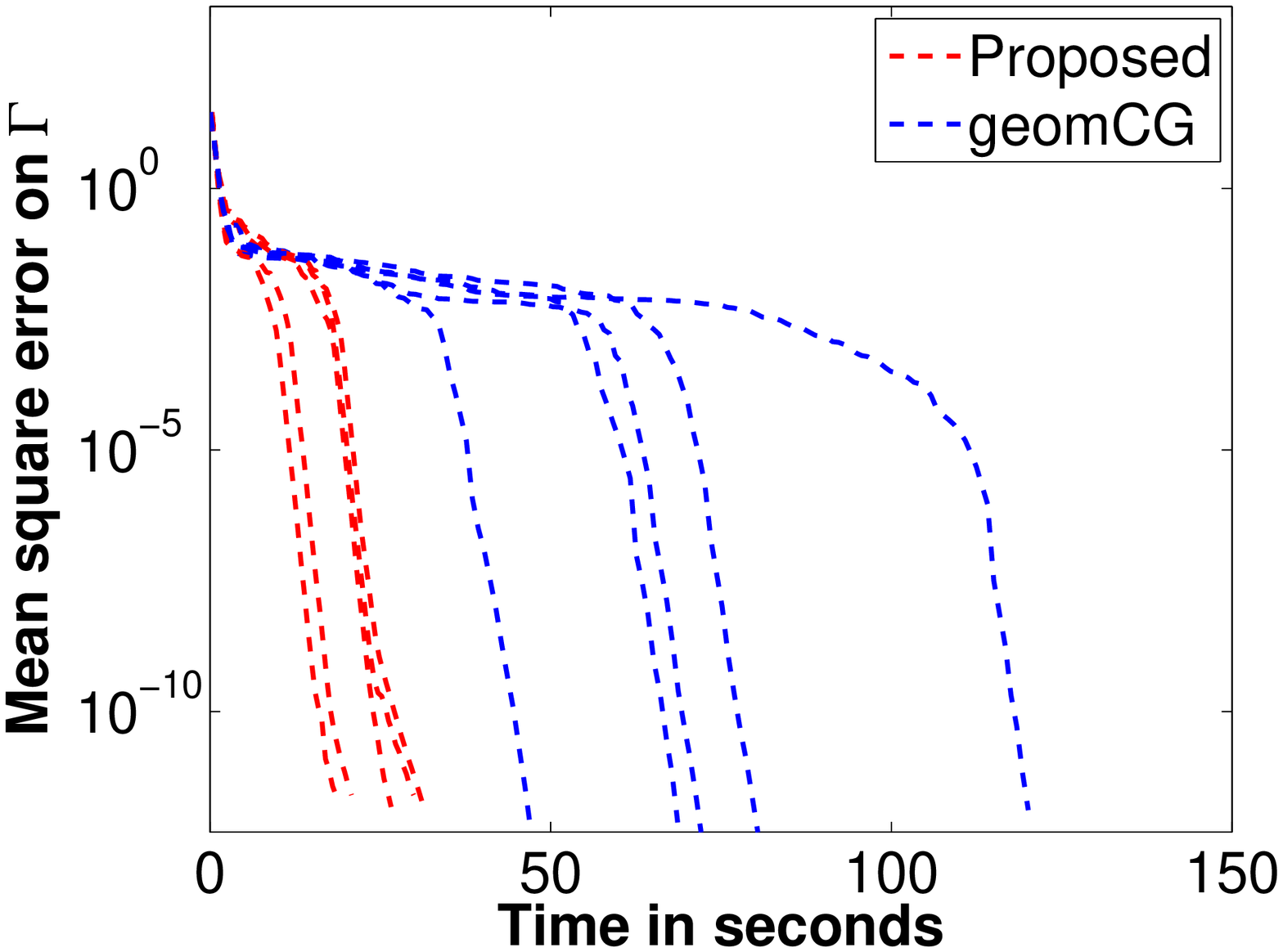}\\
{\scriptsize(b) $5000\times 5000\times 5000$, \\ {\bf r} = ($5\times 5\times 5$).}
\end{center}
\end{minipage}
\begin{minipage}{0.32\hsize}
\begin{center}
\includegraphics[width=\hsize]{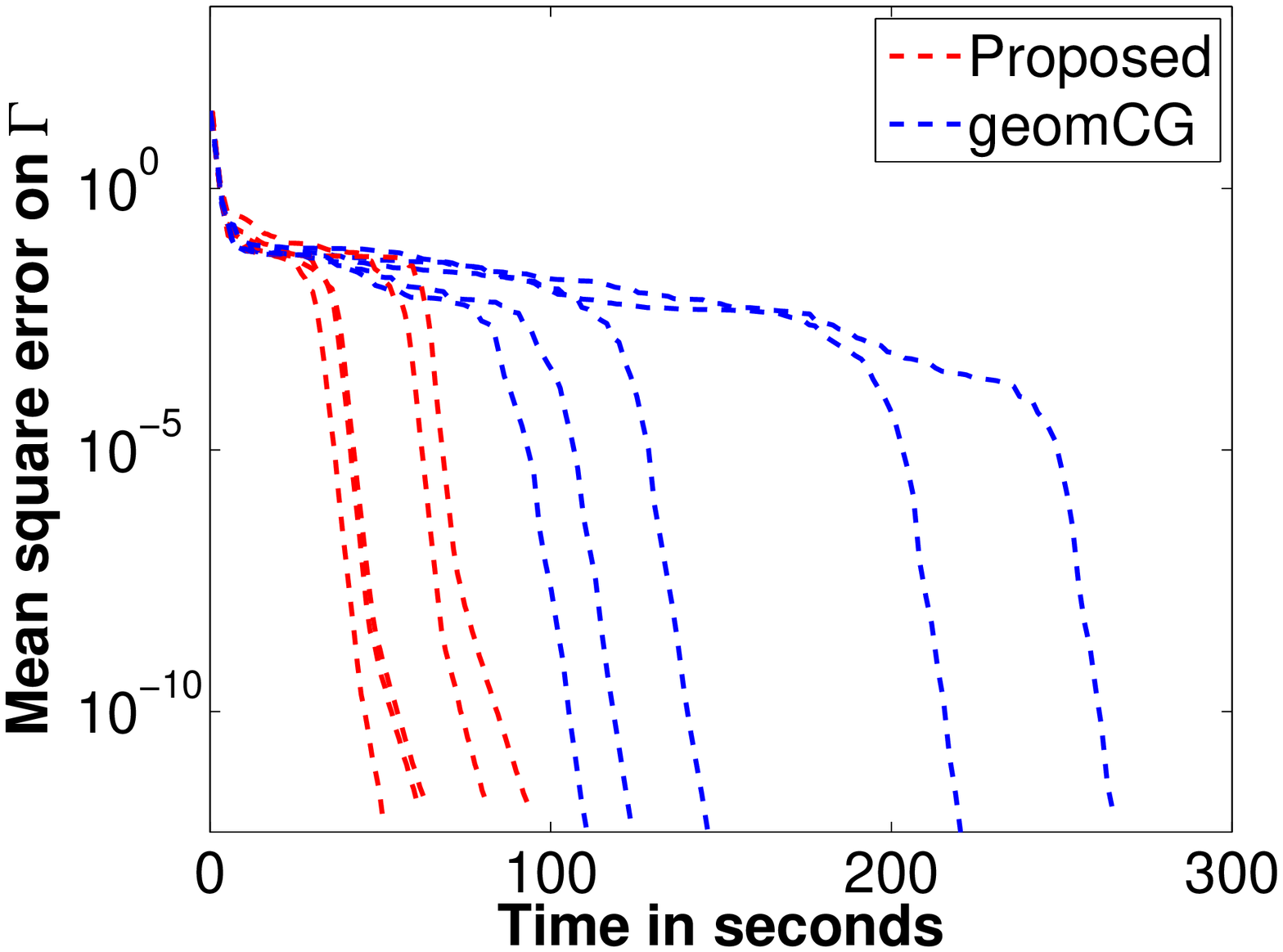}\\
{\scriptsize(c) $10000\times 10000\times 10000$, \\ {\bf r} = ($5\times 5\times 5$).}
\end{center}
\end{minipage}\\
\begin{minipage}{0.32\hsize}
\begin{center}
\includegraphics[width=\hsize]{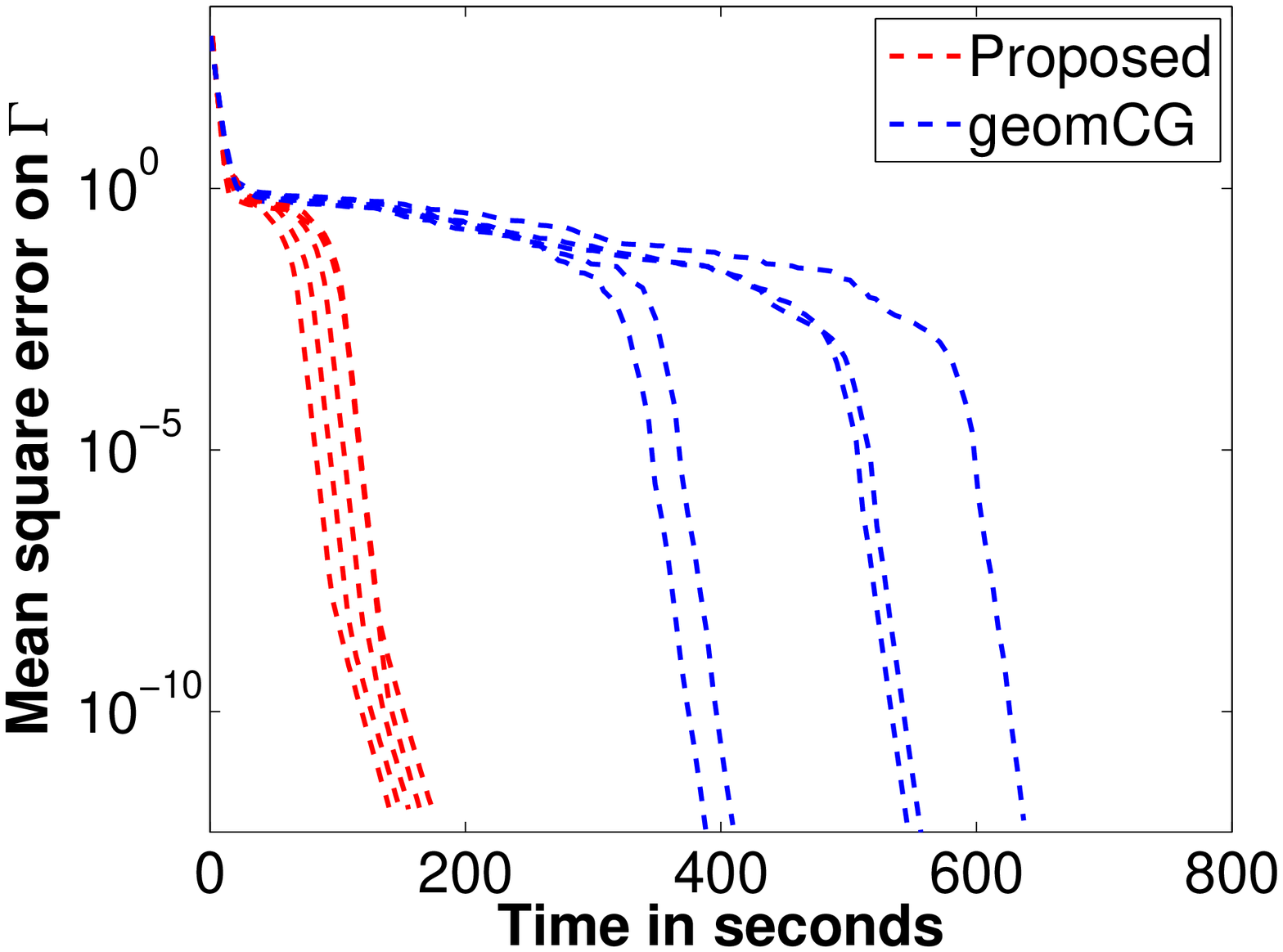}\\
{\scriptsize(d) $3000\times 3000\times 3000$, \\ {\bf r} = ($10\times 10\times 10$).}
\end{center}
\end{minipage}
\begin{minipage}{0.32\hsize}
\begin{center}
\includegraphics[width=\hsize]{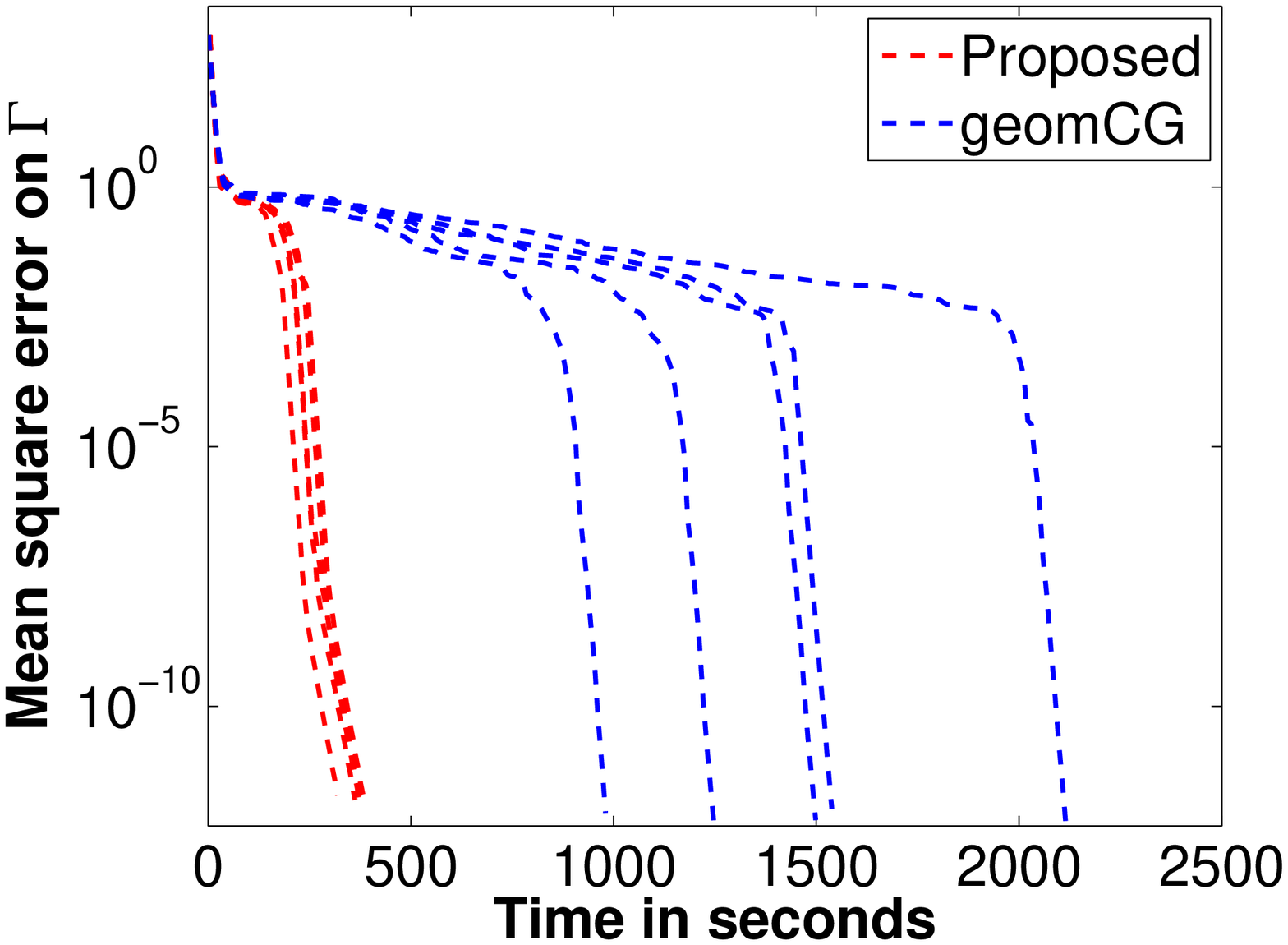}\\
{\scriptsize(e) $5000\times 5000\times 5000$, \\ {\bf r} = ($10\times 10\times 10$).}
\end{center}
\end{minipage}
\begin{minipage}{0.32\hsize}
\begin{center}
\includegraphics[width=\hsize]{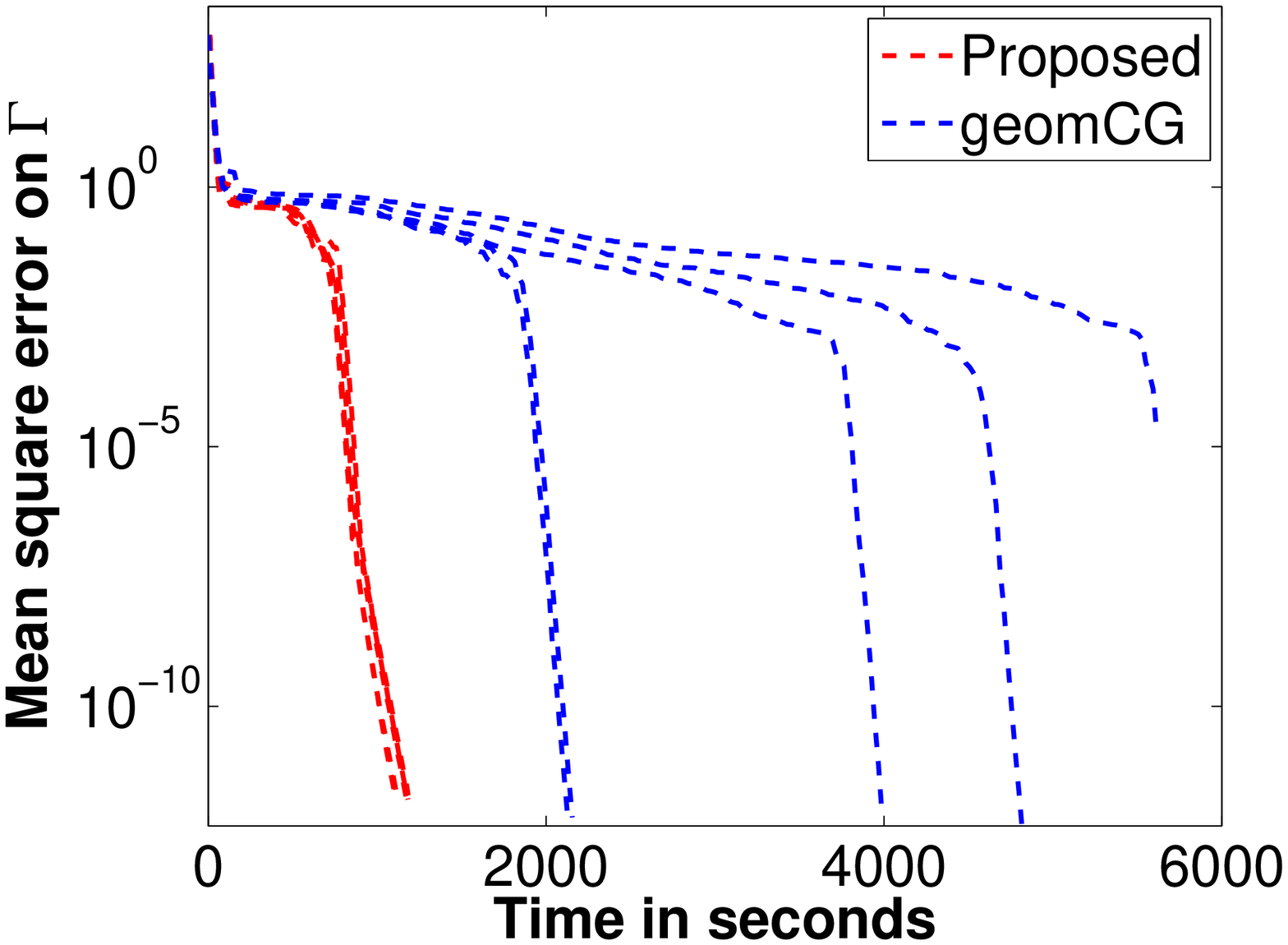}\\
{\scriptsize(f) $10000\times 10000\times 10000$, \\ {\bf r} = ($10\times 10\times 10$).}
\end{center}
\end{minipage}
\end{tabular}
\vspace{-0.1cm}
\caption{\changeHK{\bf Case S3:} large-scale comparisons.}
\label{appnfig:large-scale}
\end{figure}

{\bf Case S3: large-scale instances.} We consider large-scale tensors of size $3000 \times 3000 \times 3000$, $5000 \times 5000 \times 5000$, and $10000 \times 10000 \times 10000$ and ranks {\bf r}=$(5,5,5)$ and $(10,10,10)$. OS is $10$. We compare our proposed algorithm to geomCG. Figure \ref{appnfig:large-scale} shows the convergence behavior of the algorithms. The proposed algorithm outperforms geomCG in all the cases. 

\begin{figure}[htbp]
\vspace{-0.1cm}
\begin{tabular}{ccc}
\begin{minipage}{0.32\hsize}
\begin{center}
\includegraphics[width=\hsize]{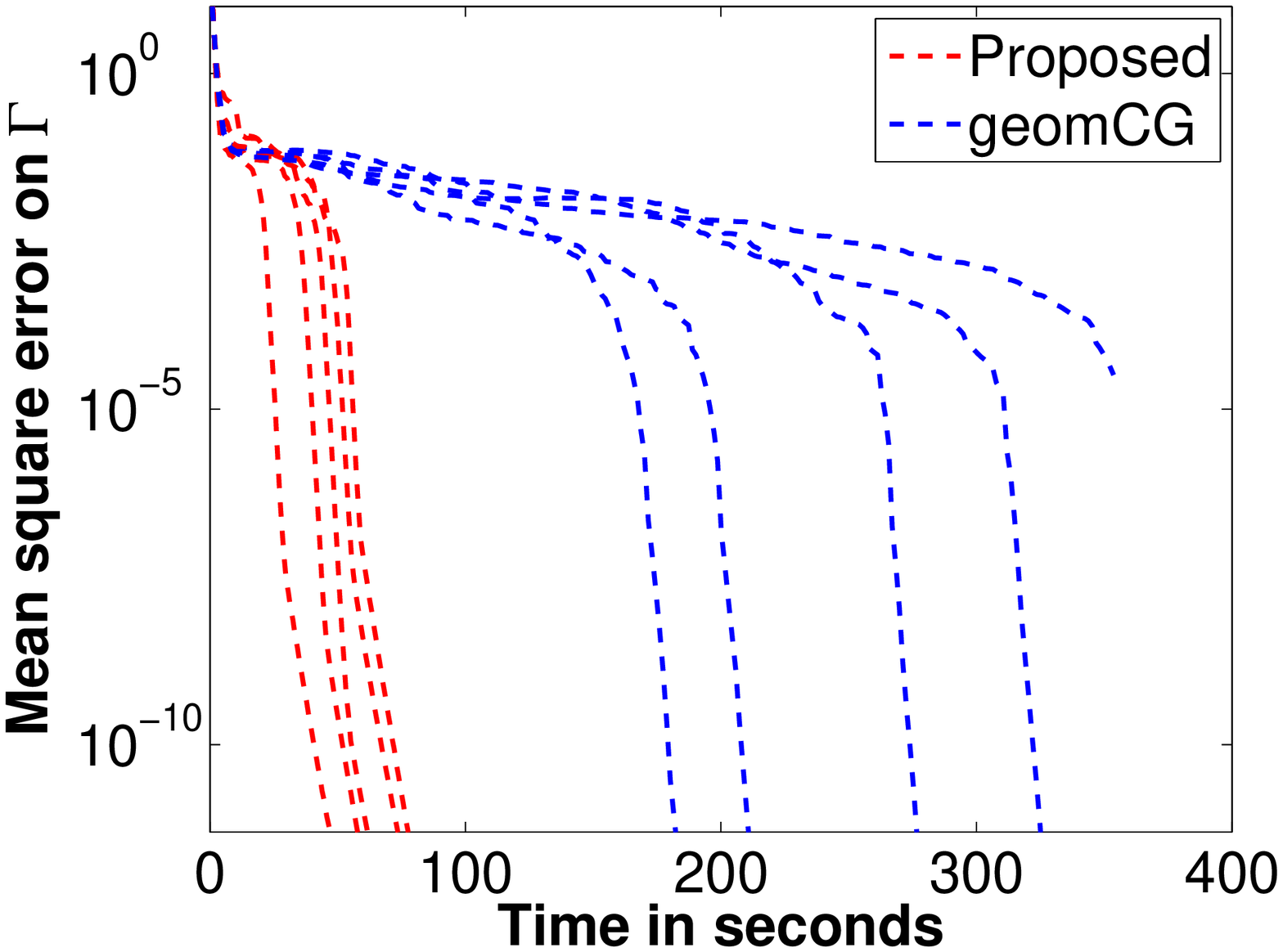}\\
{\scriptsize(a) OS = \changeHK{$8$}.}
\end{center}
\end{minipage}
\begin{minipage}{0.32\hsize}
\begin{center}
\includegraphics[width=\hsize]{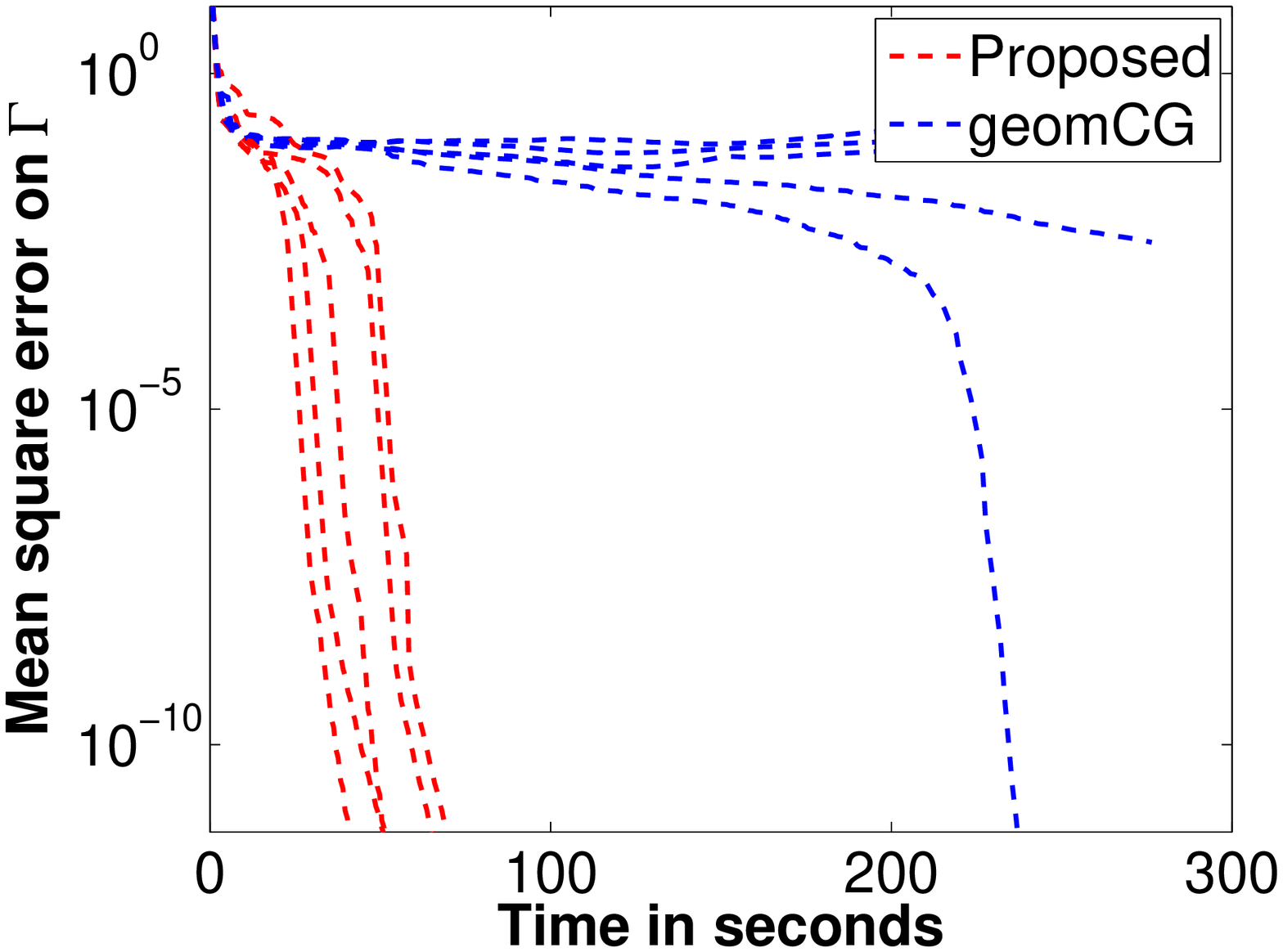}\\
{\scriptsize(b) OS = \changeHK{$6$}.}
\end{center}
\end{minipage}
\begin{minipage}{0.32\hsize}
\begin{center}
\includegraphics[width=\hsize]{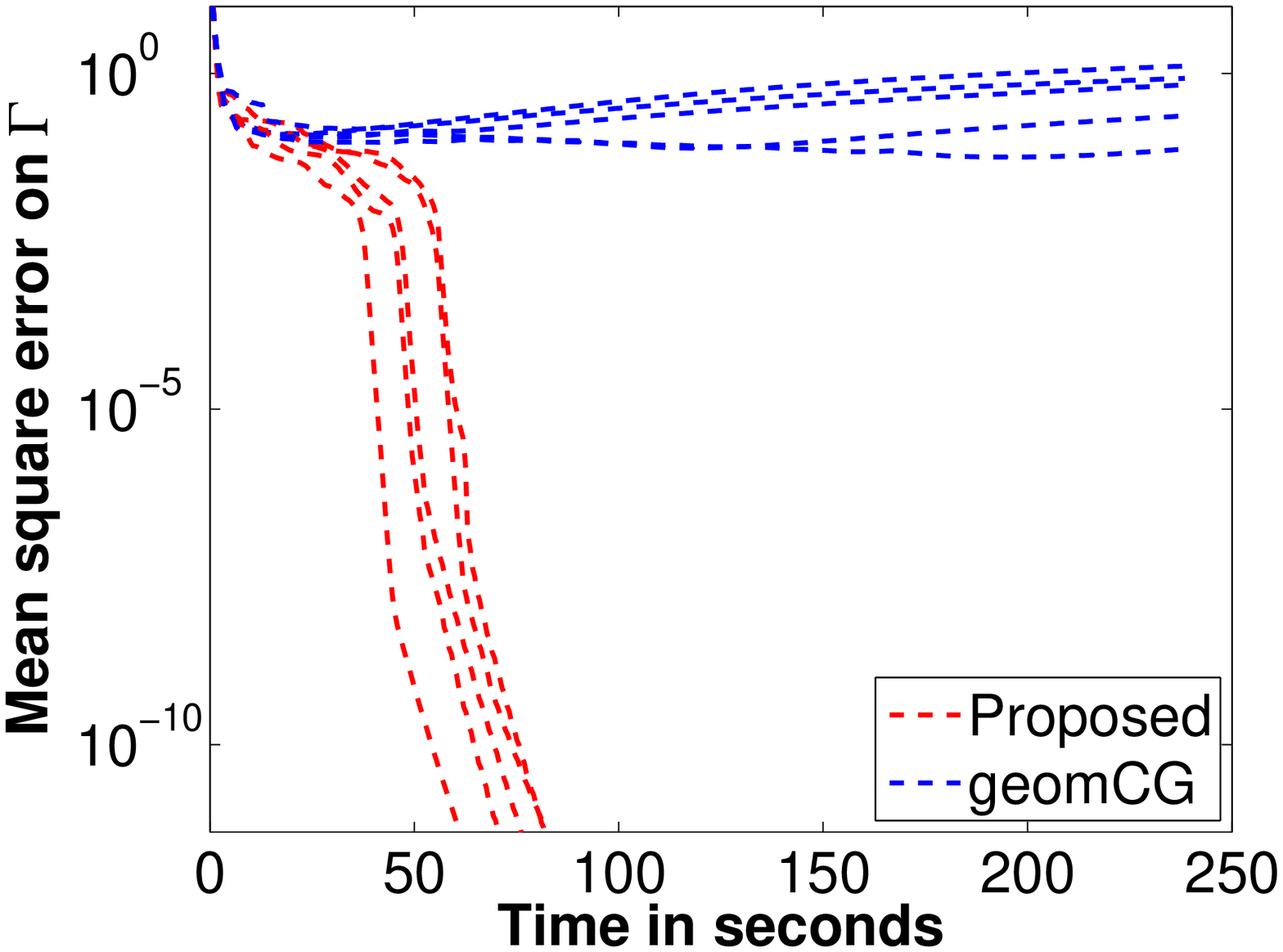}\\
{\scriptsize(c) OS = \changeHK{$5$}.}
\end{center}
\end{minipage}
\end{tabular}
\vspace{-0.2cm}
\caption{\changeHK{\bf Case S4:} low-sampling comparisons.}
\label{appnfig:low-sampling}
\end{figure}

{\bf Case S4: influence of low sampling.} We look into problem instances which result from scarcely sampled data. The test requires completing a tensor of size $10000 \times 10000 \times 10000$ and rank \changeHK{{\bf r}=$(5,5,5)$. Figure \ref{appnfig:low-sampling} shows the convergence behavior when OS is $\{8,6,5\}$.} The case of  $\text{OS}=5$ is particularly interesting. In this case, while the mean square error on $\Gamma$ increases for geomCG, the proposed algorithm stably decreases the error in all the five runs.  

\begin{figure}[htbp]
\vspace{-0.01cm}
\begin{tabular}{ccc}
\begin{minipage}{0.32\hsize}
\begin{center}
\includegraphics[width=\hsize]{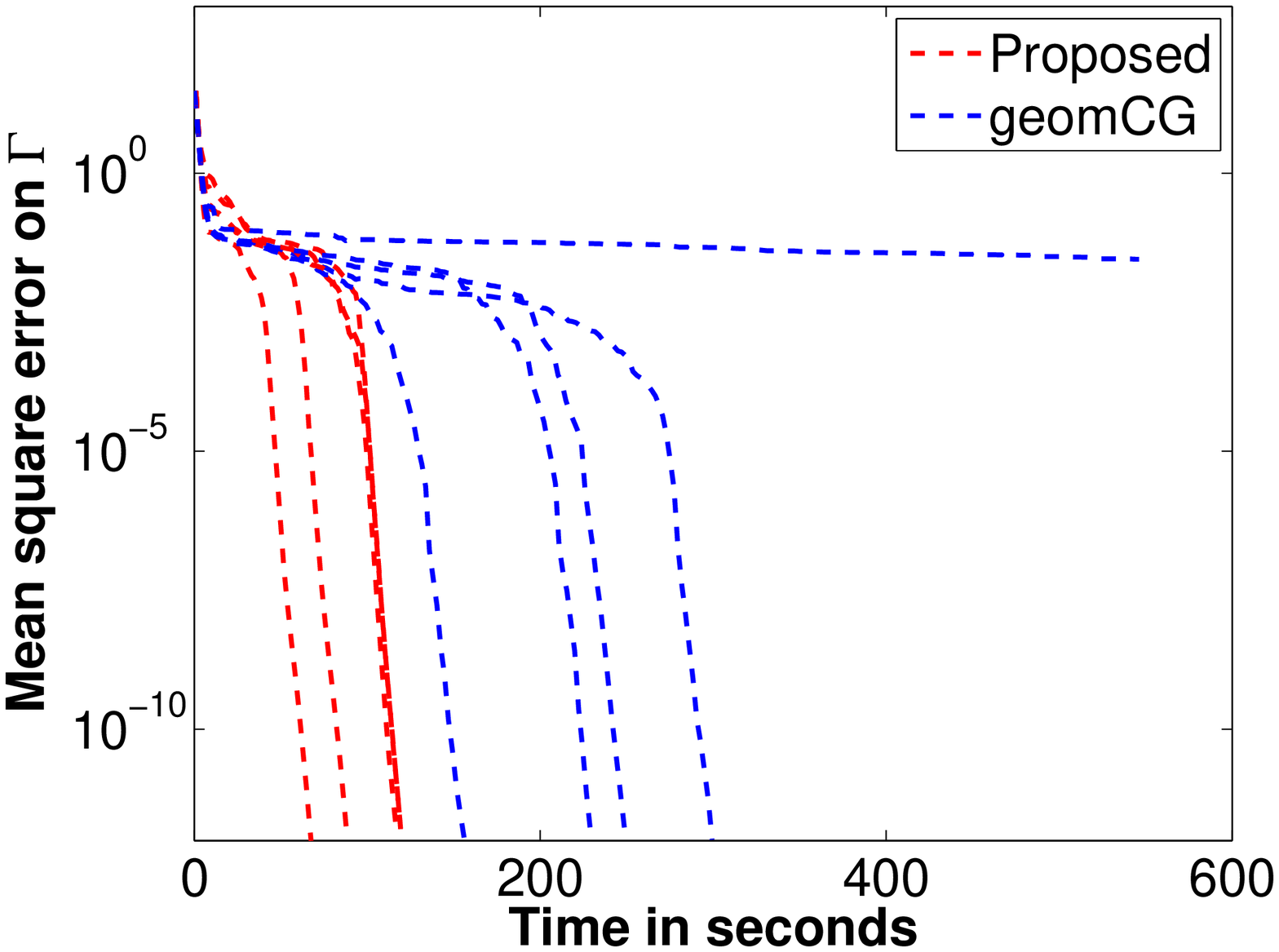}\\
{\scriptsize(a) $20000\times 7000\times 7000$, \\
{\bf r} = ($5\times 5\times 5$).}
\end{center}
\end{minipage}
\begin{minipage}{0.32\hsize}
\begin{center}
\includegraphics[width=\hsize]{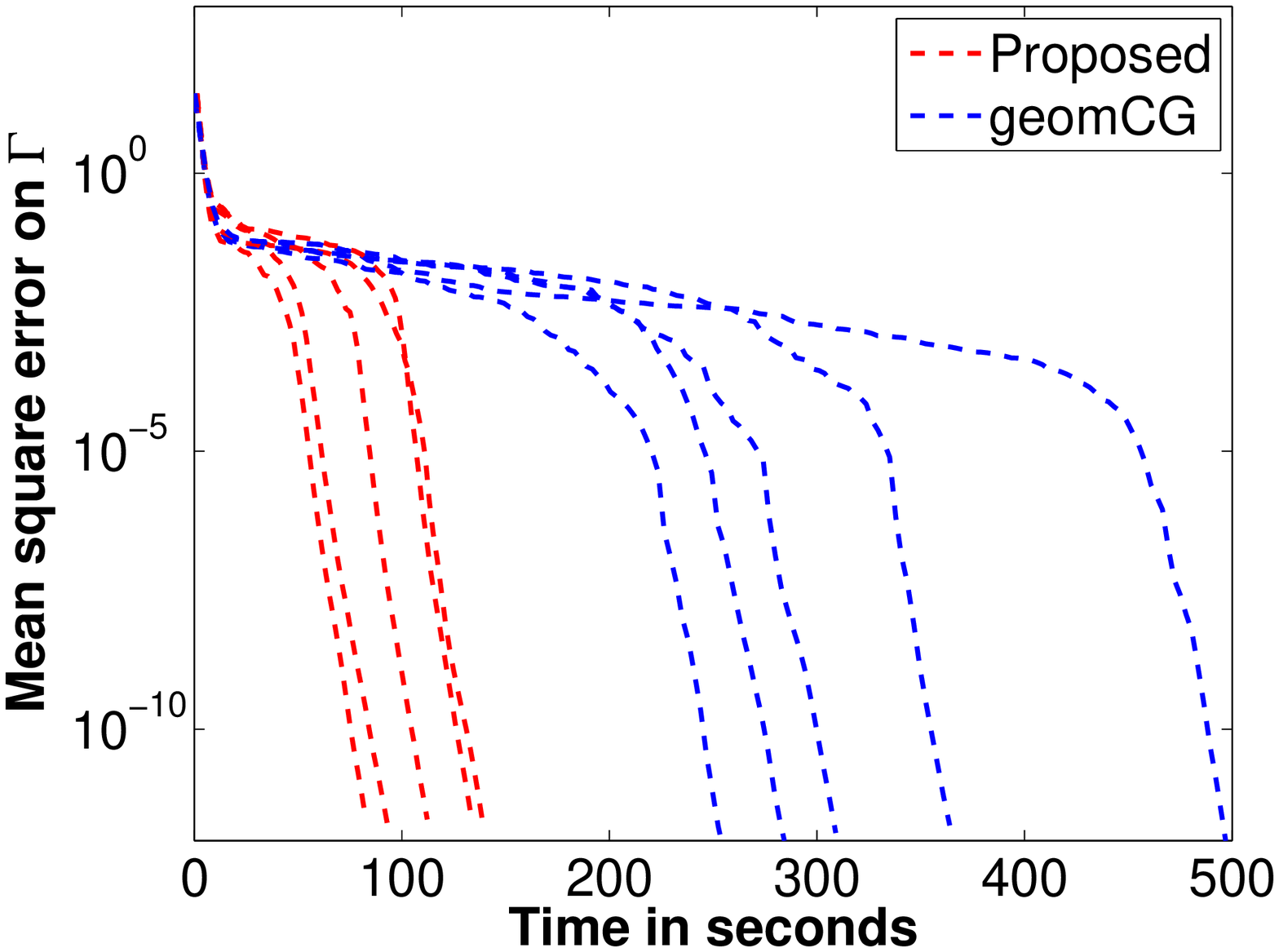}\\
{\scriptsize(b) $30000\times 60000\times 60000$, \\
{\bf r} = ($5\times 5\times 5$).}
\end{center}
\end{minipage}
\begin{minipage}{0.32\hsize}
\begin{center}
\includegraphics[width=\hsize]{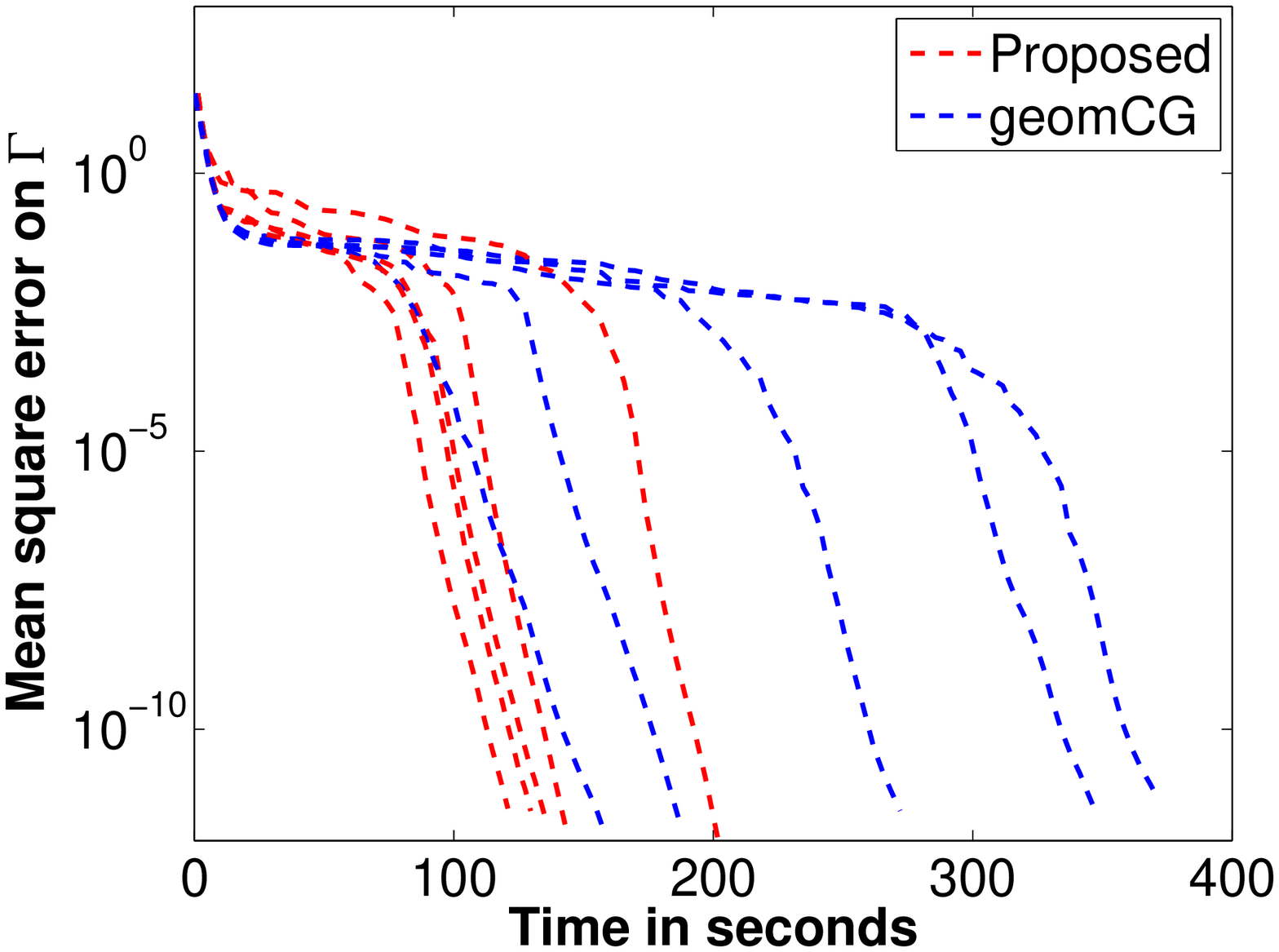}\\
{\scriptsize(c) $40000\times 5000\times 5000$, \\
{\bf r} = ($5\times 5\times 5$).}
\end{center}
\end{minipage}\\

\begin{minipage}{0.32\hsize}
\begin{center}
\includegraphics[width=\hsize]{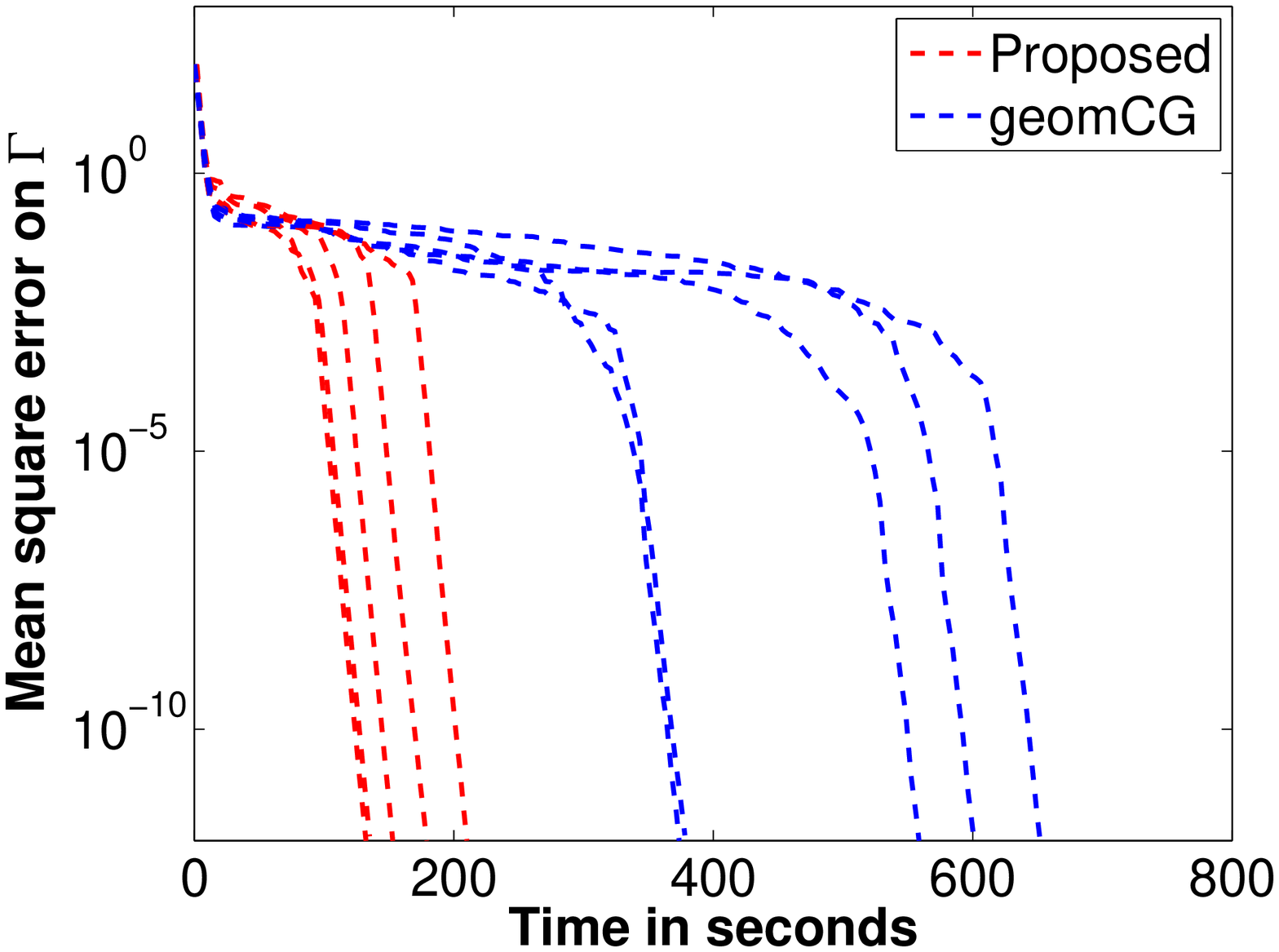}\\
{\scriptsize(d) {\bf r} = ($7\times 6\times 6$),\\
$10000\times 10000\times 10000$.}
\end{center}
\end{minipage}
\begin{minipage}{0.32\hsize}
\begin{center}
\includegraphics[width=\hsize]{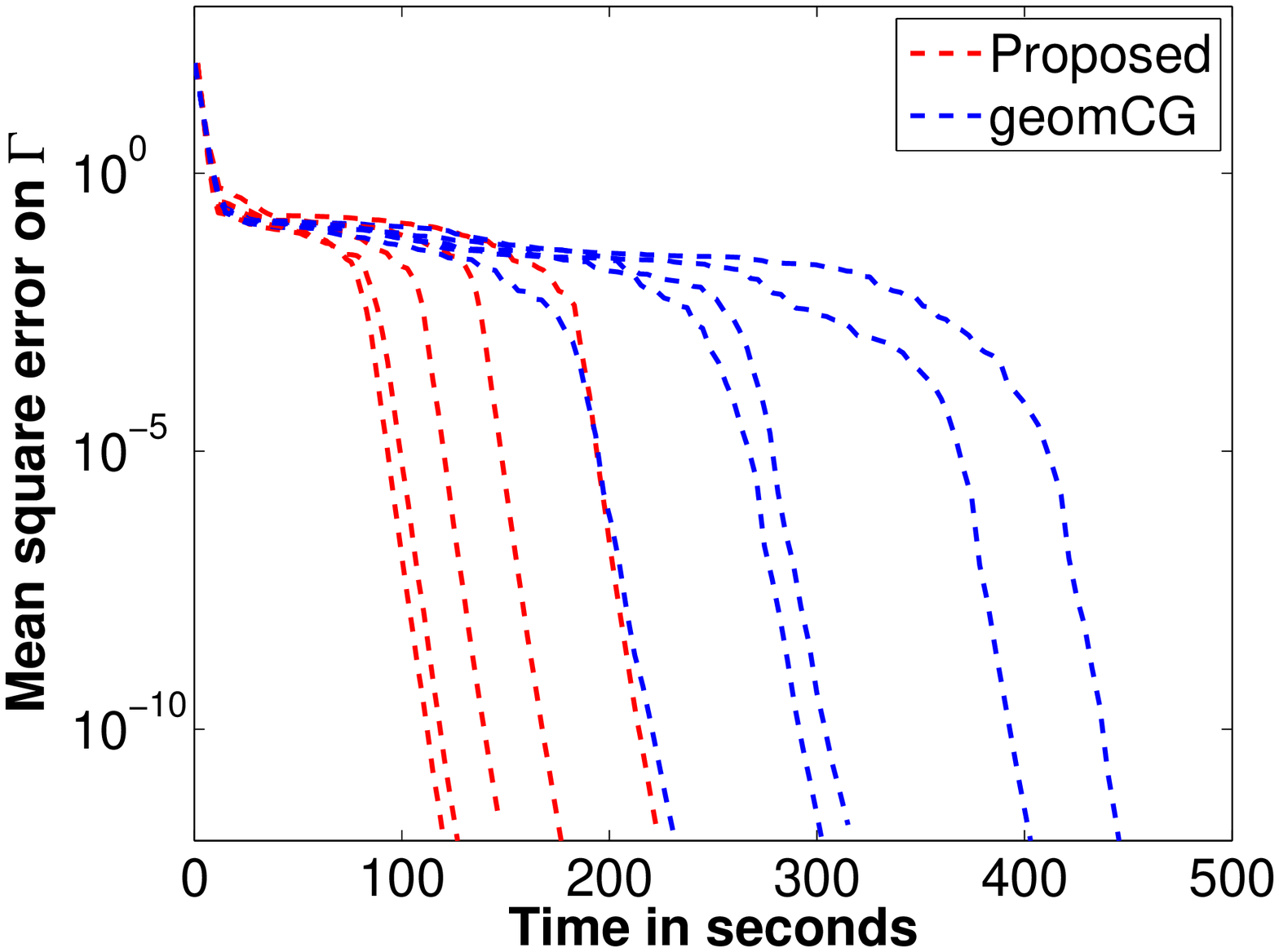}\\
{\scriptsize(e) {\bf r} = ($10\times 5\times 5$),\\
$10000\times 10000\times 10000$.}
\end{center}
\end{minipage}
\begin{minipage}{0.32\hsize}
\begin{center}
\includegraphics[width=\hsize]{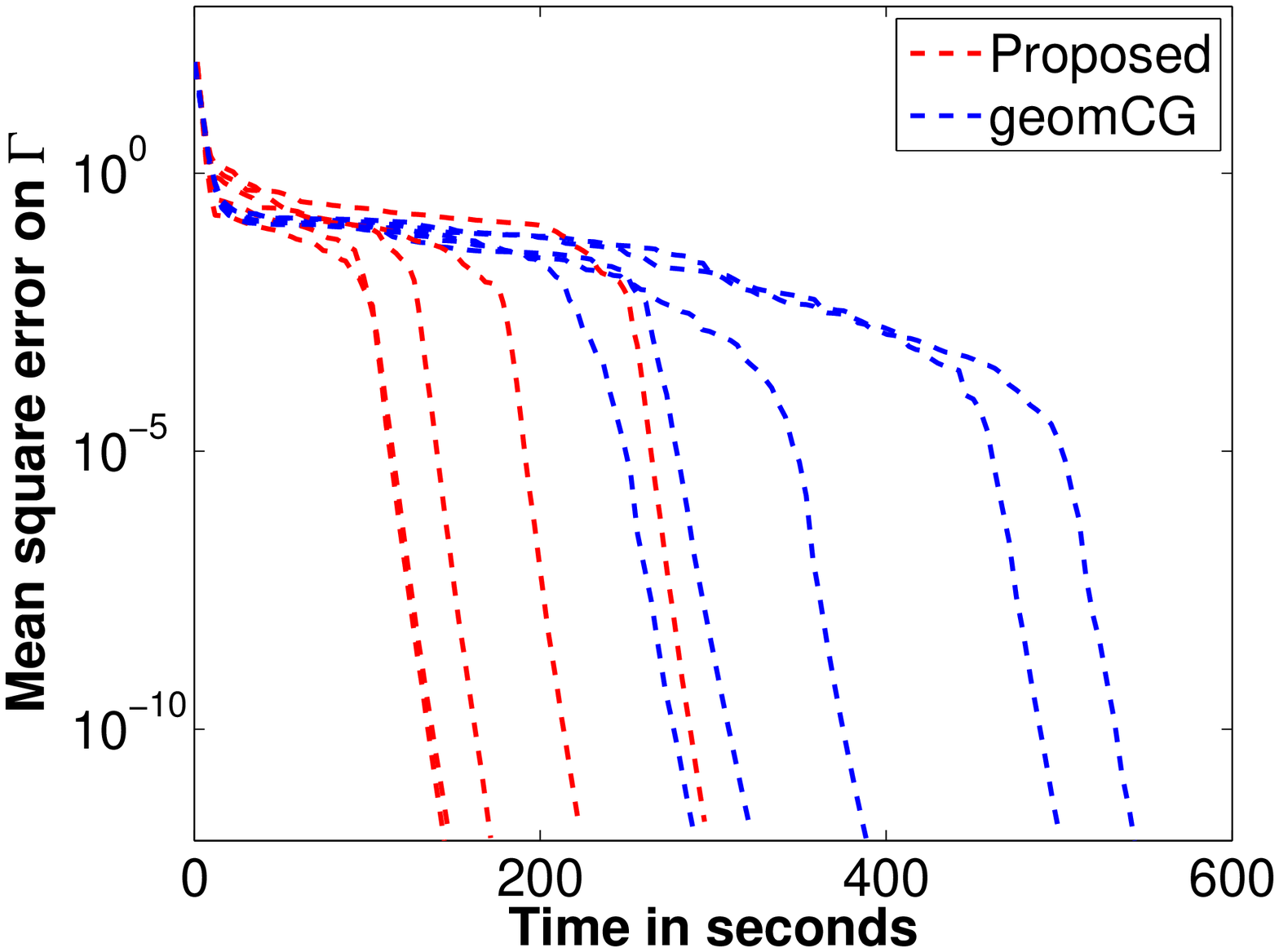}\\
{\scriptsize(f) {\bf r} = ($15\times 4\times 4$),\\
$10000\times 10000\times 10000$.}
\end{center}
\end{minipage}\\
\end{tabular}
\vspace{-0.1cm}
\caption{\changeHK{\bf Case S7:} asymmetric comparisons.}
\label{appnfig:asymmetric}
\end{figure}

\changeHK{
{\bf Case S7: asymmetric instances.} We consider instances where dimensions and ranks along certain modes are different than others. Two cases are considered. Case (7.a) considers tensors size $20000 \times7000 \times 7000$, $30000 \times 6000 \times 6000$, and $40000 \times 5000\times 5000$ and rank ${\bf r}=(5,5,5)$. Case (7.b) considers a tensor of size $10000 \times10000 \times10000$ with ranks ${\bf r}=(7,6,6)$, $(10,5,5)$, and $ (15,4,4)$. Figures \ref{appnfig:asymmetric}(a)-(c) show that the convergence behavior of our proposed algorithm is superior to that of geomCG. Our proposed algorithm also outperforms geomCG for the asymmetric rank cases as shown in Figure \ref{appnfig:asymmetric}(d)-(f).
}

\begin{figure}[htbp]
\vspace{-0.05cm}
\begin{tabular}{cc}
\begin{minipage}{0.32\hsize}
\begin{center}
\includegraphics[width=\hsize]{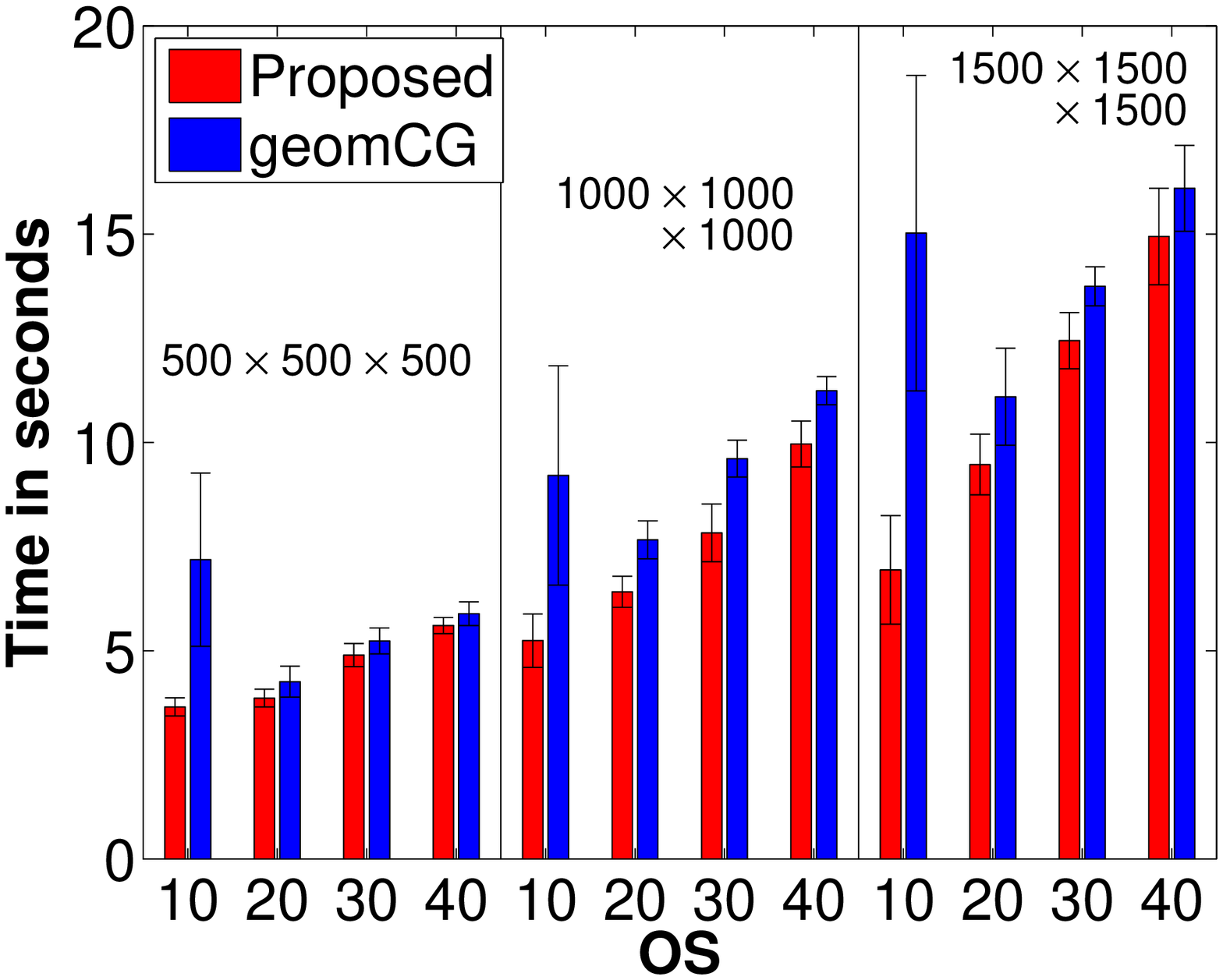}\\
{\scriptsize(a) {\bf r} = ($5\times 5\times 5$).}
\end{center}
\end{minipage}
\begin{minipage}{0.32\hsize}
\begin{center}
\includegraphics[width=\hsize]{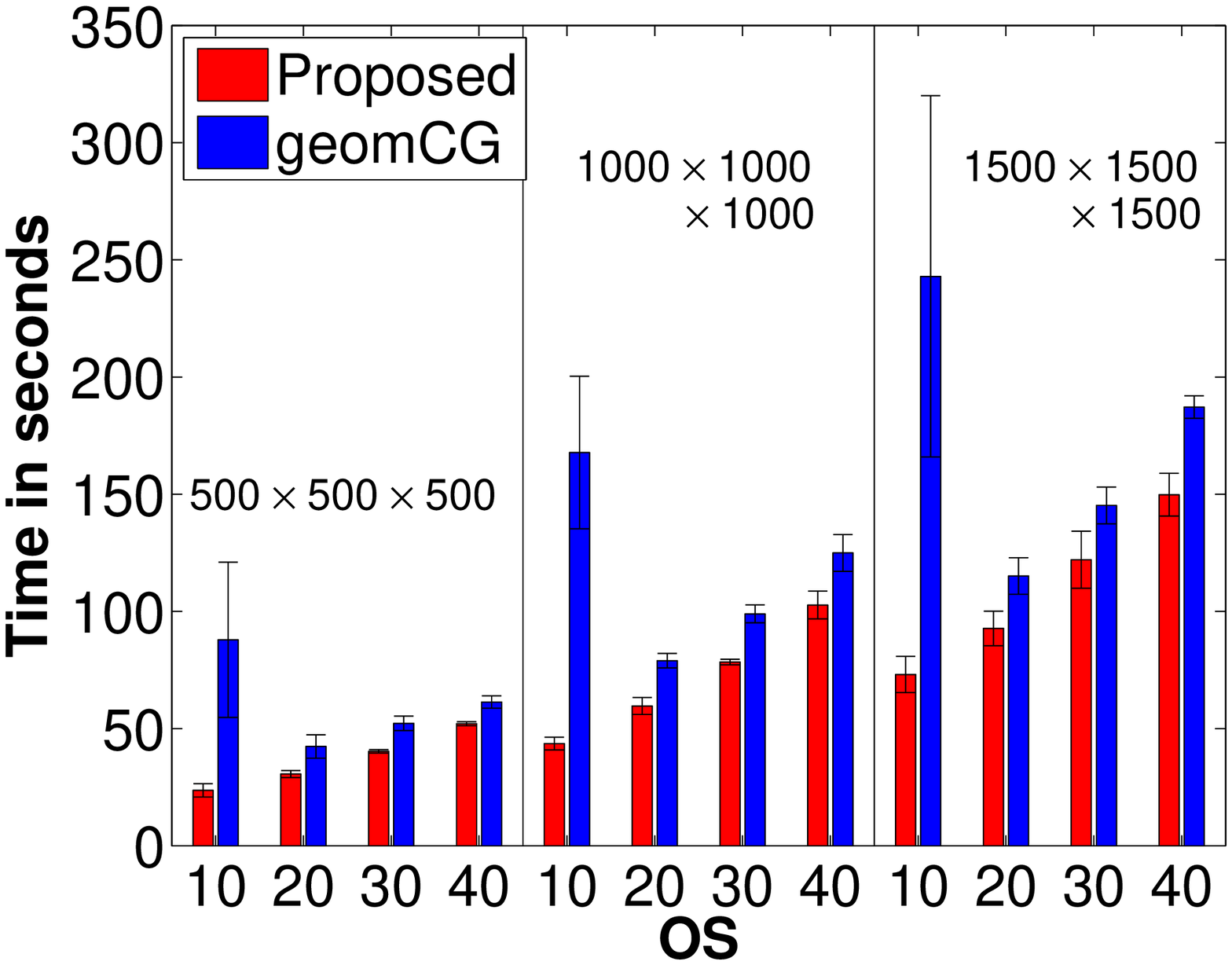}\\
{\scriptsize(b) {\bf r} = ($10\times 10\times 10$).}
\end{center}
\end{minipage}
\begin{minipage}{0.32\hsize}
\begin{center}
\includegraphics[width=\hsize]{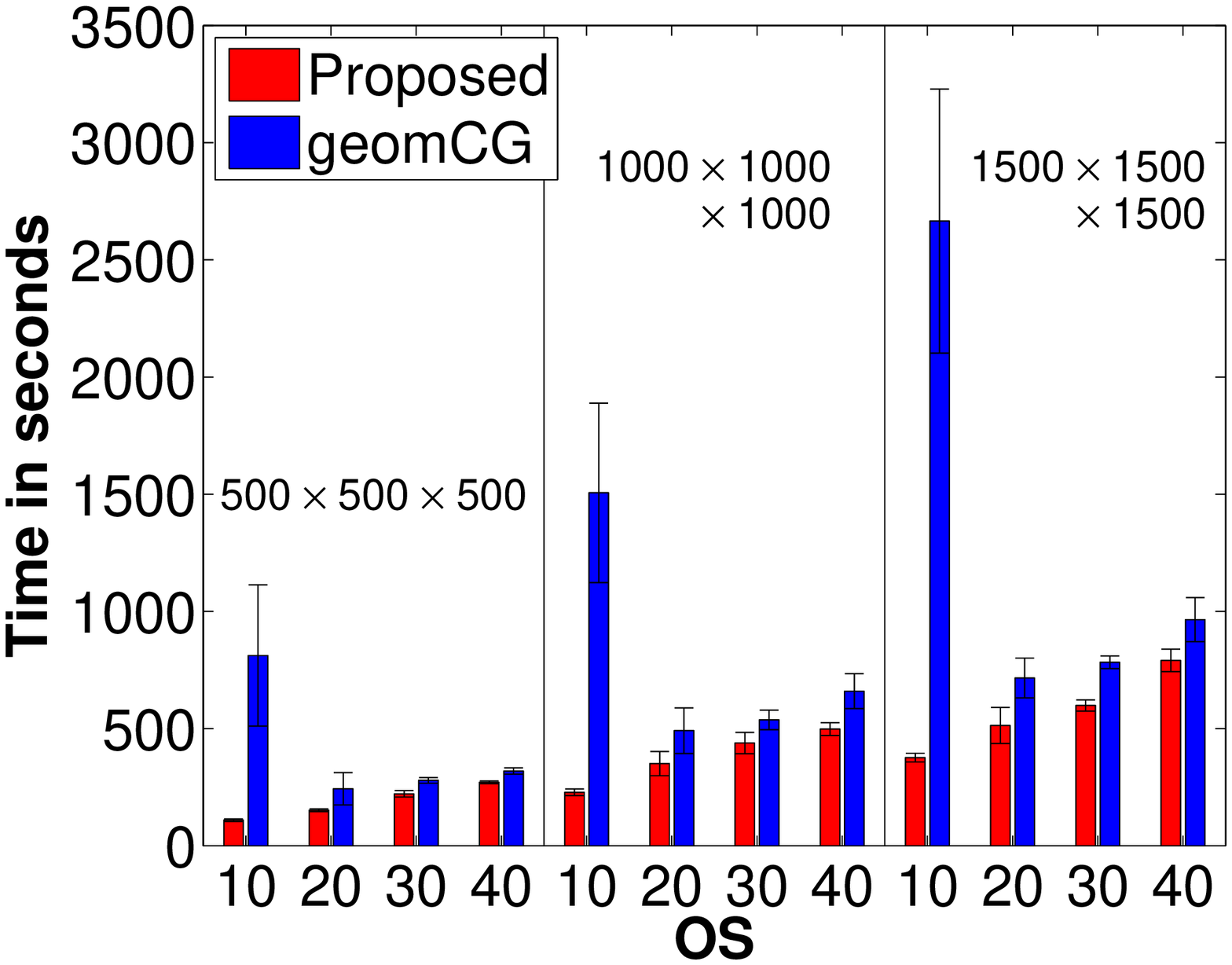}\\
{\scriptsize(c) {\bf r} = ($15\times 15\times 15$).}
\end{center}
\end{minipage}\\
\begin{minipage}{0.32\hsize}
\begin{center}
\includegraphics[width=\hsize]{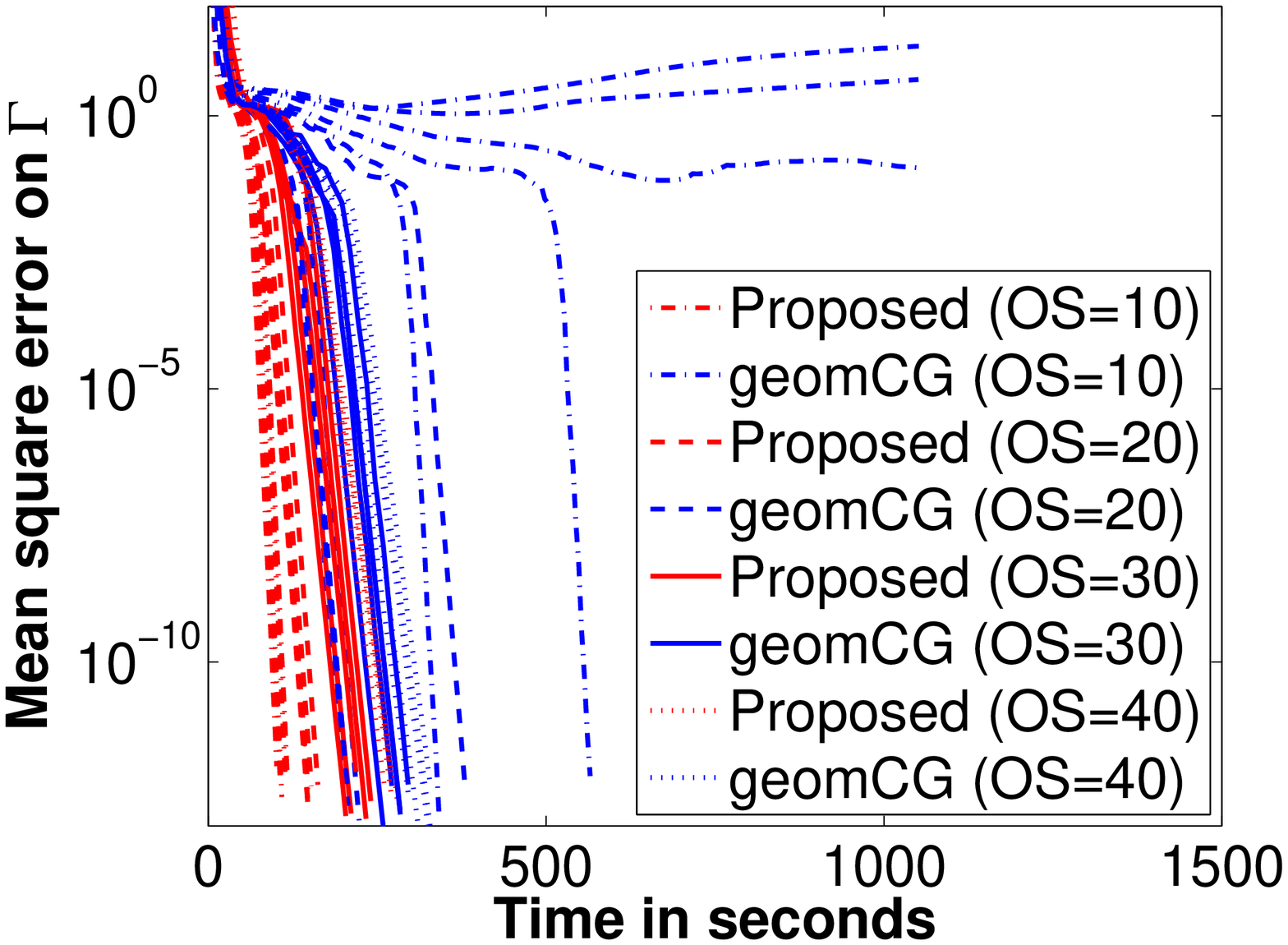}\\
{\scriptsize(d) $500\times 500\times 500$,\\
{\bf r} = ($15\times 15\times 15$).}
\end{center}
\end{minipage}
\begin{minipage}{0.32\hsize}
\begin{center}
\includegraphics[width=\hsize]{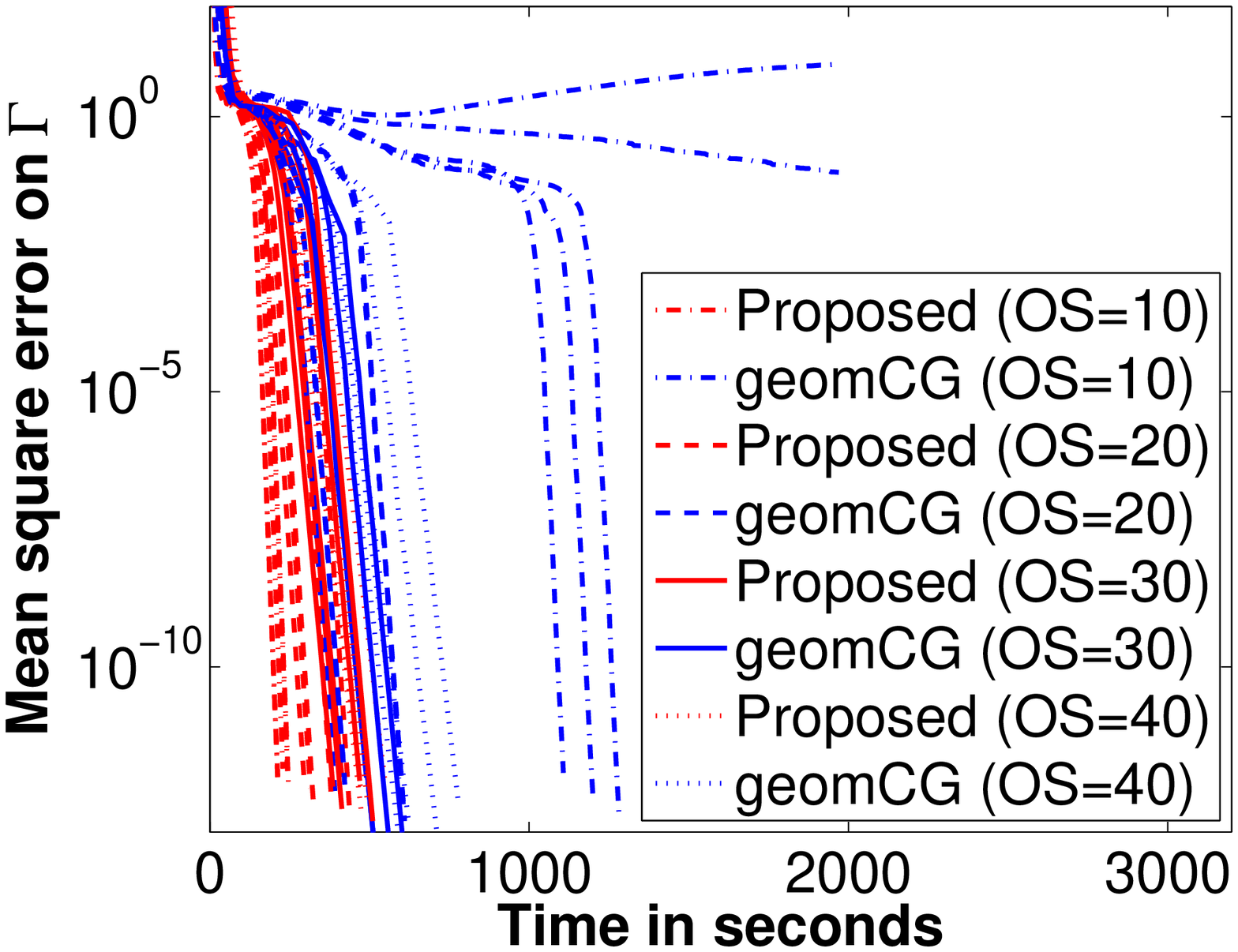}\\
{\scriptsize(e) $1000\times 1000\times 1000$,\\
{\bf r} = ($15\times 15\times 15$).}
\end{center}
\end{minipage}
\begin{minipage}{0.32\hsize}
\begin{center}
\includegraphics[width=\hsize]{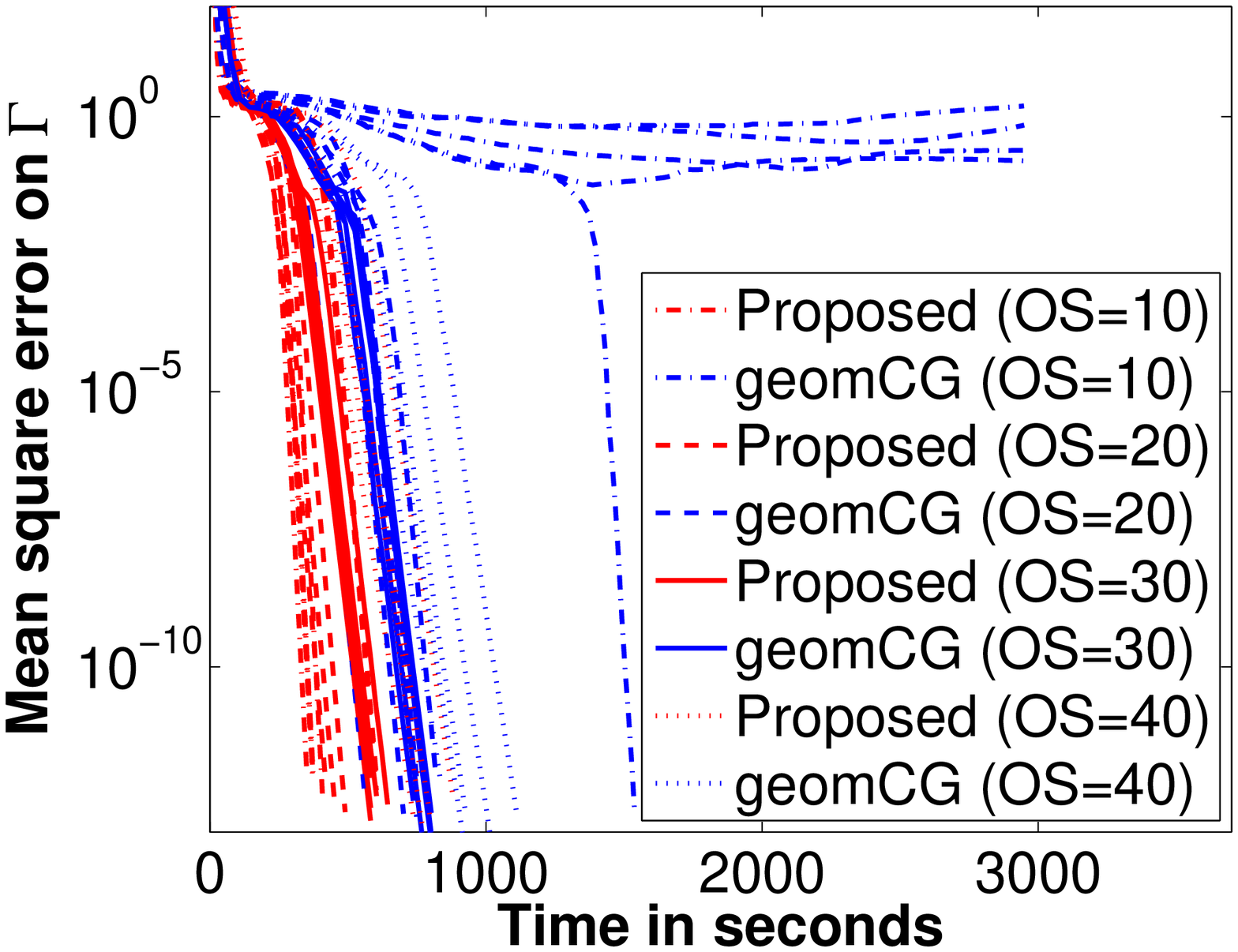}\\
{\scriptsize(f) $1500\times 1500\times 1500$,\\
{\bf r} = ($15\times 15\times 15$).}
\end{center}
\end{minipage}\\
\end{tabular}
\vspace{-0.1cm}
\caption{\changeHK{\bf Case S8:} medium-scale comparisons.}
\label{appnfig:middle-scale}
\end{figure}

{\bf Case S8: medium-scale instances.} We additionally consider medium-scale tensors of size $500 \times 500 \times 500$, $1000 \times 1000 \times 1000$, and $1500 \times 1500 \times 1500$ and ranks ${\bf r}=(5,5,5), (10,10,10)$, and $(15,15,15)$. OS is $\{10,20,30,40\}$. Our proposed algorithm and geomCG are only compared as the other algorithms cannot handle these scales efficiently. \changeHK{Figures \ref{appnfig:middle-scale}(a)-(c) show the convergence behavior. Figures \ref{appnfig:middle-scale}(d)-(f) also show the mean square error on $\Gamma$ of rank ${\bf r}=(15,15,15)$ in all the five runs.} The proposed algorithm performs better than geomCG in all the cases.

{\bf Experiments on real-world datasets:}

\begin{figure}[htbp]
\begin{tabular}{cc}
\begin{minipage}{0.48\hsize}
\begin{center}
\includegraphics[width=\hsize]{figures/caseR1_riebeira_small_OS_11_meansquaretesterror.eps}\\
{\scriptsize(a) OS = $11$.} 
\end{center}
\end{minipage}
\begin{minipage}{0.48\hsize}
\begin{center}
\includegraphics[width=\hsize]{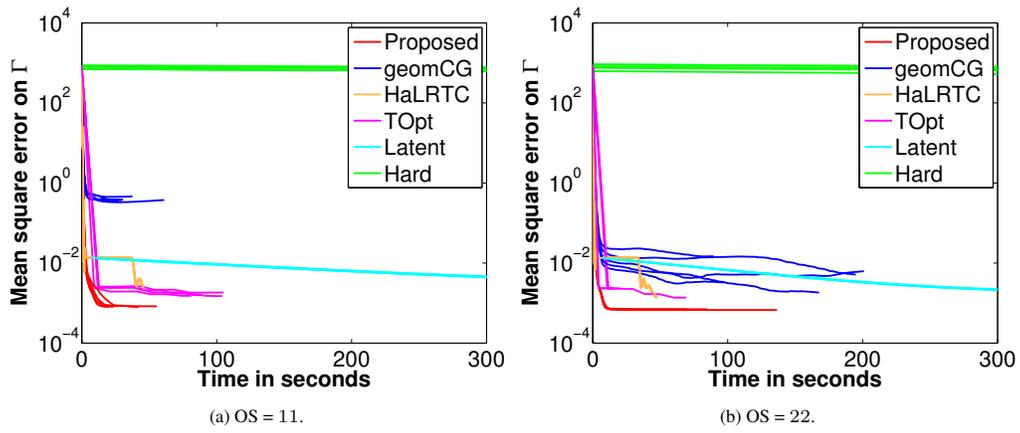}\\
{\scriptsize(b) OS = $22$.} 
\end{center}
\end{minipage}
\end{tabular}
\caption{\changeHK{{\bf Case R1:} mean square error on $\Gamma$.}}
\label{appnfig:R1}
\end{figure}

\begin{figure}[htbp]
\begin{tabular}{cccc}
\begin{minipage}{0.24\hsize}
\begin{center}
\includegraphics[width=\hsize, bb=0 0 268 203]{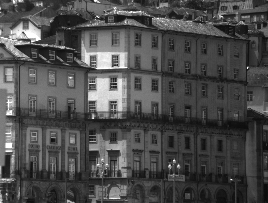}

{\scriptsize (a) Original.}
\label{fig:winter}
\end{center}
\end{minipage}
\begin{minipage}{0.24\hsize}
\begin{center}
\includegraphics[width=\hsize, bb=0 0 268 203]{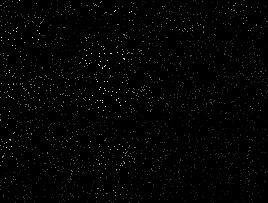}

{\scriptsize (b) Sampled ($4.98$\% observed).}
\label{fig:fall}
\end{center}
\end{minipage}

\begin{minipage}{0.24\hsize}
\begin{center}
\includegraphics[width=\hsize, bb=0 0 268 203]{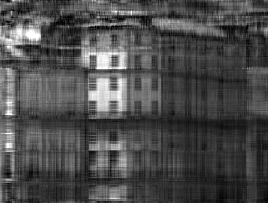}
{\scriptsize (c) Proposed.}
\label{fig:winter}
\end{center}
\end{minipage}
\begin{minipage}{0.24\hsize}
\begin{center}
\includegraphics[width=\hsize, bb=0 0 268 203]{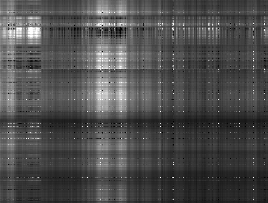}
{\scriptsize (d) geomCG.}
\label{fig:fall}
\end{center}
\end{minipage}\vspace*{0.2cm}\\

\begin{minipage}{0.24\hsize}
\begin{center}
\includegraphics[width=\hsize, bb=0 0 268 203]{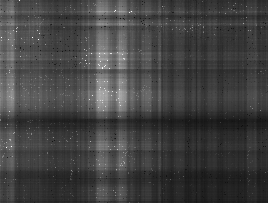}

{\scriptsize (e) HaLRTC.}
\label{fig:winter}
\end{center}
\end{minipage}
\begin{minipage}{0.24\hsize}
\begin{center}
\includegraphics[width=\hsize, bb=0 0 268 203]{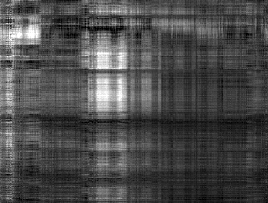}

{\scriptsize (f) TOpt.}
\label{fig:fall}
\end{center}
\end{minipage}

\begin{minipage}{0.24\hsize}
\begin{center}
\includegraphics[width=\hsize, bb=0 0 268 203]{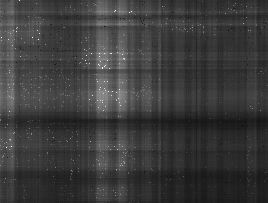}

{\scriptsize (g) Latent.}
\label{fig:winter}
\end{center}
\end{minipage}
\begin{minipage}{0.24\hsize}
\begin{center}
\includegraphics[width=\hsize, bb=0 0 268 203]{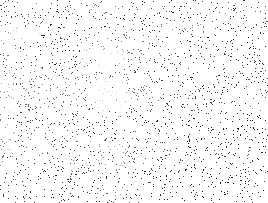}

{\scriptsize (h) Hard.}
\label{fig:fall}
\end{center}
\end{minipage}
\end{tabular}
\caption{\changeHK{{\bf Case R1:} recovery results on the hyperspectral image ``Ribeira" (frame = $16$, OS = $11$).}}
\label{appnfig:R1-reconstructedimage_OS_11}
\end{figure}

\begin{figure}[htbp]
\begin{tabular}{cccc}
\begin{minipage}{0.24\hsize}
\begin{center}
\includegraphics[width=\hsize, bb=0 0 268 203]{figures/ref_ribeira1bbb_reg1_resize_203x268x33.png}

{\scriptsize (a) Original.}
\label{fig:winter}
\end{center}
\end{minipage}
\begin{minipage}{0.24\hsize}
\begin{center}
\includegraphics[width=\hsize, bb=0 0 268 203]{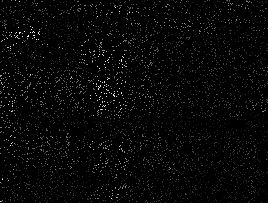}

{\scriptsize (b) Sampled ($9.96\%$ observed).}
\label{fig:fall}
\end{center}
\end{minipage}

\begin{minipage}{0.24\hsize}
\begin{center}
\includegraphics[width=\hsize, bb=0 0 268 203]{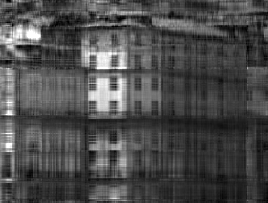}

{\scriptsize (c) Proposed.}
\label{fig:winter}
\end{center}
\end{minipage}
\begin{minipage}{0.24\hsize}
\begin{center}
\includegraphics[width=\hsize, bb=0 0 268 203]{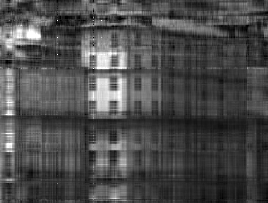}

{\scriptsize (d) geomCG.}
\label{fig:fall}
\end{center}
\end{minipage}\vspace*{0.2cm}\\

\begin{minipage}{0.24\hsize}
\begin{center}
\includegraphics[width=\hsize, bb=0 0 268 203]{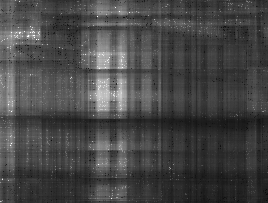}

{\scriptsize (e) HaLRTC.}
\label{fig:winter}
\end{center}
\end{minipage}
\begin{minipage}{0.24\hsize}
\begin{center}
\includegraphics[width=\hsize, bb=0 0 268 203]{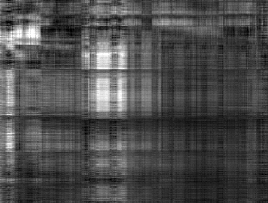}

{\scriptsize (f) TOpt.}
\label{fig:fall}
\end{center}
\end{minipage}

\begin{minipage}{0.24\hsize}
\begin{center}
\includegraphics[width=\hsize, bb=0 0 268 203]{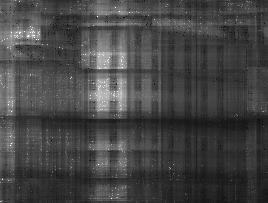}

{\scriptsize (g) Latent.}
\label{fig:winter}
\end{center}
\end{minipage}
\begin{minipage}{0.24\hsize}
\begin{center}
\includegraphics[width=\hsize, bb=0 0 268 203]{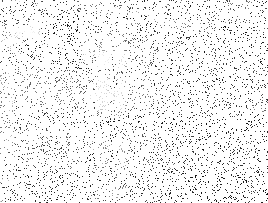}

{\scriptsize (h) Hard.}
\label{fig:fall}
\end{center}
\end{minipage}
\end{tabular}
\caption{\changeHK{{\bf Case R1:} recovery results on the hyperspectral image ``Ribeira" (frame = $16$, OS = $22$).}}
\label{appnfig:R1-reconstructedimage_OS_22}
\end{figure}

\changeHK{{\bf Case R1: hyperspectral image.} We also show the performance of our algorithm on the hyperspectral image ``Ribeira". We show the mean square error on $\Gamma$ when OS is $\{11, 22\}$ in Figure \ref{appnfig:R1}, where (a) is identical to the figure in the manuscript paper. Our proposed algorithm gives lower test errors than those obtained by the other algorithms. We also show the image recovery results. Figures \ref{appnfig:R1-reconstructedimage_OS_11} and \ref{appnfig:R1-reconstructedimage_OS_22} show the reconstructed images when OS is $\{11, 22\}$, respectively. From these figures, we find that the proposed algorithm shows a good performance, especially for the lower sampling ratio.}  

\begin{figure}[t]
\begin{tabular}{cc}
\begin{minipage}{0.48\hsize}
\begin{center}
\includegraphics[width=\hsize]{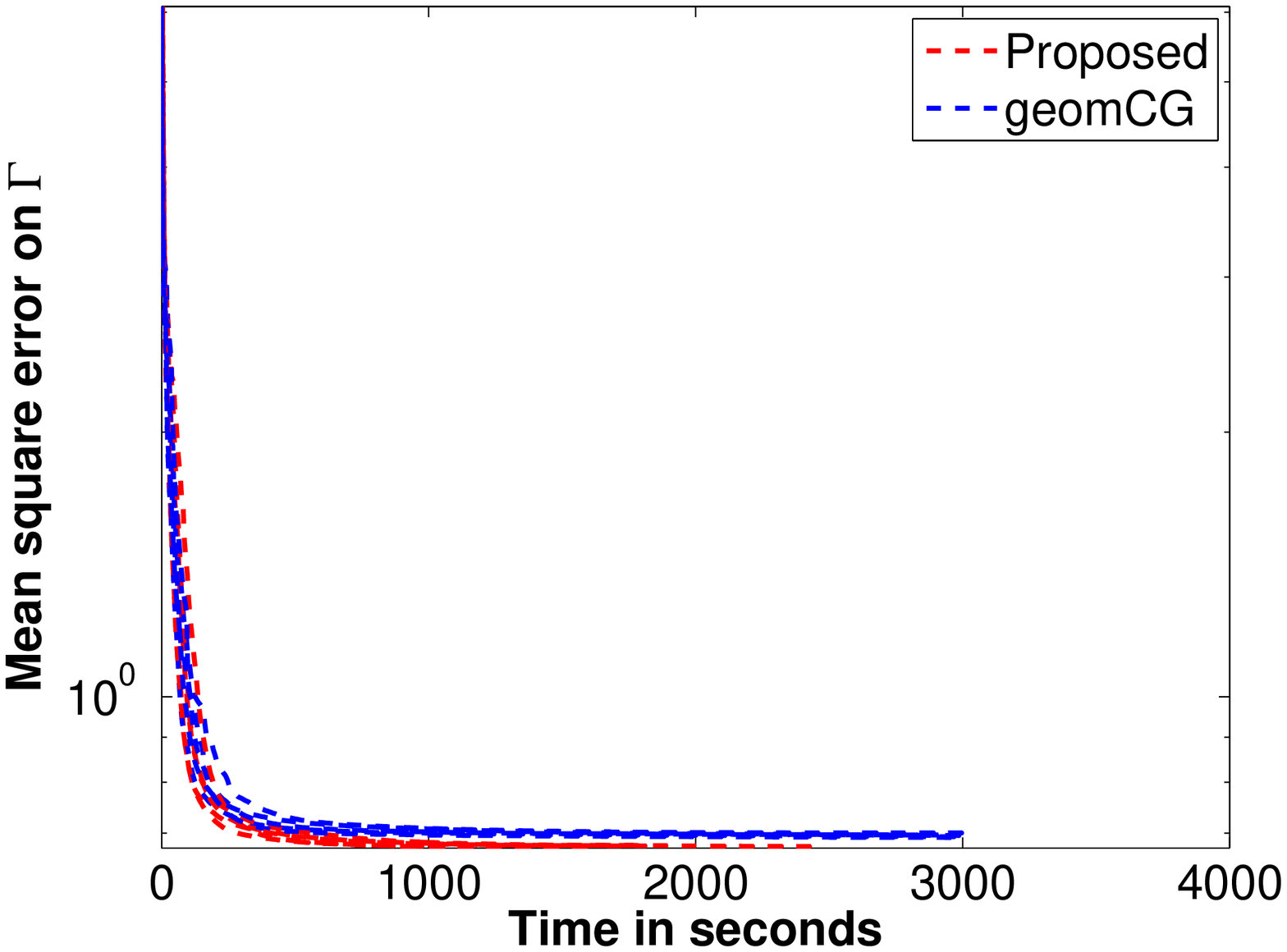}\\
{\scriptsize(a) {\bf r} = ($4\times 4\times 4$).} 
\end{center}
\end{minipage}
\begin{minipage}{0.48\hsize}
\begin{center}
\includegraphics[width=\hsize]{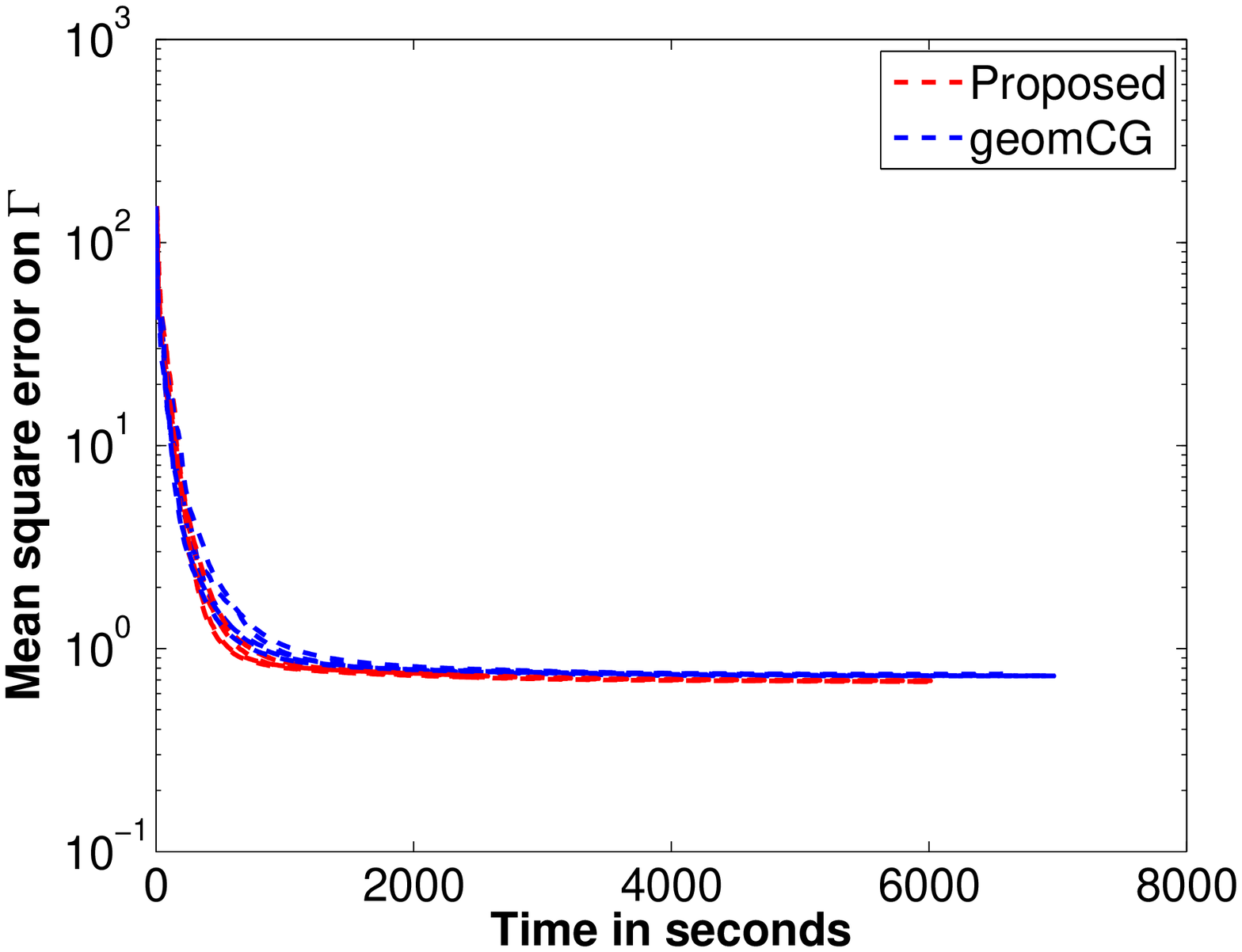}\\
{\scriptsize(b) {\bf r} = ($6\times 6\times 6$).} 
\end{center}
\end{minipage}\\

\begin{minipage}{0.48\hsize}
\begin{center}
\includegraphics[width=\hsize]{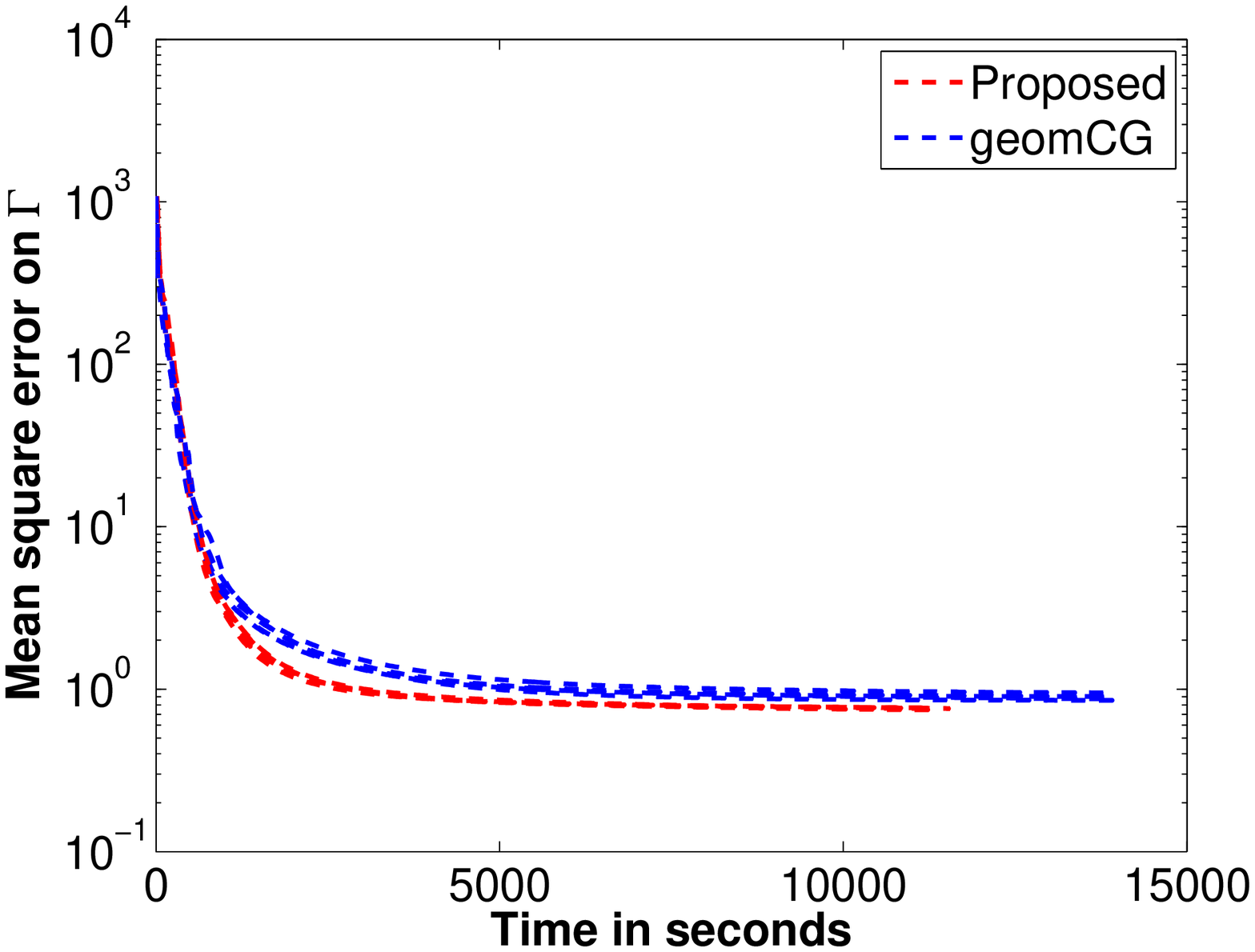}\\
{\scriptsize(c) {\bf r} = ($8\times 8\times 8$).} 
\end{center}
\end{minipage}
\begin{minipage}{0.48\hsize}
\begin{center}
\includegraphics[width=\hsize]{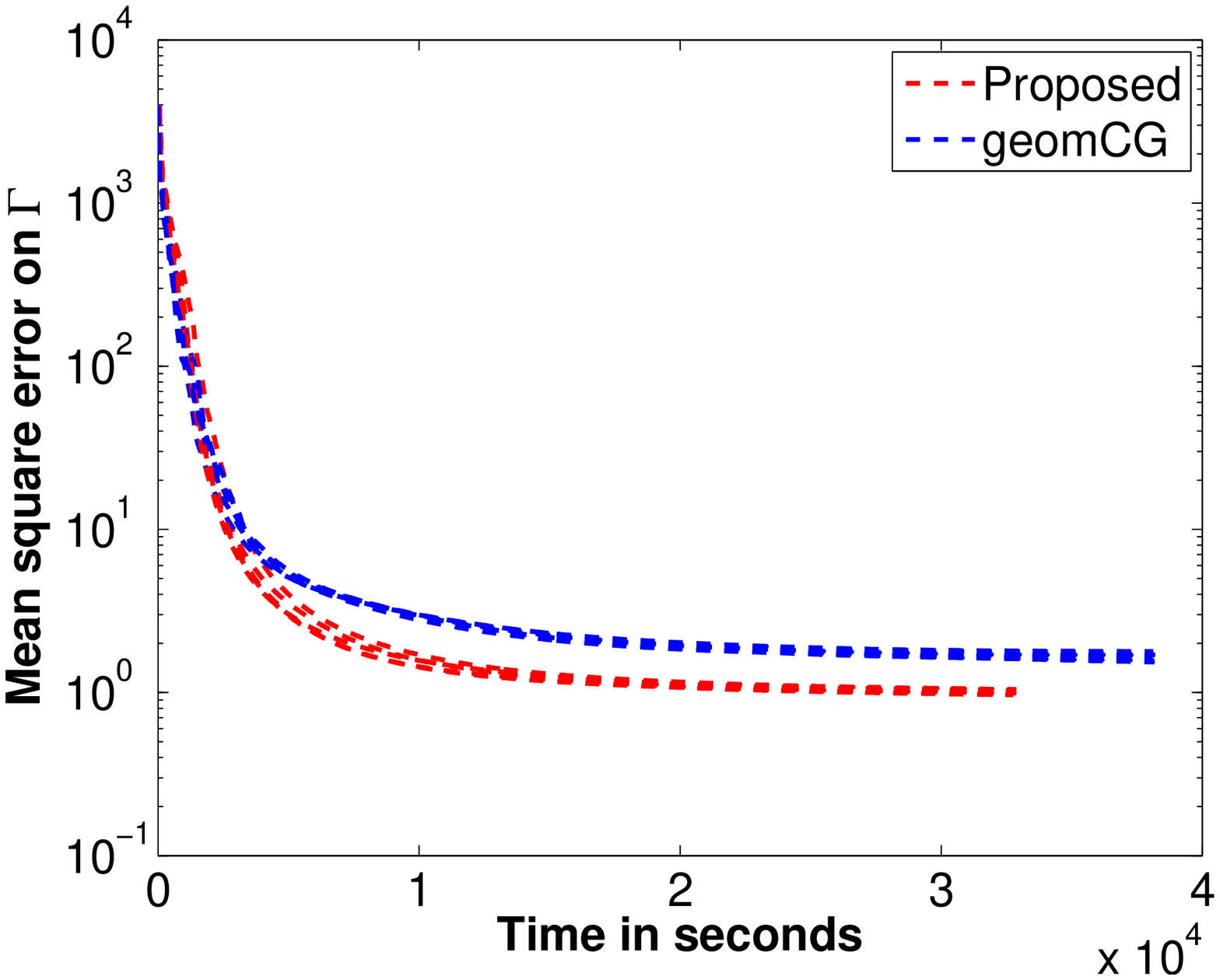}\\
{\scriptsize(d) {\bf r} = ($10\times 10\times 10$).} 
\end{center}
\end{minipage}
\end{tabular}
\caption{\changeHK{{\bf Case R2:} mean square error on $\Gamma$.}}
\label{appnfig:R2}
\end{figure}

\changeHK{
{\bf Case R2: MovieLens-10M.} Figure \ref{appnfig:R2} shows the convergence plots for all the five runs \changeHK{of ranks ${\bf r}=(4,4,4)$, $(6,6,6)$, $(8,8,8)$ and $(10,10,10)$}. 
These figures show the superior performance of our proposed algorithm.
}

\end{document}